\documentclass[a4paper,12pt]{article}
\usepackage[T2A]{fontenc}
\usepackage[cp1251]{inputenc}
\begin{document}

\author{S.V. Ludkovsky.}

\title{Twisted cohomologies of wrap groups over quaternions and octonions.}

\date{5 May 2008}
\maketitle

\begin{abstract}
This article is devoted to the investigation of wrap groups of
connected fiber bundles over the fields of real $\bf R$, complex
$\bf C$ numbers, the quaternion skew field $\bf H$ and the octonion
algebra $\bf O$. Cohomologies of wrap groups and their structure are
investigated. Sheaves of wrap groups are constructed and studied.
Moreover, twisted cohomologies and sheaves over quaternions and
octonions are investigated as well.
\end{abstract}

\section{Introduction.}
\par Geometric loop groups of circles were first introduced by
Lefshetz in 1930-th and then their construction was reconsidered by
Milnor in 1950-th. Lefshetz has used the $C^0$-uniformity on
families of continuous mappings, which led to the necessity of
combining his construction with the structure of a free group with
the help of words. Later on Milnor has used the Sobolev's
$H^1$-uniformity, that permitted to introduce group structure more
naturally \cite{milmorse}. \par  The construction of Lefshetz is
very restrictive, because it works with the $C^0$ uniformity of
continuous mappings in compact-open topology. Even for spheres $S^n$
of dimension $n>1$ it does not work directly, but uses the iterated
loop group construction of circles. Then their constructions were
generalized for fibers over circles and spheres with parallel
transport structures over $\bf C$. Smooth Deligne cohomologies were
studied on such groups \cite{gajer}.

\par Wrap groups of quaternion and octonion fibers
as well as for wider classes of fibers over $\bf R$ or $\bf C$ were
defined and various examples were given together with basic theorems
in \cite{luwrgfbqo}. Studies of their structure were begun in
\cite{lulaswgof}. This paper continues previous works of the author
on this theme, where generalized loop groups of manifolds over $\bf
R$, $\bf C$ and $\bf H$ were investigated, but neither for fibers
nor over octonions \cite{ludan,lugmlg,lujmslg,lufoclg}.

\par Holomorphic functions of quaternion and octonion variables were
investigated in \cite{luoyst,luoyst2,lufsqv}. There specific
definition of super-differentiability was considered, because the
quaternion skew field has the graded algebra structure. This
definition of super-differentiability does not impose the condition
of right or left super-linearity of a super-differential, since it
leads to narrow class of functions. There are some articles on
quaternion manifolds, but practically they undermine a complex
manifold with additional quaternion structure of its tangent space
(see, for example, \cite{museya,yano} and references therein).
Quaternion manifolds were defined in another way in \cite{lufsqv}.
Applications of quaternions in mathematics and physics can be found
in \cite{emch,guetze,hamilt,lawmich}.
\par  Geometric loop groups have important applications in modern physical
theories (see \cite{ish,mensk} and references therein). Groups of
loops are also intensively used in gauge theory. Wrap groups can be
used in the membrane theory which is the generalization of the
string (superstring) theory. Fiber bundles and sheaves and
cohomologies over quaternions and octonions are interesting in such
a respect, that they take into account spin and isospin structures
on manifolds, because there is the embedding of the Lie group $U(2)$
into the quaternion skew field $\bf H$.
\par This article is devoted to constructions and investigations
of cohomologies and sheaves of wrap groups. Moreover, over
quaternions and octonions twisted cohomologies and sheaves are
studied. Twisted analogs of bar resolutions of sheaves and smooth
Deligne cohomology are investigated as well. This is done over
twisted multiplicative groups. Previously the complex case and with
loop groups of fiber bundles on spheres was only studied. This
article treats the quaternion and octonion cases and wrap groups of
general fiber bundles.

\par All main results of this paper are obtained for the first time
and they are given in Theorems 34, 36, 44, 48.1, 55, 58, 60,
Propositions 6, 14, 15, 19, 26, 27, 29, 32, Corollaries 7, 8, 33, 45
and 47. Here the notations and definitions and results from previous
works
\cite{gajer,gajinvm,luwrgfbqo,lulaswgof,ludan,lugmlg,lujmslg,lufoclg}
are used.

\section{Cohomologies of wrap groups}
\par {\bf 1. Remarks and Definitions.} Consider a
triangulated compact polyhedron $M$ may be embedded into ${\cal
A}_r^n$ and its sub-polyhedron $S_M$ of codimension not less than
two, $codim (S_M)\ge 2$, where $M\setminus S_M$ is a $C^{\infty }$
smooth manifold such that $M\setminus S_M$ is dense in $M$. If the
covering dimension (see Chapter 7 \cite{eng}) of $M\setminus S_M$ is
$dim (M\setminus S_M)=b$, then by the definition $M$ is of dimension
$b$. Then $S_M$ is called the singularity of $M$. A pseudo-manifold
$M$ is oriented, if $M\setminus S_M$ is oriented (see also \S 1.3.1
\cite{luwrgfbqo}).
\par If $M\setminus S_M$ is without boundary, then the triangulated
pseudo-manifold $M$ is called a pseudo-manifold cycle. If
$(Y,\partial Y)$ is the pair consisting of a triangulated
pseudo-manifold $Y$ and a boundary $\partial Y$, such that
$Y\setminus S_Y$ is a manifold with boundary $\partial Y\setminus
S_Y$, $\partial Y$ is a pseudo-manifold cycle with singularity
$S_Y\cap \partial Y$, then $(Y,\partial Y)$ is called the
triangulated pseudo-manifold with boundary.

\par A pre-sheaf $F$ on a topological space $X$ is a
contra-variant functor $F$ from the category of open subsets in $X$
and their inclusions into a category of groups or rings (all either
alternative or associative) such that $F(U)$ is a group or a ring
for each $U$ open in $X$ and for each $U\subset V$ with $U$ and $V$
open in $X$ there exists a homomorphism $s_{U,V}: F(V)\to F(U)$ such
that $s_{U,U}=1$ and $s_{U,V}s_{V,Y}=s_{U,Y}$ for each $U\subset
V\subset Y$ with $U, V, Y$ open in $X$.

\par Let $F_x$ denotes the family of all elements $f\in
F(U)$ for all $U$ open in $X$ with $x\in U$. Elements $f\in F(U)$
and $g\in F(V)$ are called equivalent if there exists an open
neighborhood $Y$ of $x$ such that $s_{Y,U}(f)=s_{Y,V}(g)$. This
generates an equivalence relation and a class of all equivalent
elements with $f$ is called a germ $f_x$ of $f$ at $x$. A set ${\cal
F}_x$ of all germs of the pre-sheaf $F$ at a point $x\in X$ is the
inductive limit ${\cal F}_x = ind-\lim F(U)$ taken by all open
neighborhoods $U$ of $x$ in $X$. \par In the set $\cal F$ of all
germs ${\cal F}_x$ take a base of topology consisting of all sets $
\{ f_x\in {\cal F}_x: x\in U \} $, where $f\in F(U)$. This induces a
sheaf $\cal S$ generated by a pre-sheaf $F$.

\par A sheaf of groups or rings (all either alternative or associative)
on $X$ is a pair $({\cal S}, h)$ satisfying Conditions $(S1-S4)$:
\par $(S1)$ $\cal S$ is a topological space;
\par $(S2)$ $h: {\cal S}\to X$ is a local homeomorphism;
\par $(S3)$ for each $x\in X$ the set ${\cal F}_x = h^{-1}(x)$ is
a group called a fiber of the sheaf $\cal S$ at a point $x$;
\par $(S4)$ the group or the ring operations are continuous, that is,
${\cal S}\Delta {\cal S}\ni (a,b)\mapsto ab^{-1}\in \cal S$ or
${\cal S}\Delta {\cal S}\ni (a,b)\mapsto ab\in \cal S$ and ${\cal
S}\Delta {\cal S}\ni (a,b)\mapsto a+b\in \cal S$ are continuous
respectively, where ${\cal S}\Delta {\cal S} := \{ (a,b): a, b\in
{\cal S}, h(a)=h(b) \} $.
\par We can consider pre-sheafs and sheafs
of different classes of smoothness, for example, $H^t$ or $H^t_p$,
when the corresponding defining sheaf and pre-sheaf mappings
$s_{U,V}$, $h$ and group operations are such and $\cal S$ and $F$
are $H^t$ or $H^t_p$ differentiable spaces respectively (see also \S
1.3.2 \cite{luwrgfbqo}).

\par Consider a sheaf ${\cal S}_{N,G}$ generated by a pre-sheaf
$U\mapsto \{ f\in Hom^t_p((W^ME)_{t,H},G): supp (f)\subset U \} $,
where $U$ is open in $N$ and $supp (f)\subset U$ means that there
exists $y\in N$ and ${\hat {\eta }}\in H^t_p({\hat M},N)$ with
${\hat {\eta }}({\hat s}_{0,q})=x$ for each $q=1,...,k$ and ${\hat
{\eta }}({\hat s}_{0,q})=y$ for each $q=k+1,...,2k$ and ${\hat
{\gamma }}\in H^t_p({\hat M}, \{ {\hat s}_{0,q}: q=1,...,2k \} ;
N,y)$ such that ${\hat {\gamma }} = \gamma \circ \Xi $ and $f =
<{\hat {\eta }}\vee {\hat {\gamma }}>_{t,H}$, where the wrap group
$(W^ME)_{t,H}$ is taken for a marked point $y\in N$, $\Xi : {\hat
M}\to M$ is the quotient mapping as in \cite{luwrgfbqo}. \par In
particular, we can take $G={\cal A}_r^*$, and call ${\cal
S}_{N,{\cal A}_r^*}$ the sheaf of infinitesimal holonomies, where
$1\le r\le 3$.
\par In view of Property $(P4)$ \cite{luwrgfbqo} for each non-singular
points $y\in N$ and $u\in E_y$ in the fiber $E_y$ of $E$ over $y$
there exists an ${\cal A}_r$ vector subspace $H_u$ of the tangent
bundle $T_uE$ at $u$ called a horizontal subspace of $T_uE$ such
that $\pi _*|_{H_u}: H_u\to T_yN$ is an isomorphism, where $\pi
(u)=y$, $t'\ge [dim (E)/2]+2$ or $t'=\infty $, since there exist
generalized derivatives in the Sobolev space $H^{t'}$ (see \S III.3
\cite{miha}). This is the case for all $y\in N$ and $u\in E_y$ when
$N$ and $E$ are of class $H^{t'}$ instead of $H^{t'}_p$.
\par Due to $(P1)$ the family $ \{ H_u \} $ of horizontal subspaces
of $TE$ depends smoothly on $u$. Suppose that $Y$ is a vector field
in $TE$ corresponding to a vector field $X$ in $TN$ such that $\pi
_*(Y)=X$, then
\par $(CD1)$ $T_uE=H_u\oplus V_u$, where $V_u=\pi_*^{-1}(0)\subset
T_uE$ is the space of vectors tangent to $E_u$ at $u$. In accordance
with $(P3)$ the horizontal spaces are $G$-equivariant, that is,
\par $(CD2)$ $H_{uz} = (R_z)_*H_u$, where $R_z$ is the
diffeomorphism of $E$ given by the multiplication on $z$ from the
right and $(R_z)_*$ corresponds to the tangent mapping $TR_z$ for
the tangent fiber bundle $TE$.
\par A family $H = \{ H_u\subset T_uE: u\in E, \pi (u)=y\in N \} $
is called the connection distribution of the principal fiber bundle
$E(N,G,\pi ,\Psi )$, if $H_u$ depends smoothly on $u$ and the
Conditions $(CD1,CD2)$ are satisfied.

\par {\bf 2. Definitions and Notes.} Two smooth principal $G$ fiber bundles
$E$ and $E'$ with connection distributions $(E,H)$ and $(E',H')$ are
called isomorphic if there exists an isomorphism $f: E\to E'$ of
smooth principal $G$ fiber bundles $f: E\to E'$ such that
$f_*(H)=H'$.
\par A connection distribution $H$ on $E$ determines a parallel
transport structure ${\bf P}^H$ on $E$ posing ${\bf P}^H_{{\hat
{\gamma }},u}\in H^t_p({\hat M},E)$ with ${\bf P}^H_{{\hat {\gamma
}},u}({\hat s}_{0,q})=u$ for each $q=1,...,k$ and $\pi \circ {\bf
P}^H_{{\hat {\gamma }},u}=\hat {\gamma }$ such that $T_x{\bf
P}^H_{{\hat {\gamma }},u}=: ({\bf P}^H_{{\hat {\gamma }},u}(x),
D{\bf P}^H_{{\hat {\gamma }},u}(x))$ for each $x\in \hat M$, where
${\hat {\gamma }}\in H^t_p({\hat M}, \{ {\hat s}_{0,q}: q=1,...,2k
\}; N,y_0)$, $D{\bf P}^H_{{\hat {\gamma }},u}(x)\in H_{{\hat {\gamma
}}(x)}$, $T{\bf P}$ is the tangent mapping of $\bf P$ (see
\cite{kling}).

\par Thus there exists a bijective correspondence between
parallel transport structures and connection distributions on $E$.
Therefore, the mapping $H\mapsto {\bf P}^H$ induces a bijective
correspondence between isomorphism classes of parallel transport
structures and connection distributions.
\par Using the exponential function on ${\cal A}_r$ gives
$\exp ({\cal A}_r)= {\cal A}_r^*$ for $1\le r\le 3$ (see \S 3
\cite{luoyst,luoyst2}).

\par If $E$ is a principal
${\cal A}_r^*$ fiber bundle with $1\le r\le 3$, then for each $v\in
V_u$ there exists a unique $z(v)\in {\cal A}_r$ such that $v= [d(y~
\exp (b~ z(v))/db]|_{b=0}$, where $b\in \bf R$. Therefore, for each
connection distribution $ \{ H_u: u \in E \} $ on $E$ a differential
$1$-form $w$ over ${\cal A}_r$ exists such that $w(X_h + X_v)=
z(X_v)$ for each $X=X_h + X_v\in H_u\oplus V_u=T_uE$ and $w$ is
$G$-equivariant: $(R_z)^*w=w$ due to the $G$-equivariance of $ \{
H_u: u\in E \} $, here $G={\cal A}_r^*$.
\par A differential $1$-form $w$ on $E$ so that it is $G$-equvariant
and $w(X_v)=z(X_v)$ for each $X_v\in V_u$ is called a connection
$1$-form.

\par Two smooth principal $G$ fiber bundles with connections
$(E,w)$ and $(E',w')$ are called isomorphic, if there exists an
isomorphism $f: E\to E'$ of smooth principal $G$ fiber bundles such
that $f^*(w')=w$.

\par For $w$ there exists a connection distribution $H^w$ on $E$
for which $H^w_u = ker (w_u)\subset T_uE$, that induces a bijective
correspondence between differential $1$-forms and connection
distributions on $E$. Hence $w\mapsto H^w$ produces a bijective
correspondence between isomorphism classes of connections and
connection distributions.

\par Then there exists a wrap group $(W^ME;N,{\cal A}_r^*,\nabla
)_{t,H}$, where a parallel transport structure ${\bf P}$ is
associated with the covariant differentiation $\nabla $ of the
connection $w$.
\par The curvature $2$-form $\Omega $ over ${\cal A}_r$, $1\le r\le
3$, of a connection $1$-form $w$ on a smooth principal fiber bundle
$E(N,G,\pi ,\Psi )$ over ${\cal A}_r$ is given by $\Omega (X,Y) =
dw(hX,hY)$, where $hX$ and $hY$ are horizontal components of the
vectors $X$ and $Y$.

\par {\bf 3. Remark.} If $\eta \in H^t_p(K,E)$, $t\ge 1$,
and $\nu $ is a differential form on $E$, then there exists its
pull-back $\eta ^*\nu $ which is a differential form on $K$, where
$K$ is an $H^t_p$-pseudo-manifold. For orientable $K$ and $E$ and an
$H^t_p$ diffeomorphism $\eta $ of $K$ onto $E$ and $\nu $ with
compact support $\int_K \eta ^* \nu = \epsilon \int_E \nu $, where
$\epsilon =1$ if $\eta $ preserves an orientation, $\epsilon = -1$
if $\eta $ changes an orientation (see
\cite{dubnovfom,zorich,lufsqv}). In particular, $K = E(M,G,\pi _M,
\Psi ^M)$ can be considered, $\eta = (\eta _0, \eta _1)$, $\eta _0:
M\to N$, $\eta : E(M)\to E(N)$, $\pi _N\circ \eta = \eta _0\circ \pi
_M$, $\eta _1\circ pr_2 = pr_2\circ \eta $, $pr_2$ is a projection
in charts of $E$ from $E$ into $G$, $\eta _1 = id$ may be as well.

\par Suppose that $M$ and $E$ are an ${\cal A}_r$ holomorphic manifold
and principal fiber bundle, such that $E$ is orientable and
$2^r-1$-connected, which is not very restrictive due to Propositions
13 and Note 14 \cite{lulaswgof}. If ${\hat \gamma }\in H^t_p({\hat
M}, \{ {\hat s}_{0,q}: q=1,...,k \}; N,y_0)$, then consider a path
$l_q$ joining the point ${\hat s}_{0,q}$ with ${\hat s}_{0,q+k}$,
where $1\le q\le k$, ${\hat l}_q: [0,1]\to {\hat M}$. Therefore,
${\hat p}_q := {\hat \gamma }\circ {\hat l}_q: [0,1]\to N$ and
${\hat p}_q ^*w := ({\hat p}_q, id)^*w$ is a differential form on
$[0,1]$, where $w$ is an ${\cal A}_r$ holomorphic connection
one-form on $E$. We get that ${\hat \gamma }^*w$ is a differential
one-form on $\hat M$ and there exists its restriction $\nu _{\gamma
,q} := {\gamma ^*w}|_{{\hat l}_q[0,1]}$.
\par Then we have also $\gamma \in H^t_p(M, \{ s_{0,q}: q=1,...,k
\}; N,y_0)$ and $l_q: S^1\to M$ and $p_q: S^1\to N$ respectively,
where ${\hat \gamma } =\gamma \circ \Xi $ (see \cite{luwrgfbqo}),
$S^1$ is the unit circle in $\bf C$ with the center at zero, while
${\bf C}={\bf C}_M$ is embedded into ${\cal A}_r$ as ${\bf R}\oplus
M{\bf R}$ with $M\in {\cal A}_r$, $Re (M)=0$, $|M|=1$, when $2\le
r\le 3$.

\par Since $w$ is ${\cal A}_r$ holomorphic, then $\int_{\phi } w$ does
not depend on a rectifiable curve $\phi $ but only on the initial
and final points $\phi (0)$ and $\phi (1)$, $\phi : [0,1]\to E$ (see
Theorems 2.15 and 3.10 in \cite{luoyst,luoyst2} and \cite{lufoclg}).
\par Consider now the principal fiber bundle $E$ with the structure
group ${\cal A}_r^*$, where $1\le r\le 3$. Then the pull-back
$p_q^*E$ of the bundle $E$ is a trivial ${\cal A}_r^*$-bundle over
$S^1$. The latter bundle carries a pull-back connection differential
one-form $p_q^*w$. Take the pull-back $\rho ^*(p_q^*w)$ one form,
where $\rho : S^1\to p_q^*E$ is a trivialization of the fiber bundle
$p_q^*E\to S_1$.
\par The parallel transport structure ${\bf P}_{{\hat {\gamma
}},u}(x)$ for $(M,E)$ with $x\in \hat M$ induces the parallel
transport structures ${\bf P}_{{\hat p}_q,u^*}(s)$ for
$(S^1,p_q^*E)$ with $s\in [0,1]$ for each $q=1,...,k$, where
$p_q(u^*)=u$. Then the holonomy along $\gamma $ is given by
\par $(H)$ $h(\gamma ) = (h_1,...,h_k)\in G^k$ with $h_q =
h_q(\gamma )= \exp [-\int_{S^1} \rho ^*(p_q^*w)]$ for each
$q=1,...,k$.

\par If $\zeta : S^1\to p_q^*E$ is another trivialization and $f: S^1\to
{\bf C}_M^*$ satisfies $\zeta =f\rho $, so that $f(v) = \exp (M2\pi
\theta (v))$, where $\theta (v)\in \bf R$, $M2\pi d\theta (v) = d Ln
(f(v))$, $v\in S^1$, $\int _{S^1}d\theta $ is an integer number,
since $\bf R$ is the center of the algebra ${\cal A}_r$, where $Ln$
is the natural logarithmic function over ${\cal A}_r$ (see \S 3.7
and Theorem 3.8.3 \cite{luoyst2} and \cite{luoyst,lufejms}).
Therefore, Formula $(H)$ is independent of a trivialization $\rho $,
since $\zeta ^*(p_q^*w) = \rho ^*(p_q^*w) + d Ln (f)$, but $\exp
[\int_{S^1} d Ln (f)]=1$.

\par {\bf 4. Non-associative bar construction.}
Let $G$ be a topological group not necessarily associative, but
alternative:
\par $(A1)$ $g(gf)=(gg)f$ and $(fg)g=f(gg)$ and $g^{-1}(gf)=f$
and $(fg)g^{-1}=f$ for each $f, g\in G$ \\
and having a conjugation operation which is a continuous
automorphism of $G$ such that \par $(C1)$ $conj (gf)=conj (f)
conj(g)$ for each $g, f\in G$, \par $(C2)$ $conj (e)=e$ for the unit
element $e$ in $G$.
\par If $G$ is of definite class of smoothness, for example, $H^t_p$
differentiable, then $conj$ is supposed to be of the same class. For
commutative group in particular it can be taken the identity mapping
as the conjugation. For $G= {\cal A}_r^*$ it can be taken $conj
(z)={\tilde z}$ the usual conjugation for each $z\in {\cal A}_r^*$,
where $1\le r\le 3$.
\par Denote by $\Delta ^n :=
\{ (x_0,...,x_n)\in {\bf R}^{n+1}: x_j\ge 0, x_0+x_1+...+x_n=1 \} $
the standard simplex in ${\bf R}^{n+1}$. Consider $(AG)_n$ as the
quotient of the disjoint union $\bigcup_{k=0}^n (\Delta ^k\times
G^{k+1})$ by the equivalence relations \par $(1)$
$(x_0,...,x_k,g_0,...,g_k) \sim
(x_0,...,x_j+x_{j+1},...,x_k,g_0,...,{\hat g}_j,...,g_k)$ for
$g_j=g_{j+1}$ or $x_j=0$ with $0\le j<k$; $(x_0,...,x_k,g_0,...,g_k)
\sim (x_0,...,x_{k-1}+x_k,g_0,...,g_{k-1})$ for $g_{k-1}=g_k$ or
$x_k=0$.

\par Consider non-homogeneous coordinates $0\le t_1\le t_2\le ...
\le t_k\le 1$ on the simplex $\Delta ^k$ related with the
barycentric coordinates by the formula $t_j=x_0+x_1+...+x_{j-1}$ and
$h_0 := g_0$, $h_j= g_{j-1}^{-1}g_j$ for $j>0$ on $G^{k+1}$. Hence
$h_0h_1 = g_0(g_0^{-1}g_1) =g_1$, $(h_0h_1)h_2 = g_1(g_1^{-1}g_2) =
g_2$ and by induction $((...(h_0h_1)...)h_{k-1})h_k
=g_{k-1}(g_{k-1}^{-1}g_k) =g_k$.

\par Then equivalence relations
$(1)$ take the form:
\par $(2)$ $(t_1,...,t_k,h_0[h_1|...|h_k])\sim
(t_2,...,t_k, (h_0 h_1) [h_2|...|h_k])$ for $t_1=0$ or $h_0=e$;
\par $(t_1,...,t_k,h_0[h_1|...|h_k])\sim
(t_1,...,{\hat t}_j,...,t_k,h_0[h_1|...|h_j h_{j+1}|...|h_k])$ for
$t_j=t_{j+1}$ or $h_j=e$;
\par $(t_1,...,t_k,h_0[h_1|...|h_k])\sim
(t_1,...,t_{k-1},h_0[h_1|...|h_{k-1}])$ for $t_k=1$ or $h_k=e$.
\par Denote by $|x_0,...,x_k,g_0,...,g_k|$ the equivalence class
of the sequence \\ $(x_0,...,x_k,g_0,...,g_k)$; by
$|t_1,...,t_k,h_0[h_1|...|h_k]|$ denote the equivalence class of the
sequence $(t_1,...,t_k,h_0[h_1|...|h_k])$.

\par Then the space $AG$ is the quotient of $\bigcup_{k=0}^{\infty }
\Delta ^k\times G^{k+1}$ by the above equivalence relations $(1)$,
where $(\Delta ^k\times G^{k+1})\cap (\Delta ^m\times G^{m+1})$ is
empty for $k\ne m$.

\par Introduce in $G^{n+1}$ the equivalence relation $\cal Y$:
\par $(3)$ $(g_0,...,g_n){\cal Y} (q_0,...,q_n)$ if and only if
there exist $p_1,...,p_k\in G$ with $k\in \bf N$ such that
$g_j=p_k(p_{k-1}...(p_2(p_1q_j))....)$ for each $j=0,...,n$.
\par Evidently this relation is reflexive: $(g_0,...,g_n){\cal Y}
(g_0,...,g_n)$ with $p_1=e$ and $k=1$. It is symmetric due to the
alternativity of $G$, since $g_j=p_k(p_{k-1}...(p_2(p_1q_j))....)$
is equivalent with
$q_j=p_1^{-1}(p_2^{-1}...(p_{k-1}^{-1}(p_k^{-1}g_j))...)$ for each
$j=0,...,n$. This relation is transitive: $(g_0,...,g_n){\cal Y}
(q_0,...,q_n)$ and $(q_0,...,q_n){\cal Y} (f_0,...,f_n)$ implies
$(g_0,...,g_n){\cal Y} (f_0,...,f_n)$, since from
\par $g_j=p_k(p_{k-1}...(p_2(p_1q_j))....)$ and
$q_j=s_l(s_{l-1}...(s_2(s_1f_j))....)$ \\ it follows
$g_j=p_k(p_{k-1}...(p_2(p_1(s_l(s_{l-1}...(s_2(s_1f_j))....)))....)$
for each $j=0,...,n$, where $k, l\in \bf N$, $p_1,...,p_k,
s_1,...,s_l\in G$. In a particular case of an associative group $G$
parameters $k=1$ and $l=1$ can be taken.
\par Consider in $\bigcup_{k=0}^n \Delta ^k\times G^k$ the
equivalence relations: \par $(4)$ $(x_0,...,x_k,[g_0:...:g_k]) \sim
(x_0,...,x_j+x_{j+1},...,x_k,[g_0:...:{\hat g}_j:...:g_k])$ for
$g_j=g_{j+1}$ or $x_j=0$ with $0\le j<k$; $(x_0,...,x_k,g_0,...,g_k)
\sim (x_0,...,x_{k-1}+x_k,[g_0:...:g_{k-1}])$ for $g_{k-1}=g_k$ or
$x_k=0$, where $[g_0:...:g_k] := \{ (q_0,...,q_k)\in G^{k+1}:
(q_0,...,q_k){\cal Y}(g_0,...,g_k) \} $ denotes the equivalence
class of $(g_0,...,g_k)$ by the equivalence relation $\cal Y$. Put
$(BG)_n$ to be the quotient of $\bigcup_{k=0}^n \Delta ^k\times G^k$
by equivalence relations $(4)$.

\par Using the inhomogeneous coordinates on $(BG)_n$ rewrite
the equivalence relation $(4)$ in the form:
\par $(5)$ $(t_1,...,t_k,[h_1|...|h_k])\sim
(t_2,...,t_k, [h_2|...|h_k])$ for $t_1=0$ or $h_0=e$;
\par $(t_1,...,t_k,h_0[h_1|...|h_k])\sim
(t_1,...,{\hat t}_j,...,t_k,[h_1|...|h_j h_{j+1}|...|h_k])$ for
$t_j=t_{j+1}$ or $h_j=e$;
\par $(t_1,...,t_k,[h_1|...|h_k])\sim
(t_1,...,t_{k-1},[h_1|...|h_{k-1}])$ for $t_k=1$ or $h_k=e$.

\par Denote by $|x_0,...,x_k,[g_0:...:g_k]|$ the equivalence class
of the sequence $(x_0,...,x_k,[g_0:...:g_k])$; by
$|t_1,...,t_k,[h_1|...|h_k]|$ denote the equivalence class of the
sequence $(t_1,...,t_k,[h_1|...|h_k])$. Then $BG$ is the quotient of
the disjoint union $\bigcup_{k=0}^{\infty } \Delta ^k\times G^k$ by
the equivalence relations $(4)$.
\par Then there exists the projection $\pi ^A_B: AG\to BG$ by the formula:
\par $(6)$ $\pi ^A_B: |x_0,...,x_k,g_0,...,g_k|\mapsto
|x_0,...,x_k,[g_0:...:g_k]|$ or in the non-homogeneous coordinates
by $\pi ^A_B: |t_1,...,t_k, h_0[h_1|...|h_k]|\mapsto
|t_1,...,t_k,[h_1|...|h_k]|$.

\par The conjugation in $G$ induces that of in $AG$ and $BG$
such that:
\par $conj (t_1,...,t_k, h_0[h_1|...|h_k]) :=
(t_1,...,t_k, conj(h_0) [conj (h_1)|...|conj (h_k)])$ and
\par $conj (t_1,...,t_k, [h_1|...|h_k]) :=
(t_1,...,t_k, [conj (h_1)|...|conj (h_k)])$.

\par Suppose that \par $(A2)$ ${\hat G} =
{\hat G}_0i_0\oplus {\hat G}_1i_1\oplus ... \oplus {\hat
G}_{2^r-1}i_{2^r-1}$ such that $G$ is a multiplicative group of a
ring $\hat G$ with the multiplicative group structure, where
$G_j={\hat G}_j\setminus \{ 0 \} $, ${\hat G}_0,...,{\hat
G}_{2^r-1}$ are pairwise isomorphic commutative associative rings
and $ \{ i_0,...,i_{2^r-1} \} $ are generators of the Cayley-Dickson
algebra ${\cal A}_r$, $1\le r\le 3$ and
$(y_li_l)(y_si_s)=(y_ly_s)(i_li_s)$ is the natural multiplication of
any pure states in $\hat G$ for $y_l\in {\hat G}_l$. If a group $G$
and a ring $\hat G$ satisfy Conditions $(A1,A2,C1,C2)$, then we call
it a twisted group and a twisted ring over the set of generators $
\{ i_0,...,i_{2^r-1} \} $, where $1\le r\le 3$. The unit element of
$G$ denote by $e$. For example, $G=({\cal A}_r^*)^n$ and ${\hat G} =
{\cal A}_r^n$.

\par {\bf 5. Definitions.} Let $N_.$ be a family $\{ N_n:
n\in {\bf N} \} $ of either $C^{\infty }$ smooth or $H^{t'}_p$
manifolds together with either $C^{\infty }$ or $H^{t'}_p$ mappings
$\partial _j: N_n\to N_{n-1}$ and $s_j: N_n\to N_{n+1}$ for each
$j=0,1,...,n$ satisfying the identities: \par $(1)$ $\partial
_k\partial _j = \partial _{j-1}\partial _k$ for each $k<j$,
\par $(2)$ $s_ks_j=s_{j+1}s_k$ for each $k\le j$,
\par $(3)$ $\partial _k s_j = s_{j-1}\partial _k$ for
$k<j$, $\partial _k s_j = id|_{N_n}$ for $k=j, j+1$, $\partial _k
s_j = s_j\partial _{k-1}$ for $k>j+1$, then $N_.$ is called a
simplicial either $C^{\infty }$ smooth or $H^{t'}_p$ manifold.
\par The geometric realization $|N_.|$ of $N_.$
consists of $\coprod_{n\ge 0} \Delta ^n\times N_n/ {\cal E}$, where
$\cal E$ is the equivalence relation generated by $(\partial
^jx,y){\cal E} (x,\partial _jy)$ for $(x,y)\in \Delta ^{n-1}\times
N_n$, $(s^jx,y){\cal E} (x,s_jy)$ for $(x,y)\in \Delta ^{n+1}\times
N_n$, where $\coprod $ denotes the disjoint union of sets, the maps
$\partial ^j: \Delta ^{n-1}\to \Delta ^n$ and $s^j: \Delta ^{n+1}\to
\Delta ^n$ are such that $\partial
^j(x_0,...,x_{n-1})=(x_0,...,x_{j-1},0,x_j,...,x_{n-1})$ and
$s^j(x_0,...,x_{n+1})=(x_0,...,x_{j-1},x_j+x_{j+1},x_{j+2},...,x_{n+1})$
in barycentric coordinates.

\par A $C^{\infty }$ or $H^{t'}_p$ space structure on the geometric
realization $|N_.|$ of $N_.$ consists of all continuous $C^{\infty
}$ $\bf R$-valued or $H^{t'}_p$ ${\cal A}_r$ valued functions $f$ on
$|N_.|$ respectively, that is the composition $\coprod_{n\ge 0}
(\Delta ^n -
\partial \Delta ^n)\times N_n\hookrightarrow
\coprod_{n\ge 0} \Delta ^n\times N_n \stackrel{q}{\longrightarrow }
|N_.| \stackrel{f}{\longrightarrow } {\cal A}_r$ is either
$C^{\infty }$ or $H^{t'}_p$, where $q$ denotes the quotient mapping,
$r=0$ or $1\le r\le 3$ correspondingly, ${\cal A}_0 = {\bf R}$,
${\cal A}_1 = {\bf C}$, ${\cal A}_2 = {\bf H}$, ${\cal A}_3 = {\bf
O}$.

\par {\bf 6. Proposition.} {\it If a group $G$ satisfies Conditions
4$(A1,A2,C1,C2)$, then sets $AG$ and $BG$ can be supplied with group
structures and they are twisted for $2\le r\le 3$. If $G$ is a
topological Hausdorff or $H^t_p$ differentiable alternative for
$r=3$ or associative for $0\le r\le 2$ group, then $AG$ and $BG$ are
topological Hausdorff or $C^{\infty }$ or $H^t_p$ differentiable
alternative for $r=3$ or associative for $0\le r\le 2$ groups
respectively.}
\par {\bf Proof.} Define on $AG$ and $BG$ group structures.
Introduce a homeomorphism pairing: $\Delta ^n\times \Delta ^k\to
\Delta ^{n+k}$, where $\sigma $ is a permutation of the set $\{
1,2,...,n+m+1 \} $ such that $t_{\sigma (1)}\le t_{\sigma (2)}\le
...\le t_{\sigma (n+k+1)}$, $\sigma \in S_{n+k+1}$, $S_m$ denotes
the symmetric group of all permutations of the set $ \{ 1,...,m \}
$. Define the multiplication for pure states in $AG$:
\par $(1)$ $|t_1,...,t_n,h_0[h_1|...|h_n]| *
|t_{n+1},...,t_{n+k+1}, h_{n+k+2}[h_{n+1}|...|h_{n+k+1}]| : =$ \\ $
|t_{\sigma (1)},...,t_{\sigma (n+k+1)}, (-1)^{q(\sigma )}
(h_0h_{n+k+2})[h_{\sigma (1)}|...|h_{\sigma (n+k+1)}]|$, \\ where
$h_l = y_li_{j(l)}$, $y_l\in G_{j(l)}$ for each $l=0,...,2^r-1$,
$q(\sigma )\in \bf Z$ is such that \par $(-1)^{q(\sigma )}
i_{j(0)}(i_{j(1)}...(i_{j(n+k+1)} i_{j(n+k+2)})...) = (i_{j(\sigma
(0))} i_{j(\sigma (n+k+2))})(i_{j(\sigma (1))}...$ \\ $(i_{j(\sigma
(n+k)))} i_{j(\sigma (n+k+1))})...)$ in ${\cal A}_r$; while in $BG$:

\par $(2)$ $|t_1,...,t_n,[h_1|...|h_n]| * |t_{n+1},...,t_{n+k+1},
[h_{n+1}|...|h_{n+k+1}]| : =$ \\ $ |t_{\sigma (1)},...,t_{\sigma
(n+k+1)}, (-1)^{p(\sigma )} [h_{\sigma (1)}|...|h_{\sigma
(n+k+1)}]|$, \\ where $h_l=y_li_{j(l)}$, $y_l\in G_{j(l)}$,
$p(\sigma )\in \bf Z$ is such that \par $(-1)^{p(\sigma )}
i_{j(1)}(i_{j(2)}...(i_{j(n+k)}i_{j(n+k+1)})...) =$
\par $ i_{j(\sigma (1))}(i_{j(\sigma (2))}...(i_{j(\sigma (n+k))}
i_{j(\sigma (n+k+1))})...)$ in ${\cal A}_r$.

\par Define also an addition in $A{\hat G}$:
\par $(1')$ $|t_1,...,t_n,h_0[h_1|...|h_n]| +
|t_{n+1},...,t_{n+k+1}, h_{n+k+2} [h_{n+1}|...|h_{n+k+1}]|$ \\ $ :=
|t_{\sigma (1)},...,t_{\sigma (n+k+1)}, (-1)^{q(\sigma )}
(h_0+h_{n+k+2}) [h_{\sigma (1)}|...|h_{\sigma (n+k+1)}]|$
\\ with $j(0) = j(n+k+2)$; as well as an addition in $B{\hat G}$:
\par $(2')$ $|t_1,...,t_n, [h_1|...|h_n]| + |t_{n+1},...,t_{n+k+1},
[h_{n+1},...,h_{n+k+1}]|$ \\  $= |t_{\sigma (1)},...,t_{\sigma
(n+k+1)}, (-1)^{p(\sigma )} [h_{\sigma (1)}|...|h_{\sigma (n+k+1)}]|$ \\
for pure states $h_l = y_li_{j(l)}$, $y_l\in {\hat G}_l$ for each
$l=0,...,2^r-1$. The multiplications $(1,2)$ extend to that of rings
in the natural way, when some pure states are zero, hence due to the
distributivity on the entire ring as well.

\par Since $\hat G$ is the ring, then these multiplications have unique
extensions on $AG$ and $BG$. Verify, that $AG$ and $BG$ become
groups with multiplications $(1)$ and $(2)$ respectively.

\par Due to $(1,2)$ we  get
\par $(3)$ $v * conj (v) = |t_1,...,t_k,(h_0 conj(h_0))[(h_1 conj
(h_1))|...|(h_k conj (h_k))]|$ for each
$v=|t_1,...,t_k,h_0[h_1|....|h_k]|$ in $AG$, while \par $w * conj
(w) = |t_1,...,t_k,[(h_1 conj (h_1))|...|(h_k conj (h_k))]|$
\\ for each $w=|t_1,...,t_k,[h_1|....|h_k]|$ in $BG$, where $h conj
(h)\in {\hat G}_0$ for each $h\in G$, but ${\hat G}_0$ is the center
of the ring $\hat G$: $ab=ba$ for each $a\in {\hat G}_0$ and $b\in
\hat G$. The Moufang identities in ${\cal A}_r$ for $r=3$ (see
\cite{harvey}) induces that of in $G$ such that
\par $(4)$ $(xyx)z=x(y(xz))$ and $(x^{-1}yx)z=x^{-1}(y(xz))$;
\par $(5)$ $z(xyx) = ((zx)y)x$ and $z(x^{-1}yx) = ((zx^{-1})y)x$;
\par $(6)$ $(xy)(zx)=x(yz)x$ and $(x^{-1}y)(zx)=x^{-1}(zy)x$,
since \par $(7)$ $x^{-1} = conj (x) (x~ conj (x))^{-1}$, \\
where $(x conj (x))\in G_0$.

\par The unit element in $AG$ is \par ${\bf e} :=
\{ |t_1,...,t_k,e[e|...|e]|\in (AG)_k: k=0,1,... \} $, \\ where
$i_0=1$, since \par $|t_1,...,t_n,h_0[h_1|...|h_n]| *
|t_{n+1},...,t_{n+k+1}, e[e|...|e]| = $ \\ $|t_{n+1},...,t_{n+k+1},
e[e|...|e]| * |t_1,...,t_n,h_0[h_1|...|h_n]| =$ \\
$|t_1,...,t_n,h_0[h_1|...|h_n]|$ due to equivalence relations
4$(2)$, $(1,...,1, e[e|...|e])\in |t_1,...,t_k, e[e|...|e]|$.

\par The ring $\hat G$ is $\bf Z_2$ graded in the sense that elements
$y_lj_l\in {\hat G}_lj_l$ are even for $l=0$ and odd for
$l=1,...,2^r-1$: $(y_0i_0)(y_ly_l)= (y_li_l)(y_0i_0) =
(y_0y_l)i_l\in {\hat G}_li_l$ for each $0\le l\le 2^r-1$,
$(y_li_l)^2= - y_l^2i_0\in {\hat G}_0i_0$, $(y_li_l)(y_ki_k)=
-(y_ki_k)(y_li_l)= (y_ly_k)i_s \in {\hat G}_si_s$ for $1\le l\ne
k\le 2^r-1$, where $i_s= i_li_k$. For each pure states
$g_0,...,g_k\in {\hat G}$ their product $(...(g_0g_1)g_2...)g_k$ is
a pure state, consequently, sets $AG$ and $BG$ are $\bf Z_2$ graded
analogously to $\hat G$ having even and odd elements such that \par
$(8)$ $A{\hat G} = (A{\hat G}_0)i_0\oplus (A{\hat G}_1)i_1\oplus ...
\oplus (A{\hat G}_{2^r-1})i_{2^r-1}$ and \par $(9)$ $B{\hat G} =
(B{\hat G}_0)i_0\oplus (B{\hat G}_1)i_1\oplus ...\oplus (B{\hat
G}_{2^r-1})i_{2^r-1}$. Each $AG_j$ and $BG_j$ is an associative
topological Hausdorff or $H^t_p$ differentiable group isomorphic
with $AG_0$ or $BG_0$ correspondingly for each $j$, since $G_j$ are
commutative and associative (see also Appendix B4 in
\cite{gajinvm}), where $G_0$ denotes the multiplicative group of the
ring ${\hat G}_0$. Therefore, $AG$ and $BG$ are the multiplicative
groups of the rings $A{\hat G}$ and $B{\hat G}$.

If $a\in AG_0$ or $a\in BG_0$, then $ab=ba$ for each $b\in AG$ or
$BG$ respectively. From Definition 6 it follows, that they are
$C^{\infty }$ or $H^{t}_p$ groups, when such is $G$.
\par The inverse element is
\par $(10)$ $\{ |t_1,...,t_k,h_0[h_1|...|h_k]|: k \}^{-1} =
\{ |t_1,...,t_k,h_0^{-1}[h_1^{-1}|...|h_k^{-1}]|: k \} $ \\ due to
$(2,6)$, since \par $(h_0(h_1...(h_{k-1}h_k)...))
((...(h_k^{-1}h_{k-1}^{-1})...h_1^{-1})h_0^{-1}) =$ \\
$(...(((h_0h_0^{-1})(h_1h_1^{-1}))(h_2h_2^{-1})...)(h_kh_k^{-1})=e$
\\ for pure states for each $k$ in view of the Moufang identities
$(4-6)$. In general it follows from $(3,8,9)$, since $v*conj (v)\in
AG_0$ or $BG_0$ for $v\in AG$ or $v\in BG$ respectively, hence
$v^{-1}= conj (v)(v*conj (v))^{-1} = \{
|t_1,...,t_k,h_0[h_1|...|h_k]|: k \}^{-1}$.
\par In view of $(8,9)$ and the existence of an inverse element
we get, that $AG$ is alternative, since $1\le r \le 3$. Putting
$h_0=1$ and applying the equivalence relation $\cal Y$ we get, that
$BG$ is also an alternative group, since the multiplicative group $
\{ i_0,...,i_7 \} $ is alternative. If $G$ is associative, for
example, when $1\le r\le 2$, then $AG$ and $BG$ are associative,
since the multiplicative group $ \{ i_0, i_1, i_2, i_3 \} $ is
associative.

\par Thus, groups $AG$ and $BG$
are $\bf Z_2$ graded, hence they are twisted over $ \{
i_0,...,i_{2^r-1} \} $. Consider for $\hat G$ multiplication and
addition operations, then they induce them for $AG$ and $BG$ as
above. It follows, that $E{\hat G}$ and $B{\hat G}$ are twisted
rings.

\par {\bf 7. Corollary.} {\it Let suppositions of Proposition
6 be satisfied, then $AB^mG$ and $B^mG$ are topological or
$C^{\infty }$ or $H^t_p$ differentiable groups respectively for each
$m\ge 1$. Moreover, all maps in the short exact sequence $e\to
B^aG\to AB^aG\to B^{a+1}G\to e$ are continuous or $C^{\infty }$ or
$H^t_p$ correspondingly.}
\par {\bf Proof.} Define differentiable space structure by
induction. Suppose that it is defined on $B^aG$ and $\Delta ^k\times
(B^aG)^m$ for $k, m\ge 0$, where $a\ge 1$. Then $f: AB^aG\to {\cal
A}_r$ is $C^{\infty }$ or $H^t_p$ if the composition $\coprod_{n\ge
0} (\Delta ^n - \partial \Delta ^n)\times (B^aG)^{n+1}
\stackrel{q_A}{\longrightarrow } AB^aG \stackrel{f}{\longrightarrow
} {\cal A}_r$ is either $C^{\infty }$ or $H^{t'}_p$, while $f:
B^{a+1}G\to {\cal A}_r$ is $C^{\infty }$ or $H^t_p$ if the
composition $\coprod_{n\ge 0} (\Delta ^n - \partial \Delta ^n)\times
(B^aG)^n \stackrel{q_B}{\longrightarrow } B^{a+1}G
\stackrel{f}{\longrightarrow } {\cal A}_r$ is either $C^{\infty }$
or $H^{t'}_p$, where $0\le r\le 3$.
\par A function $f: \Delta ^k\times (B^{a+1}G)^m\to {\cal A}_r$ is
$C^{\infty }$ or $H^t_p$ if the composition $\Delta ^k\times
(\coprod_{n\ge 0} (\Delta ^n - \partial \Delta ^n)\times (B^aG)^n)^m
\stackrel{id\times (q_B)^m}{\longrightarrow } \Delta ^k\times
(B^{a+1}G)^m \stackrel{f}{\longrightarrow } {\cal A}_r$ is either
$C^{\infty }$ or $H^{t'}_p$. From this it follows that all maps in
the short exact sequences are of the same class of smoothness.\par
Then the mappings $B^aG\times B^aG\to B^aG$ and $AB^aG\times
AB^aG\to AB^aG$ of the form $(f,g)\mapsto fg^{-1}$ are $C^0$ or
$C^{\infty }$ or $H^t_p$ in respective cases due to Formulas
6$(1-3,8-10)$ (see also \S 1.3.2 \cite{luwrgfbqo} and \S 1 and
Appendix B in \cite{gajinvm}).
\par {\bf 8. Corollary.} {\it If a group $G$ satisfies Conditions
4$(A1,A2,C1,C2)$, then there exist $H^t_p$ groups $AB^a(W^{M,\{
s_{0,q}: q=1,...,k \} }E)_{t,H}$ and $B^a(W^{M, \{ s_{0,q}:
q=1,...,k \} }E)_{t,H}$ for each $a\in \bf N$.}
\par {\bf Proof.} The wrap group $(W^{M,\{ s_{0,q}: q=1,...,k \} }
E;N,G,{\bf P})_{t,H}$ is a principal $G^k$ bundle over $(W^{M,\{
s_{0,q}: q=1,..,k \} }N)_{t,H}$, where $(W^{M,\{ s_{0,q}: q=1,...,k
\} }N)_{t,H}$ is commutative and associative (see Proposition 7(1) in
\cite{lulaswgof}).
\par If $g\in {\hat G}$, then $g$ has the decomposition
$g=g_0i_0+...+g_{2^r-1}i_{2^r-1}$ with $g_j\in {\hat G}_j$ for each
$j=0,1,...,2^r-1$ and
\par $(1)$ $g_0 = (g + (2^r-2)^{-1} \{ - g + \sum_{s=1}^{2^r-1}
i_s(gi_s^*) \} )/2$ and
\par $(2)$ $g_j = (i_j(2^r-2)^{-1} \{ - g +
\sum_{s=1}^{2^r-1} i_s(gi_s^*) \} - gi_j )/2$ \\ for each
$j=1,...,2^r-1$. Therefore, each $g_0,...,g_{2^r-1}$ has analytic
expressions through $g$ due to Formulas $(1,2)$. Fix this
representations. Then the $H^t_p$ differentiable parallel transport
structure $\bf P$ with the groups $G$ induces the $H^t_p$
differentiable parallel transport structures $\mbox{ }_j{\bf P}$
with groups $G_j$. \par Since ${\hat G}^k$ is isomorphic with
$\bigoplus_{0\le j(1),...,j(k)\le 2^r-1} ({\hat
G}_{j(1)}i_{j(1)},...,{\hat G}_{j(k)}i_{j(k)})$ which is isomorphic
with $\bigoplus_{0\le j(1),...,j(k)\le 2^r-1} {\hat
G}_0^k(i_{j(1)},...,i_{j(k)})$, then \\ $(W^{M,\{ s_{0,q}:
q=1,...,k \} }E;N,G,{\bf P})_{t,H}$ is isomorphic with a group \\
$\{ f=(f_1,...,f_k)\in \bigoplus_{0\le j(1),...,j(k)\le 2^r-1}
[(W^{M,\{ s_{0,q}: q=1,...,k \} }E;N,G_0,{\bf P})_{t,H}\cup \{ 0
\}]$ \\ $ (i_{j(1)},...,i_{j(k)}): f_1\ne 0,...,f_k\ne 0 \} $, where
$(i_{j(1)},...,i_{j(k)}) \in ({\cal A}_r^*)^k$ and $({\cal
A}_r^*)^k$ has the embedding into the family of all $k\times k$
matrices with entries in ${\cal A}_r$ as diagonal matrices,
$(W^{M,\{ s_{0,q}: q=1,...,k \} }E;N,G_j,{\bf P})_{t,H}$ is
commutative for each $j=0,...,2^r-1$ due to Theorem 6
\cite{luwrgfbqo}. \par The construction of Proposition 5 above has
the natural generalization for $G^k$ instead of $G$ such that
\par $A{\hat G}^k= \bigoplus_{0\le j(1),...,j(k)\le 2^r-1}
(A{\hat G}_{j(1)}i_{j(1)},...,A{\hat G}_{j(k)}i_{j(k)})$ \\ which is
isomorphic with $\bigoplus_{0\le j(1),...,j(k)\le 2^r-1} A{\hat
G}_0^k(i_{j(1)},...,i_{j(k)})$, also \\ $B{\hat G}^k$ is isomorphic
with $\bigoplus_{0\le j(1),...,j(k)\le 2^r-1} B{\hat
G}_0^k(i_{j(1)},...,i_{j(k)})$, consequently, \\ $A(W^{M,\{ s_{0,q}:
q=1,...,k \} }E;N,G,{\bf P})_{t,H}$ is isomorphic with a group \\
$\{ v\in \bigoplus_{0\le j(1),...,j(k)\le 2^r-1} [A(W^{M,\{ s_{0,q}:
q=1,...,k \} }E;N,G_0,{\bf P})_{t,H}^k \cup \{ 0 \} ]
(i_{j(1)},...,i_{j(k)}): v_n = |t_1,...,t_n,h_0[h_1|...|h_n]|,
h_j\ne 0 \forall j, \forall n \} $ and
\\ $B(W^{M,\{ s_{0,q}: q=1,...,k \} }E;N,G,{\bf P})_{t,H}$ is
isomorphic with \\ $\{ v\in \bigoplus_{0\le j(1),...,j(k)\le 2^r-1}
[B(W^{M,\{ s_{0,q}: q=1,...,k \} }E;N,G_0,{\bf P})_{t,H}^k\cup \{ 0
\} ] (i_{j(1)},...,i_{j(k)}): v_n = |t_1,...,t_n, [h_1|...|h_n]|,
h_j\ne 0 \forall j, \forall n \} $. Continuing this by induction on
$a$ and using Corollary 7 we get the statement of this corollary for
each $a\in \bf N$.

\par {\bf 9. Lemma.} {\it Let $N$ be a $C^{\infty }$ or $H^t_p$ manifold
over ${\cal A}_r$ with $0\le r \le 3$ and $G$ a $C^{\infty }$ or
$H^t_p$ differentiable group. If $f: N\to BG$ is a mapping such that
for each $y\in N$ there exists an open neighborhood $V$ of $y$ in
$N$ such that $f|_V=|f_0,f_1,...,f_n,[g_1|...|g_n]|$ with
$f_0,...,f_n$ being $C^{\infty }$ or $H^t_p$ differentiable
mappings, then $f$ is either $C^{\infty }$ or $H^t_p$ differentiable
mapping correspondingly.}
\par {\bf Proof.} If $h: BG\to {\cal A}_r$ is a $C^{\infty }$
or $H^t_p$ mapping, then for each $n\ge 1$ the composition $\Delta
^n\times G^n \stackrel{q_B}{\longrightarrow } BG
\stackrel{h}{\longrightarrow } {\cal A}_r$ is of the corresponding
class. For the commutative diagram consisting of
$N\stackrel{f}{\longrightarrow } BG \stackrel{h}{\longrightarrow }
{\cal A}_r$ and $N\stackrel{\bar f}{\longrightarrow } \Delta
^n\times G^n \stackrel{q_B}{\longrightarrow } BG$ and $f=q_b\circ
{\bar f}$, where ${\bar f} := (f_0,...,f_n,h_1,...,h_n)$ both $\bar
f$ and $h\circ q_B$ are continuous $C^{\infty }$ or $H^t_p$. Then
the composition $h\circ f=h\circ q_B\circ {\bar f}$ is continuous
and either $C^{\infty }$ or $H^t_p$, where as usually $h\circ f(y):=
h(f(y))$. Thus $f: N\to BG$ is continuous either $C^{\infty }$ or
$H^t_p$ respectively.

\par {\bf 10. Twisted bar resolution and hypercohomologies.}
For a twisted group $G$ satisfying Conditions
4$(A1,A2,C1,C2)$ the composition of the short exact sequences \par
$(1)$ $e\to B^aG\to AB^aG\to B^{a+1}G\to e$ induces the long exact
sequence
\par $(2)$ $e\to G\to AG \stackrel{\sigma }{\longrightarrow}
ABG \stackrel{\sigma }{\longrightarrow} AB^2G \stackrel{\sigma
}{\longrightarrow} ... \stackrel{\sigma }{\longrightarrow} AB^aG
\stackrel{\sigma }{\longrightarrow} ...$, \\
where for each $a\ge 0$ the homomoprhism $\sigma : AB^aG\to
AB^{a+1}G$ is the composition $AB^aG\to B^{a+1}G\to AB^{a+1}G$ of
the surjection $AB^aG\to B^{a+1}G$ and the monomorphism $B^{a+1}G\to
AB^{a+1}G$. \par In view of Corollary 7 the short exact sequence
$(2)$ is a $C^{\infty }$ or $H^t_p$ $B^aG$-extension of $B^{a+1}G$.
Hence the long exact sequence $(2)$ induces the long exact sequence
of twisted sheaves \par $(3)$ $e\to G_N\to AG_N \stackrel{\sigma
}{\longrightarrow} ABG_N \stackrel{\sigma }{\longrightarrow} AB^2G_N
\stackrel{\sigma }{\longrightarrow} ... \stackrel{\sigma
}{\longrightarrow} AB^aG_N \stackrel{\sigma }{\longrightarrow} ...
$, \\ which we will call the (twisted) bar resolution of the sheaf
$G_N$, where $G_N$ denotes the sheaf of $C^{\infty }$ or $H^t_p$
functions on $N$ with values in $G$.

\par Suppose that ${\cal S}^*$ and ${\cal F}^*$ are complexes of sheaves
of ${\cal G}$-modules, where $\cal G$ is a sheaf of rings, where
${\cal S}^*$ and ${\cal F}^*$ and $\cal G$ may be simultaneously
twisted over $ \{ i_0,...,i^{2^r-1} \} $. Then a homomorphism
mapping $\sigma : {\cal S}^*\to {\cal F}^*$ of such complexes
induces a mapping of cohomology sheaves ${\sf H}^j(\sigma ): {\sf
H}^j({\cal S}^*)\to {\sf H}^j({\cal F}^*)$, where ${\sf H}^j({\cal
S}^*)$ is the sheaf associated with the pre-sheaf $U\mapsto [ker
(\Gamma (U,{\cal F}^j))\to \Gamma (U,{\cal F}^{j+1}))]/ im [(\Gamma
(U,{\cal F}^{j-1}))\to \Gamma (U,{\cal F}^j))]$, where $\Gamma
(U,{\cal S}^j)$ denotes the group of sections of the sheaf ${\cal
S}^j$ for a subset $U$ open in $X$ (see also \S 1). Then $\sigma $
is called a quasi-isomorphism, if ${\sf H}^j(\sigma )$ is an
isomorphism for each $j$. \par We consider complexes bounded below,
that is there exists $j_0$ such that ${\cal S}^j =0 $ for each
$j<j_0$.
\par A mapping $\sigma : {\cal S}^*\to {\cal T}^*$ is called an
injective resolution of ${\cal S}^*$ if ${\cal T}^*$ is a complex of
$\cal G$-modules bounded below, $\sigma $ is a quasi-isomorphism and
the sheaves ${\cal T}^b$ are injective, which means that $Hom ({\cal
B},{\cal T}^b)\to Hom ({\cal K},{\cal T}^b)$ is surjective for each
injective mapping ${\cal K}\to \cal B$ of sheaves of $\cal
G$-modules.

\par Let ${\cal G}$ be a constant sheaf of rings, may be twisted
over $ \{ i_0,...,i_{2^r-1} \} $. Suppose that ${\cal S}^*$ is a
complex of ${\cal G}$-modules bounded below. The hypercohomology
group $\mbox{ }_h{\sf H}^b(X,{\cal S}^*)$ is defined to be the
${\cal G}$-module such that \par $\mbox{ }_h{\sf H}^b(X,{\cal S}^*)
:= [ker (\Gamma (X,{\cal T}^b)\to \Gamma (X,{\cal T}^{b+1})]/ [im
(\Gamma (X,{\cal T}^{b-1})\to \Gamma (X,{\cal T}^b)]$. \\
If $\sigma : {\cal S}^*\to {\cal F}^*$ is a quasi-isomorphism, then
$\sigma $ induces an isomorphism of the hypercohomology groups:
\par $\sigma : \mbox{ }_h{\sf H}^b(X,{\cal S}^*) \cong \mbox{ }_h{\sf
H}^b(X,{\cal F}^*)$ (see also \cite{gajinvm} and the reference
$[EV]$ in it). In view of Lemma 16 \cite{lulaswgof} the
hypercohomology groups $\mbox{ }_h{\sf H}^b(X,{\cal S}^*)$ are
twisted over $ \{ i_0,...,i_{2^r-1} \} $, when ${\cal S}^*$ and
$\cal G$ are twisted over $\{ i_0,...,i_{2^r-1} \} $.

\par {\bf 11. Proposition.} {\it The sequence 10$(3)$ is an acyclic
resolution of the sheaf $G_N$.}
\par {\bf Proof.} Each standard simplex $\Delta ^n$ with $n\ge 1$
has a $C^{\infty }$ retraction ${\hat z}: \Delta ^n\times [0,1]\to
\{ y \} $ into a point $y$ belonging to it. \par There exists a
$C^{\infty }$ deformation retraction
\par $(1)$ ${\hat f}: AG\times [0,1]\to AG$ supplied by the family
of mappings
\par $(2)$ ${\hat f}_n: (AG)_n\times [0,1]\to (AG)_{n+1}$, where
\par $(3)$ $(AG)_n := q_A(\coprod_{j\le n} \Delta ^j\times
G^{j+1})\subset AG$ and
\par $(4)$ ${\hat f}_n(|t_1,...,t_n,h_0[h_1|...|h_n]|,t) :=
|\Phi (0,t), \Phi (t_1,t),...,\Phi (t_n,t), h_0[h_1|...|h_n]|$,
where $\Phi : [0,1]^2\to [0,1]$ is defined as the composition $\Phi
(x,t) := \phi (\min (1,x+t))$ taking $\phi $ a smooth nondecreasing
function $\phi : [0,1]\to [0,1]$ such that $\phi (0)=0$ and $\phi
(1)=1$. \par Then for each $C^{\infty }$ or $H^t_p$ differentiable
mapping $v: AG\to {\cal A}_r$ we get $v\circ {\hat f}\circ
(q_n\times id) = v^{n+1}\circ {\hat h}_n$ and $v\circ q_A =
v^{n+1}$, where $q_n\times id: (\Delta ^n\times G^{n+1})\times
[0,1]\to AG\times [0,1]$, ${\hat h}_n: (\Delta ^n\times
G^{n+1})\times [0,1]\to \Delta ^{n+1}\times G^{n+2}$ is the smooth
mapping given by the formula ${\hat h}_n(t_1,...,t_n,
g_0,g_1,...,g_n,t) = (\Phi (0,t), \Phi (t_1,t),...,\Phi (t_n,t), e,
g_0, g_1,...,g_n)$. \par At the same time $\Delta ^{n+1}$ has a
$C^{\infty }$ retraction onto $\Delta ^n$ for each $n\ge 0$ while
the group $G$ is $C^{\infty }$ or $H^t_p$ differentiable and arcwise
connected. Therefore, for each $b>0$ the cohomology group ${\sf
H}^b(N,AG_N)$ is trivial (see also \S 2.2 \cite{bredon} and \S 54
below).

\par {\bf 12.}  A sheaf $\cal S$ on $X$ is called soft if for each
closed subset $Y$ of $X$ the restriction map ${\cal S}(X)\to {\cal
S}(Y)$ is a surjection.
\par {\bf 12.1. Lemma.} {\it For a $C^{\infty }$ or $H^t_p$
differentiable group $G$ satisfying conditions 4$(A1,A2,C1,C2)$ the
sheaf $AG_N$ is soft.}

\par {\bf Proof.} Consider a closed subset $Y$ of $N$ and a section
$\sigma _Y$ of $AG_N$ over $Y$. In accordance with the definition of
a section over a closed subset there exists an open set $U$ and an
extension $\sigma _U$ of $\sigma _Y$ from $Y$ onto $U$ such that
$Y\subset U\subset N$. From the paracompactness of $N$ there exists
a neighborhood $V$ of $Y$ such that $cl (V)\subset U$, where $cl
(V)$ denotes the closure of $V$ in $N$. Therefore, the extension
$\sigma $ of $\sigma _Y$ to a global section of $AG_N$ is provided
by the formula $\sigma (x) = {\hat f} (\sigma _U(x), \psi (x))$,
where ${\hat f}: AG\times [0,1]\to AG$ is a deformation retraction
(see \S 11) and $\psi : N\to [0,1]$ is either a $C^{\infty }$ or
$H^t_p$ differentiable function equal to $1$ on $V$ and equal to $0$
on $M\setminus U$.

\par {\bf 13. Remark.} Let now $N$ be a $C^{\infty }$ or $H^{\infty
}$ manifold over ${\cal A}_r$ and $T(N,G,\pi ,\Psi )$ be a tangent
bundle with $T=TN$ and the projection $\pi : T\to N$, where
connecting mappings $\phi _j\circ \phi _k^{-1}$ for $V_j\cap V_k\ne
\emptyset $ of the atlas $At (N) = \{ (V_j, \phi _j) : j \} $ of the
manifold $N$ are ${\cal A}_r$ holomorphic for $1\le r\le 3$, $\phi
_j\circ \phi _k^{-1}\in H^{\infty }$. Denote by $\cal T$ the sheaf
of germs of smooth sections of $T$. Then $B\cal T$ denotes the sheaf
associated with a the pre-sheaf assigning to each open subset $V$ of
$N$ the group of sections of the natural projection $\coprod_{y\in
V} B(\pi ^{-1}(y))\to V$, which are locally of the form
\par $(1)$ $y\mapsto |t_1(y),...,t_n(y),[\sigma _1(y)|...|\sigma _n(y)]|$,
where $t_1,...,t_n$ are $C^{\infty }$ for $r=1$ or $H^{\infty }$ for
$1\le r\le 3$ functions and $\sigma _1,...,\sigma _n$ are $C^{\infty
}$ or $H^{\infty }$ sections of the vector bundle $T(N,G,\pi ,\Psi
)$. Using constructions above we define $B^{a+1}\cal T$ for each
$a\in \bf N$ by induction. \par Then $B^{a+1}\cal T$ is the sheaf
associated with the pre-sheaf assigning to an open subset $V$ of $N$
the group of sections of the natural projection \par $\coprod_{y\in
V} B^{a+1}(\pi ^{-1}(y))\to V$ having the local form $(1)$, where
$\sigma _1,...,\sigma _n$ are sections of $B^a\cal T$ over $V$.
Similarly we define $AB^a\cal T$.
\par If now $T = \Lambda ^bT^*N$ is the $b$-th exterior power of
the cotangent bundle of $N$, then the above construction produces
the sheaves $AB^a{\cal S}^b_{N,{\cal A}_r}$ and $B^{a+1}{\cal
S}^b_{N,{\cal A}_r}$ of $AB^a{\cal A}_r$ and $B^{a+1}{\cal A}_r$
valued respectively differential $C^{\infty }$ forms on $N$, where
the index ${\cal A}_r$ may be omitted, when the Cayley-Dickson
algebra ${\cal A}_r$ is specified. In the equation
\par $(2)$ $w = \sum_J f_J(z) dx_{b_1,j_1}\wedge dx_{b_2,j_2}\wedge
... \wedge dx_{b_k,j_k}$,
\\ where $f_J: N \to AB^a{\cal A}_r$ or
$f_J: N\to B^{a+1}{\cal A}_r$, $z = (z_1,z_2,...)$ are local
coordinates in $N$, $z_b =x_{b,0}i_0 + x_{b,1} i_1 +...+
x_{b,{2^r-1}}i_{2^r-1}$, where $z_b\in {\cal A}_r$, $x_{b,j}\in \bf
R$ for each $b$ and every $j=0,1,...,2^r-1$, $J =(b_1,j_1;
b_2,j_2;...; b_k,j_k)$.
\par Since each topological vector space $Z$ over ${\cal A}_r$
with $2\le r\le 3$ has the natural twisted structure $Z=
Z_0i_0\oplus Z_1i_1\oplus ... \oplus Z_{2^r-1}i_{2^r-1}$ with
pairwise isomorphic topological vector spaces $Z_0,...,Z_{2^r-1}$
over $\bf R$, then $TN$ and $T^*N$ and $\Lambda ^bT^*N$ have twisted
structures, where $X^*$ denotes the space of all continuous ${\cal
A}_r$ additive and $\bf R$ homogeneous functionals on $X$ with
values in ${\cal A}_r$, when $2\le r\le 3$, while $X^*$ over $\bf C$
is the usual topologically dual space of continuous $\bf C$-linear
functionals on $X$. Therefore, due to Proposition 6 $B^a\cal T$ and
$AB^a\cal T$ have the induced twisted structure for each $a\in \bf
N$.
\par Each section $\sigma $ of the sheaf $AB^a{\cal S}^k_N$ can be
written in the form: $\sigma = |h_0,..,h_n,\sigma _0,...,\sigma
_n|$, where $\sigma _0,...,\sigma _n$ are smooth differential
$B^a{\cal A}_r$ valued differential $k$-forms on $V$ and $ \{ h_j:
j=0,1,2,... \} $ is a $C^{\infty }$ smooth partition of unity on
$V$. An addition of differential forms induces the additive group
structure. As a multiplication there can be taken the external
product of differential forms which is twisted over ${\cal A}_r$ for
$2\le r\le 3$.
\par The group of sections of the sheaf $AB^b{\cal S}_N^k$ for an
open subset $V$ in $N$ denote by $\Gamma (V,AB^b{\cal S}_N^k)$. The
sequence of the groups $0\to \Gamma (V,{\cal S}_N^k)\to \Gamma
(V,A{\cal S}_n^k)\to \Gamma (V,B{\cal S}_N^k)\to 0$ is exact for
each open subset $V$ in $N$, since the sequence of vector bundles
$0\to \Lambda ^kT^*N\to A\Lambda ^kT^*N\to B\Lambda ^kT^*N\to 0$ is
exact. \par The twisted structure of the groups $AB^b{\cal S}_N^k$
induces the twisted structure of $\Gamma (V,AB^b{\cal S}_N^k)$.
Therefore, the sequence of sheaves $0\to B^b{\cal S}_N^k\to
AB^b{\cal S}_N^k\to B^{b+1}{\cal S}_N^k\to 0$ is exact as well as
for each $b\ge 0$.
\par The composition of these sequences induces a long exact sequence
\par $(2)$ $0\to {\cal S}_N^k\to A{\cal S}_N^k  \stackrel{\sigma
}{\longrightarrow} AB{\cal S}_N^k \stackrel{\sigma
}{\longrightarrow} ... \stackrel{\sigma }{\longrightarrow} AB^b{\cal
S}_N^k\stackrel{\sigma }{\longrightarrow}... $, \\
where $\sigma : AB^b{\cal S}_N^k\to AB^{b+1}{\cal S}_N^k$ is the
composition of mappings $AB^b{\cal S}_N^k\to B^{b+1}{\cal S}_N^k\to
AB^{b+1}{\cal S}_N^k$. The sequence $(2)$ will be called the bar
resolution of the sheaf ${\cal S}_N^k$.

\par Let $\cal S$ be an arbitrary twisted sheaf on a topological space
$X$. Denote by $A\cal S$ and $B\cal S$ sheaves associated with the
pre-sheaves $V\mapsto A(\Gamma (V,{\cal S}))$ and $V\mapsto B(\Gamma
(V,{\cal S}))$ correspondingly. The stalks of $A\cal S$ and $B\cal
S$ are $A{\cal S}_x$ and $B{\cal S}_x$ at $x$, while the sequence
\par $(3)$ $e\to {\cal S}_x\to A({\cal S}_x)\to B({\cal S}_x)\to e$ \\
is exact, consequently, the sequence of sheaves
\par $e\to {\cal S}\to A{\cal S}\to B{\cal S}\to e$ \\
is also exact. The composition of these sequences gives the bar
resolution of $\cal S$
\par $(4)$ $e\to {\cal S}\to A{\cal S} \stackrel{\sigma
}{\longrightarrow} AB{\cal S} \stackrel{\sigma }{\longrightarrow}
... \stackrel{\sigma }{\longrightarrow} AB^b{\cal S}\stackrel{\sigma
}{\longrightarrow}... $. \par The complex of sheaves \par $(5)$
${\cal B}^*({\cal S}): A{\cal S} \stackrel{\sigma }{\longrightarrow}
AB{\cal S} \stackrel{\sigma }{\longrightarrow} ... \stackrel{\sigma
}{\longrightarrow} AB^b{\cal S}\stackrel{\sigma
}{\longrightarrow}... $ \\
is called the bar complex of $\cal S$. The bar resolution of $\cal
S$ is an acyclic resolution of $\cal S$ that is deduced analogously
to the proofs of Proposition 11 and Lemma 12.1. Thus the cohomology
of $\cal S$ is equal to the cohomology of the cochain complex
\par $(6)$ $\Gamma (N,A{\cal S}) \stackrel{\sigma }{\longrightarrow} \Gamma
(N,AB{\cal S}) \stackrel{\sigma }{\longrightarrow} ...
\stackrel{\sigma }{\longrightarrow} \Gamma (N,AB^b{\cal S})
\stackrel{\sigma }{\longrightarrow}... $. \\
The complex $(6)$ will be called the bar cochain complex of $\cal S$
and will be denoted by $C^*_B({\cal S})$.
\par Each short exact sequence of sheaves $e\to E\to F\to Y\to e$
twisted over generators $\{ i_0, i_1,...,i_{2^r-1} \} $, $2\le r\le
3$, induces a short exact sequence of complexes sheaves $e\to {\cal
B}^*(E)\to {\cal B}^*(F)\to {\cal B}^*(Y) \to e$, where ${\cal
B}^*(F): AF \stackrel{\sigma }{\longrightarrow} ABF \stackrel{\sigma
}{\longrightarrow} AB^2F \stackrel{\sigma }{\longrightarrow} ...
\stackrel{\sigma }{\longrightarrow} AB^bF \stackrel{\sigma
}{\longrightarrow} ...$ is the bar complex of $F$.

\par {\bf 14. Proposition.} {\it If a sequence of groups
$e\to K\to G\to J\to e$ is exact, where $E$ and $K, G, J$ are
arcwise connected, then the sequence $e\to (W^ME;N,K,{\bf P})_{t,H}
\to (W^ME;N,G,{\bf P})_{t,H}\to (W^ME;N,J,{\bf P})_{t,H}\to e$ is
exact.}
\par {\bf Proof.} In view of Proposition 7.1 \cite{lulaswgof}
$(W^ME;N,K,{\bf P})_{t,H}$ is the principal fiber bundle over
$(W^MN)_{t,H}$ with the structure group $K^k$, \par $\pi _{K,*}:
(W^ME;N,K,{\bf P})_{t,H}\to (W^MN)_{t,H}$, \\ $\pi
_{K,*}^{-1}<w_0>_{t,H} = <w_0>_{t,H} \times K^k=e\times K^k$, where
$e\in (W^MN)_{t,H}$ denotes the unit element. Since the sequence
$e\to K^k\to G^k\to J^k\to e$ is exact as well, then the
corresponding sequence of wrap groups is exact.

\par {\bf 15. Proposition.} {\it Let $G$ be a $C^{\infty }$ or
$H^t_p$ differentiable twisted group over $ \{ i_0, i_1,...,i_{2^r-1
} \} $ satisfying Conditions 4$(A1,A2,C1,C2)$. Then for each
$C^{\infty }$ or $H^t_p$ principal $G$ -bundle $E(N,G,\pi ,\Psi )$
there exists a $C^{\infty }$ or $H^t_p$ differentiable mapping $\phi
: N\to BG$ such that $E\to N$ is the pull-back of the universal
principal $G$ bundle by $\phi $.}
\par {\bf Proof.} Consider an open covering ${\cal V} = \{ V_j: j\in
J \} $ of $N$, where $J$ is a set, such that  for each $j\in J$
there exists a trivialization $\psi _j: \pi ^{-1} (V_j) \to
V_j\times G$. Define a mapping $g_j: E\to G$ by the formula $g_j(x)=
pr_2(\psi _j(x))$ for $x\in \pi ^{-1}(V_j)$, $g_j(x)=e$ for $x\notin
\pi ^{-1}(V_j)$, where $e$ denotes the neutral element in $G$ and
$pr_2: V_j\times G\to G$ is the projection on the second factor.
\par For a principal $G$ bundle $E(N,G,\pi ,\Psi )$ consider a family
of $H^{t'}_p$ transition functions $ \{ g_{i,j}: i, j \in J \} $
related with an open covering ${\cal V} := \{ V_j: j\in J \} $ of an
$H^{t'}_p$ manifold $N$ over ${\cal A}_r$, where $J$ is a set,
$g_{i,j}: V_i\cap V_j\to {\cal A}_r^*$, when the intersection
$V_i\cap V_j\ne \emptyset $ is non-void, $1\le r\le 3$, $n\in \bf
N$. Introduce the mapping
\par $(1)$ $g_{E,N}(x) := |f_{j(0)},f_{j(1)},...f_{j(n)},
[g_{j(0),j(1)}| g_{j(1),j(2)}|...|g_{j(n-1),j(n)}]|$ \\ such that
$g_{E,N}: N\to BG$, where $\{ f_j: j\in J \} $ is an $H^{t_1}_p$
partition of unity subordinated to $\cal U$ with $t'\le t_1\le
\infty $. Therefore, $g_{E,N}$ can be chosen of the smoothness class
$H^{t'}_p$. Thus, $E(N,G,\pi ,\Psi )$ is the pull-back of the
universal bundle $AG(BG,G,\pi ^A_B,\Psi ^A)$ by the classifying
mapping $g_{E,N}$, where $\pi ^A_B: AG\to BG$ is as in \S 4. Show it
in details.
\par Take a partition of unity of class $C^{\infty }$ or $H^t_p$
subordinated to to the covering $\cal V$ and $\Phi : A\to AG$ be the
following mapping \par $\Phi (y) := |f_{j_0}(\pi (y)), f_{j_1}(\pi
(y)),...,f_{j_n}(\pi (y)), g_{j_0}(y), g_{j_1}(y),...,g_{j_n}(y)|$,
where $j_0,...,j_n$ are indices such that $f_j(\pi (y))\ne 0$ for
each $j\in \{ j_0,...,j_n \} $. Then $\Phi $ is $G$ equivariant,
which means that $\Phi (yh)= \Phi (y)h$ for all $y$ and $h\in G$,
since $g_j(yh)=pr_2(\psi _j(yh))$ for $yh\in \pi ^{-1}(V_j)$ and
$g_j(yh)=e$ for $yh \notin \pi ^{-1}(V_j)$. Indeed, $y\in \pi
^{-1}(y)$ is equivalent to $yh\in \pi ^{-1}(V_j)$ for each $h\in G$,
since $\pi ^{-1}(V_j)=V_j\times G$, where $y=(u,q)$ with $u\in N$
and $q\in G$ and $(u,q)h=(u,qh)$ in local coordinates. Thus
$g_j(yh)=g_j(y)R_h$, where $R_h = h$ for $y\in \pi ^{-1}(V_j)$ and
$R_h=e$ for $y\notin \pi ^{-1}(V_j)$.
\par Therefore, $\Phi $ induces a morphism of principal $G$-bundles
$$
\begin{array}{ccc}
A&\stackrel{\Phi }{\longrightarrow}&AG\\
\downarrow\lefteqn{\pi }&&\downarrow \\
N&\stackrel{\phi }{\longrightarrow}&BG
\end{array}
$$
where the restriction of $\phi $ to $V_j$ is $\phi (x)|_{V_j} =
|f_{j_0}, f_{j_1}(x),...,f_{j_n}(x),[g_{j_0}(\sigma (x)):
g_{j_1}(\sigma (x)): ...: g_{j_n}(\sigma (x))]|$, $\sigma : V_j\to
\pi ^{-1}(V_j)$ is a smooth section of the restriction $\pi
^{-1}(V_j)\to V_j$ for $\pi : E\to N$. Consider equivalence classes
$q_j\sim g_j$ if and only if there exist $s_1,...,s_m\in G$ such
that $(s_m(s_{m-1}...(s_1(q_j)...)= g_j$, hence $q_jh\sim g_jh$,
since $q_jh\sim g_jh$ if and only if $h^{-1}q_j\sim h^{-1}g_j$,
which is equivalent with $h(s_m(s_{m-1}...(s_1(h^{-1}q_j)...) =
g_j$. Due to the alternativity of the group $G$ we get
$[(...(g_j^{-1}s_1^{-1})...s_{m-1}^{-1})s_m^{-1}] [(s_m(s_{m-1}...
(s_1g_l)...)]= g_j^{-1}g_l$.
\par Therefore, in the non-homogeneous
coordinates the mapping $\phi $ takes the form $\phi (x) =
|f_{j_0}(x), f_{j_1}(x),...,f_{j_n}(x),[g_{j_0,j_1}(x)|
g_{j_1,j_2}(x)|...| g_{j_{n-1},j_n}(x)]|$, where $g_{j,l}(x) =
[g_j(\sigma (x))]^{-1} g_l(\sigma (x))$ are transition functions
associated with the open covering of $N$ by the open sets $ \{ x\in
N: f_j(x)>0 \} $. Then the mapping $\phi (x)$ is independent from
the choice of $\sigma $, since $g_j$ are $G$-equivariant. All
functions $g_j$ are either $C^{\infty }$ or $H^t_p$, hence $\phi $
is either $C^{\infty }$ or $H^t_p$ correspondingly.

\par {\bf 16. Corollary.} {\it For each smooth $C^{\infty }$ or $H^t_p$
principal $B^b{\cal A}_r^*$ bundle with $1\le r\le 3$ there exists a
$C^{\infty }$ or $H^t_p$ differentiable mapping $\phi : N\to
B^{b+1}{\cal A}_r^*$ such that $E(N,G,\pi, \Psi )$ is the pull-back
of the universal principal $B^b{\cal A}_r^*$-bundle by $\phi $,
where $G=B^b{\cal A}_r^*$.}

\par {\bf 17. Lemma.} {\it Let $G$ be a differentiable (topological) group
twisted over $ \{ i_0,...,i_{2^r-1} \} $ for $1\le r\le 3$
satisfying Conditions 4$(A1,A2,C1,C2)$. Then the group of
isomorphism classes of $C^{\infty }$ smooth (continuous) principal
$G$-bundles over $N$ is isomorphic with the group $[N,BG]^{\infty }$
of smooth (or $[N, BG]^0$ of continuous respectively) homotopy
classes of smooth (continuous) maps from $N$ to $BG$.}
\par {\bf Proof.} Each principal $G$-bundle over $N$ has properties
4$(A1,A2,C1,C2)$ induced by that of $G$, where locally $\pi
^{-1}(V_j) =V_j\times G$, $conj (y) = (u, conj (g))$ for each $y=
(u,g)\in V_j\times G$, while $ \{ i_0,...,i_{2^r-1} \} $ for $2\le
r\le 3$ is the multiplicative group, which is associative for $r=2$
and alternative for $r=3$.
\par Consider a short exact sequence $e\to G_N\to AG_N\to BG_N\to
e$. In view of Lemma 16 \cite{lulaswgof} it induces the cohomology
long exact sequence
\par $... \to C^{\infty }(N,AG) \stackrel{\pi _*}\longrightarrow
C^{\infty }(N,BG)\to {\sf H}^1(N,G_N)\to {\sf H}^1(N,AG_N)\to ...$.
From ${\sf H}^1(N,AG_N)\cong e$ we get the isomorphism $C^{\infty
}(N,BG)/\pi _* C^{\infty }(N,AG) \cong {\sf H}^1(N,G_N)$. Then the
image $\pi _*C^{\infty }(N,AG)$ of the group $C^{\infty }(N,AG)$ in
$C^{\infty }(N,BG)$ consists of all smooth maps from $N$ to $BG$ for
which there exist lift mappings from $N$ to $AG$. \par On the other
hand, $f\in C^{\infty }(N,BG)$ has a lift $F: N\to AG$ if and only
if $f$ is smooth (or continuous) homotopic to a constant mapping,
since $[g_0:...:g_n]$ in $(BG)_n$ is the equivalence class $\{
(g_0,...,g_n)\sim (s_m(...(s_1g_0)...),...,(s_m...(s_1g_n)...)):
s_1,...,s_m\in G, m\in {\bf N} \} $, consequently, $C^{\infty
}(N,BG)/\pi _*C^{\infty }(N,AG) \cong [N,BG]^{\infty }$. In the
class of continuous mappings we get analogously $C^0(N,BG)/\pi
_*C^0(N,AG) \cong [N,BG]^0$.

\par {\bf 18. Notes.} In view of \S \S 4-6 there exists a short
exact sequence
\par $e\to G\to AG\to BG\to e$ \\ of $H^{t'}_p$ homomorphisms
due to the twisted structures of $G$, $AG$ and $BG$ (see Equations
4$(A2)$ and 6$(8,9)$).
\par To groups $AG$ and $BG$ are assigned simplicial topological
groups $AG_.$ and $BG_.$ with face homomorphisms $\partial _j:
AG_n\to AG_{n-1}$ given by:
\par $(1)$ $\partial _j (h_0[h_1|...|h_n]) = h_0h_1[h_2|...|h_n]$
for $j=0$, \par $\partial _j (h_0[h_1|...|h_n]) =
h_0[h_1|...|h_jh_{j+1}|...|h_n]$ for $0<j<n$, \par $\partial _j
(h_0[h_1|...|h_n]) = h_0[h_1|...|h_{n-1}]$ for $j=n$. While
$\partial _j: BG_n\to BG_{n-1}$ has the form:
\par $(2)$ $\partial _j ([h_1|...|h_n]) = [h_2|...|h_n]$
for $j=0$, \par $\partial _j ([h_1|...|h_n]) =
[h_1|...|h_jh_{j+1}|...|h_n]$ for $0<j<n$, \par $\partial _j
([h_1|...|h_n]) = [h_1|...|h_{n-1}]$ for $j=n$. \par The degeneracy
homomorphisms $s_j: AG_n\to AG_{n+1}$ are prescribed by the formula:
\par $(3)$  $s_j(h_0[h_1|...|h_n]) = h_0[e|h_1|...|h_n]$ for $j=0$,
\par $s_j(h_0[h_1|...|h_n]) = h_0[h_1|...|h_j|e|h_{j+1}|...|h_n]$ for
$0<j<n$, \par $s_j(h_0[h_1|...|h_n]) = h_0[h_1|...|h_n|e]$ for
$j=n$. While $s_j: BG_n\to BG_{n+1}$ is given by:
\par $(4)$  $s_j([h_1|...|h_n]) = [e|h_1|...|h_n]$ for $j=0$,
\par $s_j([h_1|...|h_n]) = [h_1|...|h_j|e|h_{j+1}|...|h_n]$ for
$0<j<n$, \par $s_j([h_1|...|h_n]) = [h_1|...|h_n|e]$ for $j=n$.
\par Analogous mappings are for simplices:
\par $(5)$ $\partial ^j (t_0,...,t_{n+1}) =
(t_0,...,t_j,t_j,t_{j+1},...,t_{n+1})$ and
\par $(6)$ $s^j(t_0,....,t_{n+1}) = (t_0,...,t_j,{\hat
t}_{j+1},t_{j+2},...,t_{n+1})$, where ${\hat t}_{j+1}$ means that
$t_{j+1}$ is absent.
\par The geometric realization $|AG_.|$ of the simplicial space
$AG_.$ is defined to be the quotient space of the disjoint union
$\sqcup_{n=0}^{\infty } \Delta ^n\times G^{n+1}$ by the equivalence
relations \par $(7)$  $(\partial ^jx,{\bar g}) \sim (x,\partial
_j{\bar g})$ for each $(x,{\bar g}) \in \Delta ^{n-1}\times
G^{n+1}$, while $(s^jx,{\bar g}) \sim (x,s_j{\bar g})$ for each
$(x,{\bar g})\in \Delta ^{n+1}\times G^{n+1}$. At the same time the
geometric realization $|BG_.|$ of the simplicial space $BG_.$ is the
quotient of the disjoint union $\sqcup_{n=0}^{\infty } \Delta
^n\times G^n$ by the equivalence relations \par $(8)$ $(\partial
^jx,{\bar g}) \sim (x,\partial _j{\bar g})$ for each $(x,{\bar g})
\in \Delta ^{n-1}\times G^n$, while $(s^jx,{\bar g}) \sim
(x,s_j{\bar g})$ for each $(x,{\bar g})\in \Delta ^{n+1}\times G^n$.
\par Consider a non-commutative sphere ${\cal
C}_r := \{ z\in {\cal I}_r: |z|=1 \} $ for $r=2, 3$, where ${\cal
I}_r := \{ z\in {\cal A}_r: Re (z)=0 \} $. For $r=1$ put ${\cal C}_r
= \{ i, -i \} $, where $i=(-1)^{1/2}$. Let ${\bf Z}({\cal C}_r)$
denotes the additive group ${\bf Z}^{{\cal C}_r}/{\cal Z}$, where
${\bf Z}^{{\cal C}_r} := \prod_{b\in {\cal C}_r} T_b$, $T_b ={\bf
Z}b$ for each $b\in {\cal C}_r$, $\bf Z$ is the additive group of
integers, $\cal Z$ is the equivalence relation such that $T_b\times
T_{-b}/{\cal Z} =T_b$ for each $b\in {\cal C}_r$. For $2\le r\le 3$
the group ${\bf Z}({\cal C}_r)$ is isomorphic with ${\bf Z}^{\alpha
}$, where $card (\alpha )= card ({\bf R})=: \sf c$. Particularly,
${\bf Z}({\cal C}_1)= {\bf Z}i$ for $r=1$.
\par Henceforth, we consider twisted sheaves and cohomologies over
$ \{ i_0,...,i_{2^r-1} \} $, where $2\le r\le 3$. In particular, the
complex case will also be included for $r=1$, but the latter case is
commutative over $\bf C$. So we can consider simultaneously $1\le
r\le 3$ and generally speak about twisting undermining that for
$r=1$ it is degenerate.

\par {\bf 19. Proposition.} {\it  Let $G$ be the group either
${\cal A}_r^*$ or ${\bf Z}({\cal C}_r)$, where $1\le r\le 3$. Then
for each $H^{\infty }$ smooth manifold $N$ over ${\cal A}_r$ and
each $b\ge 2$ the group ${\sf H}^b(N,{\bf Z}({\cal C}_r))$ is
isomorphic with:
\par $(1)$ the group ${\sf E}(N,B^{b-2}G)$ of isomorphism classes of
smooth principal $B^{b-2}G$-bundles over $N$; \par $(2)$ the group
$[N,B^{b-1}G]^{\infty }$ of smooth homotopy classes of smooth
mappings from $N$ to $B^{b-1}G$.}
\par {\bf Proof.} In view of Corollary 3.4 \cite{luoyst,luoyst2}
there exists the short exact sequence \par $(1)$ $0\to {\bf Z}({\cal
C}_r)\stackrel{\eta }\longrightarrow
{\cal A}_r\to {\cal A}_r^*\to 1$, \\
since $\exp (M + 2\pi kM/|M|)=\exp (M)$ for each non-zero purely
imaginary $M\in {\cal I}_r$ (with $Re(M)=0$) and every $k\in \bf Z$,
$1\le r\le 3$, where $\eta (z)=2\pi z$ for each $z\in {\cal A}_r$.
If $f: {\cal A}_r\to {\cal A}_r^*$ is a differentiable function,
then $(dLn f).h = w(h)$ is the differential one-form considering $d$
as the external differentiation over $\bf R$, where $h\in {\cal
A}_r$. In the particular case of $G={\cal A}_r^*$ with $1\le r \le
3$ there exist further short exact sequences
\par $(2)$ $1\to {\cal A}_r^*\to A{\cal A}_r^*\to B{\cal A}_r^*\to 1$
\par $(3)$ $1\to B{\cal A}_r^*\to AB{\cal A}_r^*\to B^2{\cal A}_r^*\to 1$
\par $(4)$ $1\to B^m{\cal A}_r^*\to AB^m{\cal A}_r^*\to
B^{m+1}{\cal A}_r^*\to 1$.

\par Therefore, identifying the ends of these short exact sequences
we get the long exact sequence
\par $(5)$ $0\to {\bf Z}({\cal C}_r)\to {\cal A}_r\to A{\cal A}_r^*\to
AB{\cal A}_r^*\to ... \to AB^m{\cal A}_r^*\to ... $,
\\ where $\sigma : {\cal A}_r\to A{\cal A}_r^*$,...,
$\sigma : AB^{m-1}{\cal A}_r^*\to AB^m{\cal A}_r^*$ are
homomorphisms, all terms ${\cal A}_r$, $A{\cal
A}_r^*$,...,$AB^m{\cal A}_r^*$,... are contractible spaces.
\par Suppose now that $N$ and $E$ are of class $H^{\infty }$.
Let ${\sf C}^{\infty }(N,AB^m{\cal A}_r^*)$ denotes the sheaf of
germs of $C^{\infty }$ functions from $N$ into $AB^m{\cal A}_r^*$.
Thus, we get the functor $C^{\infty }$. Then the application of
$C^{\infty }$ functor to the long exact sequence $(5)$ gives:
\par $(6)$ $0\to {\bf Z}({\cal C}_r)_N\to
{\sf C}^{\infty }(N,{\cal A}_r)\to {\sf C}^{\infty }(N,A{\cal
A}_r^*)\to {\sf C}^{\infty }(N,AB{\cal A}_r^*)\to ... \to {\sf
C}^{\infty }(N,AB^m{\cal A}_r^*)\to ...$,
\\ where $\sigma _*: {\sf C}^{\infty }(N,{\cal A}_r)\to
{\sf C}^{\infty }(N,A{\cal A}_r^*)$,...,$\sigma _*: {\sf C}^{\infty
}(N,AB^{m-1}{\cal A}_r^*)\to {\sf C}^{\infty }(N,AB^m{\cal A}_r^*)$
are induced homomorphisms.
\par The latter exact sequence is called the bar resolution of
${\bf Z}({\cal C}_r)_N$. Sheaves ${\sf C}^{\infty }(N,{\cal A}_r)$
and ${\sf C}^{\infty }(N,AB^m{\cal A}_r^*)$ are contractible, since
${\cal A}_r$ and $AB^m{\cal A}_r^*$ are contractible. Therefore, the
cohomology of the sheaf ${\bf Z}({\cal C}_r)_N$ can be computed
using the complex
\par $(7)$ $C^{\infty }(N,{\cal A}_r)\to C^{\infty }(N,A{\cal A}_r^*)\to
... \to C^{\infty }(N,AB^m{\cal A}_r^*)\to ... $ \\
with homomorphisms \\ $\sigma _*: C^{\infty }(N,{\cal A}_r)\to
C^{\infty }(N,A{\cal A}_r^*)$,...,$\sigma _*: C^{\infty
}(N,AB^{m-1}{\cal A}_r^*)\to C^{\infty }(N,AB^m{\cal A}_r^*)$. \\
The long exact sequence $(7)$ we call a bar cochain complex of ${\bf
Z}({\cal C}_r)_N$. The cohomology of ${\bf Z}({\cal C}_r)_N$
computed with the help of the bar complex is denoted by ${\sf
H}_b^*(N,{\bf Z}({\cal A}_r)_N)$ and it is called the bar cohomology
of ${\bf Z}({\cal C}_r)_N$. Then $\pi _0C^{\infty }(N,B{\cal
A}_r^*)$ is the first bar cohomology ${\sf H}_b^1(N,{\bf Z}({\cal
A}_r)_N)$ of ${\bf Z}({\cal C}_r)_N$.
\par For the generalized exponential sequence
\par $0\to {\bf Z}({\cal C}_r)_N \to AB^{<b-2}G_N\to B^{b-2}G_N[2-b]\to e$
there exists the cohomology long exact sequence
\par $... \to \mbox{ }_h{\sf H}^{b-1}(N,AB^{< b-2}G_N)\to {\sf
H}^1(N,B^{b-2}G_N)$ \par $\to {\sf H}^b(N;{\bf Z}({\cal C}_r))\to
\mbox{ }_h{\sf H}^b(N,AB^{< b-2}G_N)\to ...$, \\
where $AB^{<b}G_N$ is the complex
\par $K_N\stackrel{\sigma}\longrightarrow AG_N
\stackrel{\sigma}\longrightarrow ABG_N\stackrel{\sigma
}\longrightarrow AB^2G_N \stackrel{\sigma}\longrightarrow ...
\stackrel{\sigma }\longrightarrow  AB^{b-1}G_N$, \\
$K$ is equal to either ${\cal A}_r$ or $A{\bf Z}({\cal C}_r)$ for
$G={\cal A}_r^*$ or $G= B{\bf Z}({\cal C}_r)$ respectively, where
$\mbox{ }_h{\sf H}^b(N,{\cal B})$ denotes a hypercohomology on $N$
with coefficients in a complex of sheaves $\cal B$ (see its
definition in \cite{bredon,gajer,gajinvm} and \S 10 above).
\par We have that for each $b\ge 0$ the sheaf $AB^bG_N$ is acyclic.
Consider the groups of global sections $AB^bG_N(N)$ of the sheaves
$AB^bG_N$. Therefore, the cohomology of the complex $AB^{<b-2}G_N$
is equal to the cohomology of the cochain complex
\par $K\stackrel{\sigma}\longrightarrow AG_N(N)
\stackrel{\sigma}\longrightarrow
ABG_N(N)\stackrel{\sigma}\longrightarrow AB^2G_N(N)
\stackrel{\sigma}\longrightarrow ...
\stackrel{\sigma}\longrightarrow  AB^{b-1}G_N(N)$. \par Thus,
$\mbox{ }_h{\sf H}^m(N, AB^{<b-2}G_N)\cong {\sf H}^m(N, AB^{<
b-2}G_N(N)) \cong e$ for each $m>b-2$, consequently, the coboundary
homomorphism ${\sf H}^1(N, B^{b-2}G_N)\to {\sf H}^b(N;{\bf Z}({\cal
C}_r))$ is an isomorphism (see also Chapter 2 \S 4 in \cite{bredon}
for abelian sheafs).
\par The second statement of this proposition follows from Lemma
17.

\par {\bf 20. Lemma.} {\it Let $X$ be a topological vector space
over ${\cal A}_r$, $2\le r\le 3$. Then $AX$ and $BX$ with respect to
the additive group structure of $X$ and with respect to the
multiplication on scalars from ${\cal A}_r$ in homogeneous
coordinates are ${\cal A}_r$ vector spaces and the projection $AX\to
BX$ is $\bf R$-homogeoneous and ${\cal A}_r$ additive.}

\par {\bf Proof.} Define the multiplications by:
for ${\cal A}_r\times AX\to AX$ as \par $s|t_1,...,t_n;
v_0[v_1|...|v_n]| = |t_1,...,t_n; sv_0[sv_1|...|sv_n]|$, \\ for
$AX\times {\cal A}_r\to AX$ as \par $|t_1,...,t_n;
v_0[v_1|...|v_n]|s = |t_1,...,t_n; v_0s[v_1s|...|v_ns]|$,\\ for
${\cal A}_r\times BX\to BX$ as \par $s|t_1,...,t_n; [v_1|...|v_n]| =
|t_1,...,t_n; [sv_0|...|sx_n]|$,\\ for $BX\times {\cal A}_r\to BX$
as \par $|t_1,...,t_n; [v_1|...|v_n]|s = |t_1,...,t_n;
[v_1s|...|v_ns]|$. \\ Then if $q_j =s_m(s_{m-1}...(s_1v_j)...)$ for
each $j$, then for $z\ne 0$ we get $zq_j =
z(s_m...(s_1(z^{-1}(zv_j)))...)$ due to the alternativity of the
octonion algebra $\bf O$, while for $z=0$ we trivially get $0
=(s_m...(s_10)...)$. Thus such multiplication is compatible with the
equivalence relations, since $X= X_0i_0\oplus ... \oplus
X_{2^r-1}i_{2^r-1}$, where $X_0,...,X_{2^r-1}$ are pairwise
isomorphic topological vector spaces over $\bf R$ such that we put
$vx=xv$ for each $v\in X_j$ and $x\in \bf R$.
\par Since $\bf R$ is the center of the algebra $\bf O$, then
the projection from $AX$ to $BX$ is $\bf R$-linear. Evidently, it is
additive as the additive group homomorphism.

\par {\bf 21. Remark.} Let $N_.$ be a simplicial smooth manifold
over ${\cal A}_r$, where $0\le r\le 3$. A smooth $m$-form $w$ on the
geometric realization $|N_.|$ of $N_.$ is defined to be as a family
$ \{ w^k: k \} $ of smooth differential $m$-forms $w^k$ on $\Delta
^k\times N_k$ with values in ${\cal A}_r$ being applied to vectors,
satisfying for each $0\le j\le n$ the compatibility conditions:
\par $(1)$ $(\partial ^j\times id)^* w^n = (id\times \partial _j)^*w^{n-1}$
\par $(2)$ $(s^j\times id)^*w^n = (id\times s_j)^*w^{n+1}$, \\
where $\partial ^j\times id$, $id\times \partial _j$, $s^j\times id$
and $id\times s_j$ are the maps as follows: \par $(3)$ $id \times
\partial _j: \Delta ^{n-1}\times N_n\to \Delta ^{n-1}\times
N_{n-1}$, $\partial ^j\times id: \Delta ^{n-1}\times N_n\to \Delta
^n\times N_n$, $id\times s_j: \Delta ^{n+1}\times N_n\to \Delta
^{n+1}\times N_{n+1}$, $s^j\times id : \Delta ^{n+1}\times N_n\to
\Delta ^n\times N_n$ such that $\partial ^j$ and $s^j$ are coface
and the codegeneracy maps on $\Delta ^n$ and $\partial _j$, $s_j$
are the face and the degeneracy maps on $N_n$. We consider $w$
taking values in a vector space or an algebra over ${\cal A}_r$ as
is specified below.
\par For a Lie group $G$ either over $\bf R$ or may be twisted over ${\cal A}_r$
and its Lie algebra $\sf g$ put $g^{-1}dg$ as the canonical $\sf
g$-valued connection $1$-form on $G$ (see also Lemma 20). Under the
mapping $g\mapsto hg$ and $dg\mapsto hdg$ we have
$(g^{-1}h^{-1})(hdg) = g^{-1}dg$ due to the alternativity of $G$ and
the Maufang identity $(xy)(zx) = x(yz)x$ for each $x, y, z\in \bf O$
and $de = d(g^{-1}g) = 0 = (dg^{-1})g + g^{-1}dg =
[(dg^{-1})h^{-1}](hg) + (g^{-1}h^{-1})(hdg)$. Iterating this
relation due to the alternativity of $\bf O$ and the Moufang
identities in it we get the equivariance condition in homogeneous
coordinates over $\bf O$ as well:
$[(...(g^{-1}s_1^{-1})...s_{m-1}^{-1})s_m^{-1}]
[s_m(s_{m-1}...(s_1dg_1)...)] = g_1^{-1}dg_1$. \par The total space
$AG$ of the universal principal $g$-bundle $AG\to BG$ carries a
smooth $\sf g$-valued form $w$. The evaluation of $w$ is $w
|x_0,...,x_n,g_0,...,g_n| = x_0g_0^{-1}dg_0 + x_1g_1^{-1}dg_1 + ...
+ x_n g_n^{-1}dg_n$, where $x_0,...,x_n$ are barycentric coordinates
in $\Delta ^n$. Each term $x_ng_n^{-1}dg_n.s$ is in $\sf g$ for each
$s\in \sf g$ such that $g_j^{-1}dg_j = \pi
_j^*(g^{-1}dg|_{T_{g_j}G})$, where $\pi _j: G^{n+1}\to G$ is the
projection on the $j$-th factor and $g^{-1}dg|_{T_{g_j}G}$ is the
restriction of $g^{-1}dg$ to the tangent space $T_{g_j}G$ of $G$ at
$g_j$.
\par For $A{\cal A}_r^*$ with $2\le r\le 3$ define the canonical
connection $1$-form $A(z^{-1}dz)$ by the family of $A{\cal
A}_r$-valued $1$-forms $A(z^{-1}dz)^n$ on $\Delta ^n\times ({\cal
A}_r^*)^{n+1}$ such that $A(z^{-1}dz)^n$ evaluated on a vector
$(v_0,...,v_n)$ at a point $|t_1,...,t_n,z_0[z_1|...|z_n]|$ is given
by the formula
\par $(4)$ $(A(z^{-1}dz)^n|_{|t_1,...,t_n,z_0[z_1|...|z_n]|}.(v_0,...,v_n)
=|t_1,...,t_n,z_0^{-1}v_0[z_1^{-1}v_1|...|z_n^{-1}v_n]|$ \\ and
formally denote it by
\par $(5)$ $(A(z^{-1}dz)^n|_{|t_1,...,t_n,z_0[z_1|...|z_n]|}
=|t_1,...,t_n,z_0^{-1}dz_0[z_1^{-1}dz_1|...|z_n^{-1}dz_n]|$.
\par For $B{\cal A}_r^*$ with $2\le r\le 3$ the canonical connection
$1$-form $B(z^{-1}dz)$ on $B{\cal A}_r^*$ is defined by the family
of $B{\cal A}_r$-valued $1$-forms $B(z^{-1}dz)^n$ on $\Delta
^n\times ({\cal A}_r^*)^n$, where
\par $(6)$ $B(z^{-1}dz)^n|_{|t_1,...,t_n,[z_1|...|z_n]|}
= |t_1,...,t_n,[z_1^{-1}dz_1|...|z_n^{-1}dz_n]|$.
\par We have that \par $(\partial ^j\times id)^*A(z^{-1}dz)^n
|_{|t_1,...,t_{n-1}; z_0[z_1|...|z_n]|} =$ \par
$|t_1,...,t_j,t_j,t_{j+1},...,t_{n-1}; z_0^{-1}dz_0
[z_1^{-1}dz_1|... |z_n^{-1}dz_n]|$ and
\par $(\partial ^j\times id)^*B(z^{-1}dz)^n
|_{|t_1,...,t_{n-1},[z_1|...|z_n]|} =$ \par
$|t_1,...,t_j,t_j,t_{j+1},...,t_{n-1},[z_1^{-1}dz_1|...
|z_n^{-1}dz_n]|$ and \par $(id\times \partial _j)^*A(z^{-1}dz)^{n-1}
|_{|t_1,...,t_{n-1}; z_0[z_1|...|z_n]|} =$ \\
$|t_1,...,t_{n-1}; (z_0^{-1}dz_0) + (z_1^{-1}dz_1)[ z_2^{-1}dz_2|
... |z_n^{-1}dz_n]|$ for $j=0$,
\par $(id\times \partial _j)^*A(z^{-1}dz)^{n-1}
|_{|t_1,...,t_{n-1}; z_0[z_1|...|z_n]|} =$ \\
$|t_1,...,t_{n-1}; z_0^{-1}dz_0 [z_1^{-1}dz_1|...
|z_{j-1}^{-1}dz_{j-1} |(z_j^{-1}dz_j) + (z_{j+1}^{-1}dz_{j+1})|
z_{j+2}^{-1}dz_{j+2}|... |z_n^{-1}dz_n]|$ \\
for $0<j<n$, \par $(id\times \partial _j)^*A(z^{-1}dz)^{n-1}
|_{|t_1,...,t_{n-1}; z_0[z_1|...|z_n]|} =$ \\
$|t_1,...,t_{n-1}; z_0^{-1}dz_0 [z_1^{-1}dz_1|...
|z_{n-1}^{-1}dz_{n-1}]|$ for $j=n$, while
\par $(id\times \partial _j)^*B(z^{-1}dz)^{n-1}
|_{|t_1,...,t_{n-1}; [z_1|...|z_n]|} =$ \\
$|t_1,...,t_{n-1}; [ z_2^{-1}dz_2| ... |z_n^{-1}dz_n]|$ for $j=0$,
\par and $(id\times
\partial _j)^*B(z^{-1}dz)^{n-1} |_{|t_1,...,t_{n-1}; [z_1|...|z_n]|}
=$ \\ $|t_1,...,t_{n-1}; [z_1^{-1}dz_1|... |z_{j-1}^{-1}dz_{j-1}
|(z_j^{-1}dz_j)+ (z_{j+1}^{-1}dz_{j+1})| z_{j+2}^{-1}dz_{j+2}|...
|z_n^{-1}dz_n]|$ for $0<j<n$ \\
\par $(id\times \partial _j)^*B(z^{-1}dz)^{n-1}
|_{|t_1,...,t_{n-1}; [z_1|...|z_n]|} =$ \\
$|t_1,...,t_{n-1}; [z_1^{-1}dz_1|... |z_{n-1}^{-1}dz_{n-1}]|$ for
$j=n$ \\
and using the equivalence relations 4$(2,3,5)$, we get the
compatibility Condition $(1)$ for such differential forms, since
$\partial ^j$ corresponds to inserting $x_j=0$ and the latter
corresponds to $t_j=t_{j+1}$, because $t_j=x_0+...+x_{j-1}$, while
$h_j=e$ corresponds to $g_j=g_{j+1}$.
\par Further we get:
\par $(s^j\times id)^*A(z^{-1}dz)^n
|_{|t_1,...,t_{n+1}; z_0[z_1|...|z_n]|} =$ \\
$|t_1,...,t_j, t_{j+2},..., t_{n+1}; z_0^{-1}dz_0,[z_1^{-1}dz_1|...
|z_n^{-1}dz_n]|$ for each $j$ and
\par $(s^j\times id)^*B(z^{-1}dz)^n
|_{|t_1,...,t_{n+1}; [z_1|...|z_n]|} =$ \par
$|t_1,...,t_j,t_{j+2},...,t_{n+1}; [z_1^{-1}dz_1|...
|z_n^{-1}dz_n]|$, \\
since $s^j(x_0,...,x_{n+1}) =
(x_0,...,x_{j-1},x_j+x_{j+1},x_{j+2},...,x_{n+1})$ and $x_j=0$
corresponds to $t_j=t_{j+1}$. Then
\par $(id\times s_j)^*A(z^{-1}dz)^{n+1}
|_{|t_1,...,t_{n+1}; z_0[z_1|...|z_n]|} =$ \\
$|t_1,...,t_{n+1}; z_0^{-1}dz_0[0|z_1^{-1}dz_1|...|z_n^{-1}dz_n]|$
and
\par $(id\times s_j)^*A(z^{-1}dz)^{n+1}
|_{|t_1,...,t_{n+1}; z_0[z_1|...|z_n]|} =$ \\
$|t_1,...,t_{n+1}; z_0^{-1}dz_0,[z_1^{-1}dz_1|... |z_j^{-1}dz_j|0
|z_{j+1}^{-1}dz_{j+1}|...|z_n^{-1}dz_n]|$ for $0<j<n$ and
\par $(id\times s_j)^*A(z^{-1}dz)^{n+1}
|_{|t_1,...,t_{n+1}; z_0[z_1|...|z_n]|} =$ \\
$|t_1,...,t_{n+1}; z_0^{-1}dz_0,[z_1^{-1}dz_1|... |z_n^{-1}dz_n|0]|$
for $j=n$ and shortly write
\par $(id\times s_j)^* B(z^{-1}dz)^{n+1}
|_{|t_1,...,t_{n+1}; [z_1|...|z_n]|} =$ \\
$|t_1,...,t_{n+1}; [z_1^{-1}dz_1|... |z_j^{-1}dz_j|
0|z_{j+1}^{-1}dz_{j+1}|...|z_n^{-1}dz_n]|$ for each $j$. \\
Using the equivalence relations 4$(2,3,5)$ in $AB^bG$ and $B^{b+1}G$
we get the compatibility Condition $(2)$.
\par For a twisted group $G$ satisfying Conditions 4$(A1,A2,C1,C2)$
a smooth $k$-form on $AB^bG$ and $B^{b+1}G$ is defined by induction.
Let smooth differential $k$-forms on $B^bG$ and each $\Delta
^k\times (B^bG)^m$ for $k, m\ge 0$ be defined. Then a smooth
$k$-form $w$ on $AB^{b+1}G$ consists of a family of $k$-forms $w^n$
on $\Delta ^n\times (B^bG)^{n+1}$ satisfying the compatibility
conditions $(1,2)$. For $B^{b+1}G$ a smooth $k$-form $w$ on
$B^{b+1}G$ consists of a family of $k$-forms $w^n$ on $\Delta
^n\times (B^bG)^{n+1}$ satisfying the compatibility Conditions
$(1,2)$. Then a smooth $k$-form $w$ on $\Delta ^k\times
(B^{b+1}G)^m$ consists of a family of $k$-forms $w^n$ on $\Delta
^k\times (\Delta ^n\times (B^bG)^{n+1})^m$ satisfying the
compatibility conditions:
\par $(7)$ $id_{\Delta ^k}\times (\partial ^j\times id)^m)^* w^n =
(id_{\Delta ^k}\times (id\times \partial _j)^m)^* w^{n-1}$,
\par $(8)$ $(id_{\Delta ^k}\times (s^j\times id)^m)^* w^n =
(id_{\Delta ^k} \times (id\times s_j)^m)^* w^{n+1}$.
\par It was shown above that the groups $AB^b{\cal A}_r$ and
$B^{b+1}{\cal A}_r$ also have the structures of ${\cal A}_r$ vector
spaces. Therefore, the canonical connection $1$-form
$AB^b(z^{-1}dz)$ on $AB^b{\cal A}_r^*$ is a $1$-form on $AB^b{\cal
A}_r^*$ such that it satisfies the inductive formula:
\par $(9)$ $AB^b(z^{-1}dz)|_{|t_1,...,t_n,g_0[g_1|...|g_n]|} =
|t_1,...,t_n,B^b(g_0^{-1}g_0)
[B^b(g_1^{-1}dg_1)|...|B^b(g_n^{-1}dg_n)]|$.
\par Then the canonical connection $1$-form $B^{b+1}(z^{-1}dz)$ on
$B ^{b+1}{\cal A}_r^*$ is a $1$-form on $B^{b+1}{\cal A}_r^*$ such
that
\par $(10)$ $B^{b+1}(z^{-1}dz)|_{|t_1,...,t_n,[g_1|...|g_n]|} =
|t_1,...,t_n,[B^b(g_1^{-1}dg_1)|...|B^b(g_n^{-1}dg_n)]|$, \\
where $g_0, g_1,...,g_n\in B^b{\cal A}_r^*$ and $B^b(g_j^{-1}dg_j)$
is the canonical connection $1$-form $B^b(z^{-1}dz)$ on $B^b{\cal
A}_r^*$ evaluated at $g_j$.

\par {\bf 22. Gerbes over quaternions and octonions.}
Consider twisted groups $C, K, G$ satisfying Conditions
4$(A1,A2,C1,C2)$. If \par $(CE1)$ $e\to C_0\to K_0\to G_0\to e$ is a
topological central extension, then we say, that
\par $(CE2)$ $e\to C\to K\to G\to e$ is a topological twisted extension.
\par A gerbe on a topological space $X$ is a sheaf $\cal S$ of
categories satisfying the conditions $(G1-G3)$:
\par $(G1)$ for each open subset $V$ in $X$ the category ${\cal
S}(V)$ is a groupoid, which means that every morphism is invertible;
\par $(G2)$ each point $x\in X$ has a neighborhood $V_x$ for which
${\cal S}(V_x)$ is non-empty;
\par $(G3)$ any two objects $P_1$ and $P_2$ of ${\cal S}(V)$ are
locally isomorphic, that is, each $x\in V$ has a neighborhood $Y$
for which the restrictions $P_1|_Y$ and $P_2|_Y$ are isomorphic.
\par A gerbe $\cal S$ is called bound by a sheaf $\cal G$ of twisted
groups over ${\cal A}_r$ satisfying Conditions 5$(A1,C1,C2,7)$, if
for each open subset $V$ in $X$ and every object $P$ of ${\cal
S}(V)$ there exists an isomorphism of sheaves $\nu : Aut (P)\to
{\cal G}|_V$, where ${\cal G}|_V$ denotes the restriction of the
sheaf $\cal G$ onto $V$, while $Aut (P)$ is the sheaf of
automorphisms of $P$ so that for an open subset $Y$ in $V$ the group
$Aut (P)(Y)$ is the group of automorphisms of the restriction
$s_Y(P)$. It is supposed that such an isomorphism commutes with with
morphisms of $\cal S$ and must be compatible with restrictions to
smaller open subsets.
\par Two gerbes $\cal S$ and $\cal E$ bounded by $\cal G$ on a
manifold $N$ are equivalent, if they satisfy $(G4,G5)$:
\par $(G4)$ if $V$ is an open subset in $X$, then there exists
an equivalence of categories $\mu _V: {\cal S}(V)\to {\cal E}(V)$ so
that for each object $P$ of ${\cal S}(V)$ there is a commutative
diagram: \par $\mu _V: Aut _{{\cal S}(V)}(P)\to Aut _{{\cal
E}(V)}(P)$, \par $\nu _{\cal S}: Aut _{{\cal S}(V)}(P)\to \Gamma
(V,{\cal G})$,
\par $\nu _{\cal E}: Aut _{{\cal E}(V)}(P)\to \Gamma (V,{\cal G})$
such that \par $\nu _{\cal S} = \nu _{\cal E}(\mu _V)$;
\par $(G5)$ for each pair of open subsets $V$ and $Y$ in $N$ with
$Y\subset V$ there exists an invertible natural transformation:
$\beta : R_{\cal E}(\mu _V) = \mu _Y(R_{\cal S})$, where $R_{\cal
S}: {\cal S}(V)\to {\cal S}(Y)$ denotes the natural restriction
transformation. It is also imposed the condition, that for a triple
of open subsets $Y\subset V\subset J$ in $N$ the compatibility
conditions are satisfied.
\par If there is a principal $G$-bundle $E(B,G,\pi ,\Psi )$ and
and an extension $(CE1, CE2)$ of topological groups, then there
exists a gerbe ${\cal G}_{\pi }$ bound by $C_N$ on $B$. This gerbe
is constructed from the sheaf of sections of the bundle $E(B,G,\pi
,\Psi )$ by posing for each open subset $V$ of $B$ objects and
morphisms of ${\cal G}_{\pi }(V)$ as follows. Associate with each
section $s: V\to \pi ^{-1}(V)$ of $\pi : \pi ^{-1}(V)\to V$ the
$G$-equivariant map $t_s: \pi ^{-1}(V)\to G$ such that $t_s(z) s(\pi
(z)) =z$ for every $z\in \pi ^{-1}(V)$. We have as well the
pull-back of principal $C$-bundle $K\to G$ from $G$ to $\pi
^{-1}(V)$ due to the mapping $t_s: \pi ^{-1}(V)\to G$. \par The
composition $\pi \circ \pi _s: E_s\to V$ is a principal $K$-bundle
having a lifting of the structure group of $\pi ^{-1}(V)\to V$ to
$K$. Then pairs $(E,f)$ of principal $K$-bundles $\pi _V: E(V,K,\pi
,\Psi )\to V$ and principal $C$-bundles $f: E\to \pi ^{-1}(V)$ such
that there exists the commutative diagram with $\pi (f(*))=\pi
_V(*)$.
\par A morphism of principal $K$-bundles $\eta : E\to E_1$ from
$(E,f)$ to $(E_1,f_1)$ is described with the help of the condition
$f_1(\eta (*)) = f$ with the corresponding commutative diagram.
Therefore, the group of automorphisms of every object $(E,f)$ of
${\cal G}_{\pi }(V)$ is the group of mappings from $V$ to $C$ being
the section of the sheaf $C_N$ over $V$, consequently, ${\cal
G}_{\pi }$ is the gerbe bound by $C_N$.
\par The constructed above gerbe ${\cal G}_{\pi }$ has a global
section if and only if there exists a lifting of the structure group
from $C$ to $K$. For $G=B{\cal A}_r^*$ the extension $(CE1,CE2)$
takes the form
\par $1\to {\cal A}_r^*\to A{\cal A}_r^*\to B{\cal A}_r^*\to 1$.
Then each principal $A{\cal A}_r^*$-bundle is trivial, since $A{\cal
A}_r^*$ is contractible. Hence the gerbe ${\cal G}_{\pi }$ has a
global section if and only if $E(N,B{\cal A}_r^*,\pi ,\Psi )$ is a
trivial $B{\cal A}_r^*$-bundle.
\par Construct now another gerbe ${\cal L}_{\pi }$ of local sections
of the bundle $E(B,B{\cal A}_r^*,\pi ,\Psi )$. For each open subset
$V$ in $B$ the objects of ${\cal L}_{\pi }(V)$ are sections of $E$
over $V$ so that each local section $s: V\to \pi ^{-1}(V)$ induces a
$B{\cal A}_r^*$-equivariant mapping $t_s: \pi ^{-1}(V)\to B{\cal
A}_r^*$ that induces the mapping $\tau _s=t_s(s(*)): V\to B{\cal
A}_r^*$. \par If ${\sf E}_s$ is a principal ${\cal A}_r^*$-bundle
over $V$ induced by the mapping $\tau _s$, then a morphism between
the objects $s, s_1\in {\cal L}_{\pi }(V)$ induces the morphism
${\sf E}_s\to {\sf E}_{s_1}$ of the corresponding principal ${\cal
A}_r^*$-bundles. Then ${\cal L}_{\pi }$ is a gerbe bounded by
$({\cal A}_r^*)_N$, where $2\le r\le 3$. Therefore, the natural
transformation ${\cal L}_{\pi }(V)\to {\cal G}_{\pi }(V)$ sending a
section $s$ to the pull-back ${\sf E}_s$ of the universal principal
${\cal A}_r^*$-bundle by $t_s$ is an equivalence of categories,
which extends to an equivalence of gerbes ${\cal L}_{\pi }\to {\cal
G}_{\pi }$.
\par For a gerbe ${\cal G}$ on $N$ bounded by $({\cal A}_r^*)_N$ with
$2\le r\le 3$, assigning to each object $Q$ in ${\cal G}(V)$ an
${\cal S}^1_{N,{\cal A}_r}$-torsor ${\cal C}_{{\cal O}Q}$ on $V$
induces a connective structure. This torsor ${\cal C}_{{\cal O}Q}$
consists of a sheaf on which ${\cal S}^1_{N,r}$ acts so that for
each point $x\in N$ there exists a neighborhood $V$ having the
property that for each open subset $Y\subset V$ the group ${\cal
C}_{{\cal O}Q}(Y)$ is a principal homogeneous space under the group
$\Gamma (Y,{\cal S}^1_{N,{\cal A}_r})$. This assignment $Q\mapsto
{\cal C}_{{\cal O}Q}(V)$ need to be functorial in accordance with
restrictions from $V$ onto $Y$. Moreover, for each morphism $\phi :
Q\to J$ of objects of ${\cal G}(V)$ there exists an isomorphism
$\phi _*: {\cal C}_{{\cal O}Q}(V)\to {\cal C}_{{\cal O}J}(V)$ of
${\cal S}^1_{N,{\cal A}_r}$-torsors. Since ${\cal G}$ is a gerbe,
then $\phi $ is an isomorphism and $\phi _*$ is compatible with
compositions of morphisms and with restrictions to smaller open
subsets, $Y\subset V$. If $\phi $ is an automorphism of $Q$ induced
by an ${\cal A}_r^*$-valued function $g$ we suppose that $\phi _*$
is an automorphism $\nabla \mapsto \nabla - dLn (g)$ of the ${\cal
S}^1_{N,{\cal A}_r}$-torsor ${\cal C}_{{\cal O}Q}(V)$.
\par Consider a connection $\omega $ on a smooth principal $B{\cal
A}_r^*$-bundle $E(N,B{\cal A}_r^*,\pi ,\Psi )$ and let $V$ be an
open subset in $N$ such that ${\cal G}_{\pi }(V)$ is non-void and
let $\omega _V$ be the restriction of $\omega $ to $\pi ^{-1}(V)$.
To each element $(E,f)$ of ${\cal G}_{\pi }(V)$ it is possible
assign a set ${\cal C}_{{\cal O}E}^{\omega }(V)$ of connections on
$E$ compatible with $\omega $. If $\omega (q(*)) = f^*\omega $ for
principal ${\cal A}_r^*$-bundles $q: A{\cal A}_r\to B{\cal A}_r$ and
$f: E\to \pi ^{-1}(V)$, then $\omega $ generates an element ${\hat
{\omega }}\in {\cal C}_{{\cal O}E}(V)$. Therefore, the assignment
$\omega \mapsto {\cal C}_{\cal O}^{\omega }$ is the connective
structure on ${\cal G}_{\pi }$.
\par The equivalence of gerbes ${\cal L}_{\pi }\to {\cal G}_{\pi }$
implies an extension of a pull-back of the connective structure from
${\cal G}_{\pi }$ to ${\cal L}_{\pi }$.

\par {\bf 23. Corollary.} {\it A mapping posing to the isomorphism
class of a principal $B{\cal A}_r^*$-bundle $E(B,B{\cal A}_r^*,\pi
,\Psi )$, $\pi : E = A{\cal A}_r^* \to B{\cal A}_r^*$, the
equivalence class of the gerbe of section ${\cal L}_{\pi }$ of
$E(B,B{\cal A}_r^*,\pi ,\Psi )$ induces an isomorphism between the
group of isomorphism classes of principal $B{\cal A}_r^*$-bundles
and the group of equivalence classes of gerbes bound by ${\cal
A}_r^*$.}

\par {\bf 24. Corollary.} {\it A mapping sending to the isomorphism
class of a principal $B{\cal A}_r^*$-bundle $E(B,B{\cal A}_r^*,\pi
,\Psi )$, $\pi : E = A{\cal A}_r^* \to B{\cal A}_r^*$, with a
connection $\omega $ the equivalence class of the gerbe of section
${\cal L}_{\pi }$ of $E(B,B{\cal A}_r^*,\pi ,\Psi )$ with the
connective structure on ${\cal L}_{\pi }$ induced by $\omega $
induces an isomorphism between the group of isomorphism classes of
principal $B{\cal A}_r^*$-bundles with connection on the group of
equivalence classes of gerbes bound by ${\cal A}_r^*$ with
connective structures.}

\par {\bf Proof.} This follows from Proposition 19 and \S 22,
since the case of $2\le r\le 3$ is obtained from the complex case
(see Theorems A1, A2 \cite{gajinvm}) by additional doubling
procedure of groups with doubling generators: $\bf H$ from $\bf C$
and $\bf O$ from $\bf H$, while the considered groups satisfy
Conditions 4$(A1,A2,C1,C2)$.

\par {\bf 25. Sheaves, geometric bars and gerbes for wrap groups.}
\par If $G$ satisfies Conditions 4$(A1,A2,C1,C2)$, then wrap groups
$(W^ME)_{t,H}$ satisfy these Conditions 4$(A1,A2,C1,C2)$ as well,
since $G^k$ satisfies them being a multiplicative subgroup of the
ring ${\hat G}^k$ and $(W^{M, \{ s_{0,q}: q=1,...,k \} } E;N,G,{\bf
P})_{t,H}$ is the principal $G^k$-bundle over the commutative group
$(W^{M, \{ s_{0,q}: q=1,...,k \} }N)_{t,H}$ (see Propositions
7$(1,2)$ \cite{lulaswgof}). Thus, wrap groups can be taken as the
particular cases of groups for the sheaves, geometric bar and gerbes
constructions (see \S \S 1, 4, 11-13, 22, Corollary 9, Lemmas 16
\cite{lulaswgof}, 17, etc.).

\par More concretely this can be done as follows. For a
pseudo-manifold $X=X_1\times X_2$ over ${\cal A}_r$, where $X_1$ and
$X_2$ are $H^t_p$-pseudo-manifolds over ${\cal A}_r$, suppose that
for each points $s_{0,1},...,s_{0,k}$ in $X_1$ and every
neighborhood $U$ of $\{ s_{0,1},...,s_{0,k} \} $ in $X_1$ and a
point $y_0\in X_2$ and every neighborhood $V$ of $y_0$ in $X_2$
there exist manifolds $M$ and $N$ such that $\{ s_{0,1},...,s_{0,k}
\} \subset M\subset U$ and $y_0\in N\subset V$ for which a principal
$G$-bundle $E(N,G,\pi ,\Psi )$ exists with a marked group $G$
satisfying conditions of \S 2 \cite{luwrgfbqo}. If \par $(1)$
$J(\Lambda ) = \prod_{\alpha \in \Lambda } J_{\alpha }$ \\ is the
product of topological groups $J_{\alpha }$, where $\Lambda $ is a
set, and $\Lambda _2\subset \Lambda _1$, then there exists the
natural projection group homomorphism \par $(2)$ ${\hat s}_{\Lambda
_2,\Lambda _1} : J(\Lambda _1)\to J(\Lambda _2)$. \\ Then define a
pre-sheaf $F$ on $X$ such that \par $(3)$ $F(U\times V) =
\prod_{s_{0,1},...,s_{0,k}\in M\subset U; y_0\in N\subset V} (W^{M,
\{ s_{0,q}: q=1,...,k \} } E;N,G,{\bf P})_{t,H}$ \\ and
$s_{U_2\times V_2,U_1\times V_1}: F(U_1\times V_1)\to F(U_2\times
V_2)$, since $(W^{M_2, \{ s_{0,q}: q=1,...,k \} } E;N_2,G,{\bf
P})_{t,H}\subset (W^{M_1, \{ s_{0,q}: q=1,...,k \} } E;N_1,G,{\bf
P})_{t,H}$ for $\{ s_{0,q}: q=1,...,k \}\subset M_2\subset M_1$ and
$y_0\in N_2\subset N_1$ satisfying conditions of Theorem 10
\cite{lulaswgof}, where $U_2\subset U_1$ and $V_2\subset V_1$, while
open subsets of the form $U\times V$ contain the base of topology of
$X$.
\par If $\cal S$ is a sheaf on $X$ and ${\cal S}(U)$ satisfies
Conditions 4$(A1,A2,C1,C2)$ for each $U$ open in $X$, then we call
$\cal S$ the twisted sheaf over $ \{ i_0,...,i_{2^r-1} \} $.
\par For $k=1$ consider $x = \{ s_0; y_0 \} \in X$, but
generally, consider $x= \{ s_{0,1},...,s_{0,k}; y_0 \} \in
X_1^k\times X_2$ instead of $X_1\times X_2$. Then a set ${\cal F}_x$
of all germs of the pre-sheaf $F$ at a point $x\in X_1^k\times X_2$
is the inductive limit ${\cal F}_x = ind-\lim F(U\times V)$ taken by
all open neighborhoods $U^k\times V$ of $x$ in $X_1^k\times X_2$.
Then applying the general construction of \S 1 gives the sheaf
${\cal S}_{W,X_1,X_2}$ of wrap groups. It is twisted over $\{
i_0,...,i_{2^r-1} \}$ for the group $G$ twisted over generators $\{
i_0,...,i_{2^r-1} \}$ for $2\le r\le 3$. This sheaf is commutative,
if $G$ is commutative.
\par This sheaf is obtained from the given below generalization
taking a constant sheaf of the group $G=G(U)$ for each $U$ open in
$X_1$.
\par More generally, if there is a sheaf ${\cal G}={\cal G}_{X_1}$
on $X_1$ of groups such that for each $U$ open in $X_1$ a group
$G(U)$ satisfies conditions of \S 2 in \cite{luwrgfbqo}, then put
\par $(4)$ $F(U\times V) = \prod_{s_{0,1},...,s_{0,k}\in M\subset U;
y_0\in N\subset V} (W^{M, \{ s_{0,q}: q=1,...,k \} } E;N,G(U),{\bf
P})_{t,H}$, \\ where $s_{U_2,U_1}: G(U_1)\to G(U_2)$ is the
restriction mapping for each $U_2\subset U_1$ so that the parallel
transport structure for $M\subset U$ is defined, ${\cal G}_{X_1}$,
may be twisted for $2\le r\le 3$. Therefore, due to Theorem 10
\cite{lulaswgof} and $(1,2)$ above there exists a restriction
mapping $s_{U_2\times V_2,U_1\times V_1}: F(U_1\times V_1)\to
F(U_2\times V_2)$ for each open $U_2\subset U_1$ and $V_2\subset
V_1$. Then this presheaf induces a sheaf ${\cal S}_{W,X_1,X_2,\cal
G}$ of wrap groups. If ${\cal G}$ is a twisted over $\{
i_0,...,i_{2^r-1} \}$ for $2\le r\le 3$ sheaf, then the sheaf ${\cal
S}_{W,X_1,X_2,\cal G}$ is twisted over $\{ i_0,...,i_{2^r-1} \}$. If
the sheaf $\cal G$ is commutative, then the sheaf ${\cal
S}_{W,X_1,X_2,\cal G}$ is commutative.

\par {\bf 26. Proposition.} {\it If $h_j: X_j\to Y_j$
are $H^t_p$ differentiable mappings from $X_j$ onto $Y_j$, $j=1, 2$,
where $X=X_1\times X_2$ and $Y=Y_1\times Y_2$, $X, X_1, X_2, Y, Y_1,
Y_2$ are $H^t_p$-pseudo-manifolds over ${\cal A}_r$, $0\le r\le 3$,
$h_3: {\cal G}_{Y_1}\to {\cal G}_{X_1}$ is an $H^t_p$ sheaf
homomorphism, $t\ge [\max \{ dim (X_1), dim (X_2), dim (Y_1), dim
(Y_2) \} ]/2 +2$. Then they induce homomorphisms $(h_1,h_3)_*: {\cal
S}_{W,Y_1,X_2,{\cal G}_{Y_1}} \to {\cal S}_{W,X_1,X_2,{\cal
G}_{X_1}}$ and $h_{2,*}: {\cal S}_{W,X_1,X_2,{\cal G}_{X_1}} \to
{\cal S}_{W,X_1,Y_2,{\cal G}_{X_1}}$ of wrap sheaves.}

\par {\bf Proof.} If $M_2\subset U_2\subset Y_1$, then
$h_1^{-1}(M_2) =: M_1\subset h_1^{-1}(U_2) =: U_1\subset X_1$ and
$h_1^{-1}(U_2) =: U_1$ is open in $X_1$ for each $U_2$ open in
$Y_1$. In view of Corollary 9 \cite{luwrgfbqo} and Proposition 7.1
and Theorem 10 \cite{lulaswgof} there exists a homomorphism
$(h_1,h_3)_*: (W^{M_2, \{ v_{0,q}: q=1,...,k_2 \} }E;N,G(U_2),{\bf
P})_{t,H} \to (W^{M_1, \{ s_{0,q}: q=1,...,k_1 \} }E;N,G(U_1),{\bf
P})_{t,H}$, where $h_3: G(U_2)\to G(U_1)$ is the group homomorphism,
$h_1(s_{0,q})=v_{0,a(q)}$ for each $q=1,...,k_2$, $1\le a=a(q)\le
k_2$. Choose in particular $s_{0,q}$ such that $k_1=k_2=k$.
Therefore, there exists the presheaf homomorphism $(h_1,h_3)_*:
F_{Y_1,X_2,{\cal G}_{Y_1}}(U_2\times V)\to F_{X_1,X_2,{\cal
G}_{X_1}}(U_1\times V)$ for each $U_2$ open in $Y_1$ and $V$ open in
$X_2$. This presheaf homomorphism induces the sheaf homomorphism.
\par If $f: M_1\to N_1\subset X_2$, then $h_2\circ f: M_1\to N_2$
for $H^t_p$ pseudo-manifolds $M_1$ in $X_1$, $N_1$ in $X_2$, $N_2$
in $Y_2$. If $f$ and $h_2$ are $H^t_p$ mappings, then due to the
Sobolev embedding theorem \cite{miha} for \par  $t\ge [\max \{ dim
(X_1), dim (X_2), dim (Y_1), dim (Y_2) \} ]/2 +2$ we have that $f'$
exists and is continuous almost everywhere on $X_1$ and $h_2(f(*))$
is the $H^t_p$ mapping (see also \cite{ebi}). Then $h_{2,*} ({\bf
P}_{{\hat {\gamma }},u}(x)) := {\bf P}_{h_2\circ {\hat {\gamma
}},u}(x)$ implies $h_{2,*}<{\bf P}_{{\hat {\gamma }},u}>_{t,H} = <
{\bf P}_{h_2\circ {\hat {\gamma }},u}>_{t,H}$ for classes of
$R_{t,H}$ equivalent elements, since the group $G(U)$ and the
manifold $M$ are specified, and the same for $N_1$ and $N_2$.
Therefore, there exists the induced homomorphism \\ $h_{2,*}: (W^{M,
\{ s_{0,q}: q=1,...,k \} }E;N_1,G(U),{\bf P})_{t,H} \to (W^{M, \{
s_{0,q}: q=1,...,k \} }E;N_2,G(U),{\bf P})_{t,H}$, \\ where
$N_1\subset V_1\subset X_2$, $N_1= h_2^{-1}(N_2)$, $y_{0,1}\in N_1$,
$h_2(y_{0,1})= y_{0,2}$, $y_{0,2}\in N_2\subset V_2\subset Y_2$.
Consequently, there exists the homomorphism of pre-sheaves $h_{2,*}:
F_{X_1,X_2,{\cal G}_{X_1}}(U\times V_1)\to F_{X_1,Y_2,{\cal
G}_{X_1}}(U\times V_2)$ (see \S 25), where $V_1 = h_2^{-1}(V_2)$,
$V_2$ is open in $Y_2$. Thus $h_{2,*}$ induces the homomorphism of
the wrap sheaves.

\par {\bf 27. Proposition.} {\it Let $e\to {\cal G}_1\to {\cal
G}_2\to {\cal G}_3\to e$ be an exact sequence of sheaves on $X_1$.
Then there exists an exact sequence $e\to {\cal S}_{W,X_1,X_2,{\cal
G}_1} \to {\cal S}_{W,X_1,X_2,{\cal G}_2}\to {\cal
S}_{W,X_1,X_2,{\cal G}_3}\to e$ of wrap sheaves, where $e$ is the
unit element (see \S 25).}
\par {\bf Proof.} For each $U$ open in $X_1$ there exists a short
exact sequence of groups $e\to G_1(U)\to G_2(U)\to G_3(U)\to e$ such
that $G_3(U)$ is isomorphic with the quotient group $G_2(U)/G_1(U)$,
where $G_1(U)$ is the normal closed subgroup in $G_2(U)$. In view of
Theorem 10 \cite{lulaswgof} there exists the short exact sequence
$e\to (W^ME;N,G_1(U),{\bf P})_{t,H}\to (W^ME;N,G_2(U),{\bf
P})_{t,H}\to (W^ME;N,G_3(U),{\bf P})_{t,H}\to e$. Then this induces
the short exact sequence of wrap presheaves $e\to F_{G_1(U)}(U)\to
F_{G_2(U)}(U)\to F_{G_3(U)}\to e$ and the latter in its turn gives
the short exact sequence of wrap sheaves (see also in general
\cite{bredon}).

\par {\bf 28. Wrap sub-sheaf.} In the construction of \S 25 consider
a sub-pre-sheaf corresponding to $F_N(U)$, that is, for $V=N$ with a
fixed marked point $y_0\in N$, where \par $(1)$ $F_N(U) :=
\prod_{s_{0,1},...,s_{0,k}\in M\subset U} (W^{M, \{ s_{0,q}:
q=1,...,k \} } E;N,G(U),{\bf P})_{t,H}$, \\ where $s^G_{U_2,U_1}:
G(U_1)\to G(U_2)$ is the restriction mapping for each $U_2\subset
U_1$. In view of Theorem 10 \cite{lulaswgof} there exists a
restriction mapping $s_{U_2,U_1}: F_N(U_1)\to F_N(U_2)$ for each
open $U_2\subset U_1$. Then this presheaf induces a sheaf ${\cal
S}_{W,X_1,\cal G}(N)$ of wrap groups, which is the subsheaf of
${\cal S}_{W,X_1,X_2,\cal G}$.

\par {\bf 29. Proposition.} {\it Let $\eta : N_1\to N_2$
be an $H^{t'}_p$-retraction of $H^{t'}_p$ manifolds, $N_2\subset
N_1$, $\eta |_{N_2}=id$, $y_0\in N_2$, where $t'\ge t$, $M$ is an
$H^t_p$ manifold, $E(N_1,G,\pi ,\Psi )$ and $E(N_2,G,\pi ,\Psi )$
are principal $H^{t'}_p$ bundles with a structure group $G$
satisfying conditions of \S 2 \cite{luwrgfbqo}, then there exists a
sheaf homomorphism $\eta _*$ from ${\cal S}_{W,X_1,\cal G}(N_1)$
onto ${\cal S}_{W,X_1,\cal G}(N_2)$.}
\par {\bf Proof.} In view of Proposition 17 \cite{lulaswgof}
there exists a group homomorphism $\eta _*(U)$ from $F_{N_1}(U)$
onto $F_{N_2}(U)$ for each $U$ open in $X_1$ such that $\{ s_{0,q}:
q=1,...,k \} \subset M\subset X_1$. If $\cal B$ is a sheaf on $X$
and $\eta {\cal B}(U) = {\cal B}(\eta ^{-1}(U))$ for each $U$ open
in $X$, then there exists a sheaf $\eta {\cal B}$ which is called
the image of the sheaf $\cal B$ (see \cite{bredon}). On the other
hand, $\eta _*(U_2) \circ s_{U_2,U_1}= s_{U_2,U_1}\circ \eta
_*(U_1)$ for each open $U_2\subset U_1$ due to Condition 25$(2)$.
Then ${\cal S}_{W,X_1,\cal G}(N_2)$ is the image of ${\cal
S}_{W,X_1,\cal G}(N_1)$, that is $\eta _*{\cal S}_{W,X_1,\cal
G}(N_1)={\cal S}_{W,X_1,\cal G}(N_2)$, since there exists an
$H^{t}_p$ mapping $id\times \eta $ from $M\times N_1$ onto $M\times
N_2$ (see \S 28). This gives the sheaf homomorphism (see also \S 3
\cite{bredon}).

\par {\bf 30. Remark.} For a continuous mapping $f: X\to Y$ and a
sheaf ${\cal B}$ on $Y$ a inverse image $f^*{\cal B}$ is a sheaf on
$X$ such that $f^*{\cal B} = \{ (x,q) \in X\times {\cal B}: f(x)=\pi
(q) \} $ (see \cite{bredon}). Particularly, if $f: X\to Y$ is an
$H^t_p$ mapping such that $f=(f_1,f_2)$, $f_1: X_1\to Y_1$, $f_2:
X_2\to Y_2$, then there exists a sheaf inverse image $f^*{\cal
S}_{W,Y_1,Y_2,{\cal G}_2}$, where $f_1^*{\cal G}_2 = {\cal G}_1$.

\par {\bf 31. Corollary.} {\it Let suppositions of Proposition 26 be
satisfied, where $h_j$ are diffeomorphisms for $j=1, 2$ and an
isomorphism for $j=3$, then ${\cal S}_{W,X_1,X_2,{\cal G}_{X_1}}$
and ${\cal S}_{W,Y_1,Y_2,{\cal G}_{Y_1}}$ are isomorphic sheaves.}
\par {\bf Proof.} This follows from Proposition 26 and Remark 30.

\par {\bf 32. Proposition.} {\it Let a sheaf $\cal G$ be an inductive
limit $ind-\lim_{\alpha \in \Lambda } G_{\alpha }$ of sheaves ${\cal
G}_{\alpha }$, where $\Lambda $ is a directed set. Then the wrap
sheaf ${\cal S}_{W,X_1,X_2,{\cal G}}$ is the inductive limit
$ind-\lim_{\alpha \in \Lambda } {\cal S}_{W,X_1,X_2,{\cal G}_{\alpha
}}$.}
\par {\bf Proof.} For each $U$ open in $X_1$ and all $\alpha <\beta \in
\Lambda $ there exists a homomorphism $\pi ^{\alpha }_{\beta }:
{\cal G}_{\alpha }(U)\to {\cal G}_{\beta }(U)$. Then the sheaf
${\cal G}$ is defined as the sheaf generated by a pre-sheaf
$U\mapsto ind-\lim_{\alpha \in \Lambda } {\cal G}_{\alpha }(U)$ (see
Chapter 1 \S 5 \cite{bredon}). Each homomorphism $\pi ^{\alpha
}_{\beta }$ generates the homomorphism of principal bundles from
$E(N,G_{\alpha },\pi ,\Psi )$ into $E(N,G_{\beta },\pi ,\Psi )$. In
view of Proposition 26 for each $\alpha <\beta \in \Lambda $ and
every $U$ open in $X_1$ there exists the group homomorphism $\pi
^{\alpha }_{\beta ,*}: {\cal S}_{W,X_1,X_2,{\cal G}_{\alpha }}(U)\to
{\cal S}_{W,X_1,X_2,{\cal G}_{\beta }}(U)$ generated by $\pi
^{\alpha }_{\beta }$. Thus, there exists ${\cal S}_{W,X_1,X_2,{\cal
G}} := ind-\lim_{\alpha \in \Lambda } {\cal S}_{W,X_1,X_2,{\cal
G}_{\alpha }}$.

\par {\bf 33. Corollary.} {\it Let $X_1 = ind-\lim_{\alpha \in
\Lambda } X_{1,\alpha }$ and ${\cal G} = ind-\lim_{\alpha \in
\Lambda } {\cal G}_{\alpha }$ satisfy conditions of Proposition 26,
where ${\cal G}_{\alpha } ={\cal G}_{X_{1,\alpha }}$ and
$X_{1,\alpha }$ is an $H^t_p$ pseudo-manifold for each $\alpha $ in
a directed set $\Lambda $. Then ${\cal S}_{W,X_1,X_2,{\cal G}} =
ind-\lim_{\alpha \in \Lambda } {\cal S}_{W,X_{1,\alpha },X_2,{\cal
G}_{\alpha }}$.}
\par {\bf Proof.} For an $H^t_p$ pseudo-manifold $X_1$ its
base of topology consists of all those subsets $U$ in $X_1$ such
that $U = \bigcap_{v=1}^m U_{\alpha (v)}$ for some $m\in \bf N$ and
$\alpha (1),...,\alpha (m)\in \Lambda $, where $V_{\alpha } = \pi
_{\alpha }^{-1}(U_{\alpha (v)})$ is open in $X_{\alpha }$, $\pi
_{\alpha }: X_{1,\alpha }\to X_1$ is an embedding for each $\alpha $
in $\Lambda $. In view of Proposition 26 for each $U$ open in $X_1$
and $\alpha <\beta \in \Lambda $ there exists a group homomorphism
$\pi ^{\alpha }_{\beta ,*}: {\cal S}_{W,X_{1,\alpha },X_2,{\cal
G}_{\alpha }}(U)\to {\cal S}_{W,X_{1,\beta },X_2,{\cal G}_{\beta
}}(U)$. Due to Proposition 32 this generates ${\cal
S}_{W,X_1,X_2,{\cal G}}$ as the inductive limit of sheaves ${\cal
S}_{W,X_{1,\alpha },X_2,{\cal G}_{\alpha }}$.

\par {\bf 34. Theorem.} {\it Let $X_2 = X_{2,1}\times X_{2,2}$,
where $X_1$, $X_{1,2}$ and $X_{2,2}$ are $H^t_p$ and $H^{t'}_p$
pseudo-manifolds respectively over ${\cal A}_r$ as in \S 25. Then
the restriction of the complete tensor product of wrap sheaves
${\cal S}_{W,X_1,X_{2,1},{\cal G}_1} {\hat {\otimes }} {\cal
S}_{W,X_1,X_{2,2},{\cal G}_2}$ on $\Delta _1\times X_2$ is
isomorphic with ${\cal S}_{W,X_1,X_2,{\cal G}}$, where ${\cal G}$ is
the tensor product ${\cal G} := {\cal G}_1 {\otimes } {\cal G}_2$ of
sheaves ${\cal G}_1$ and ${\cal G}_2$ on $X_1$, $\Delta _1 := \{
(x,x): x\in X_1 \} $ is the diagonal in $X_1^2$.}

\par {\bf Proof.} If ${\cal B}_1$ and ${\cal B}_2$ are sheaves on
a topological space $X$, then ${\cal B}_1\otimes {\cal B}_2$ denotes
the sheaf on $X$ generated by the presheaf $U\mapsto {\cal
B}_1(U)\otimes {\cal B}_2(U)$, where $({\cal B}_1\otimes {\cal
B}_2)_x \cong {\cal B}_{1,x}\otimes {\cal B}_{2,x}$ is the natural
isomorphism of fibers. The sheaf ${\cal B}_1\otimes {\cal B}_2$ is
called the tensor product of sheaves.
\par Consider the natural projections $\phi _1: X_2\to X_{2,1}$ and
$\phi _2: X_2\to X_{2,2}$ having extensions $id\times \phi _1:
X_1\times X_2\to X_1\times X_{2,1}$ and $id\times \phi _2: X_1\times
X_2\to X_1\times X_{2,2}$. Therefore, there exists the sheaf ${\cal
S} := {\cal S}_{W,X_1,X_{2,1},{\cal G}_1} {\hat {\otimes }} {\cal
S}_{W,X_1,X_{2,2},{\cal G}_2} := [(id\times \phi _1)^*{\cal
S}_{W,X_1,X_{2,1},{\cal G}_1}]\otimes [(id\times \phi _2)^*{\cal
S}_{W,X_1,X_{2,2},{\cal G}_2}]$ which is the complete tensor product
of sheaves (see in general Chapter 1 \S 5 \cite{bredon}).
\par If $\gamma : M\to X_2$ is an $H^t_p$ mapping preserving marked points,
then $\gamma = (\gamma _1, \gamma _2)$, where $\gamma _j: M\to
X_{2,j}$ for $j=1, 2$, $\gamma (s_{0,q})=y_0$, $\gamma
_j(s_{0,q})=y_{j,0}$ for each $q=1,...,k$ and $j=1, 2$,
$y_0=y_{1,0}\times y_{2,0}$. Then we get a lifting ${\hat {\gamma
}}: {\hat M}\to X_2$ such that $\gamma \circ \Xi = {\hat {\gamma }}$
(see \S \S 2, 3 and 6 in \cite{luwrgfbqo}). Therefore, ${\bf
P}_{{\hat {\gamma }},u} ({\hat s}_{0,k+q}) = {\bf P}_{{\hat {\gamma
}}_1,u_1} ({\hat s}_{0,k+q})\otimes {\bf P}_{{\hat {\gamma }}_2,u_2}
({\hat s}_{0,k+q}) \in G$ for each $q=1,...,k$, with $G=G_1\otimes
G_2$ being the direct product of groups for $G_1={\cal G}_1(U_1)$
and $G_2={\cal G}_2(U_2)$ for every $U_j$ open in $X_1$, $j=1, 2$,
where $u\in E_{y_0}$, $u_j\in E_{j,y_{j,0}}$, $N=N_1\times N_2$,
$N_j\subset V_j\subset X_{2,j}$, $E=E(N,G,\pi ,\Psi )$, $E_j
=E(N_j,G_j,\pi _j,\Psi _j)$ are principal bundles,
$y_0=y_{0,1}\times y_{0,2}$, $y_{j,0}\in N_j$ are marked points,
$V_j$ is open in $X_{2,j}$ for $j=1, 2$ (see also \S 25). \par For
classes of equivalent parallel transport structures we get $<{\bf
P}_{{\hat {\gamma }},u}>_{t,H} = <{\bf P}_{{\hat {\gamma
}}_1,u_1}>_{t,H}\otimes <{\bf P}_{{\hat {\gamma }}_2,u_2}>_{t,H}$,
hence $F(U\times (V_1\times V_2))$ is isomorphic with $(\phi
_1)^*F(U\times V_1)\otimes (\phi _2)^*F(U\times V_2)$ for each $U$
open in $X_1$ and all $V_j$ open in $X_{2,j}$, $j=1, 2$, since open
sets of the form $V=V_1\times V_2$ form a base of topology in $X_2$,
where $F(U\times V_j)$ is given for the group ${\cal G}_j(U)$. Here
$U=U_1=U_2$ and $(\phi _1)^*F(U\times V_1)\otimes (\phi
_2)^*F(U\times V_2)$ is isomorphic with the restriction of
$(id\times \phi _1)^*F(U\times V_1)\otimes (id\times \phi
_2)^*F(U\times V_2)$ from $U^2\times V_1\times V_2$ onto $\Delta
(U)\times V_1\times V_2$, where $\Delta (U)$ denotes the diagonal in
$U^2$. Thus, ${\cal S}_{W,X_1,X_2,{\cal G}}$ is isomorphic with the
restriction of the complete tensor product of sheaves ${\cal
S}_{W,X_1,X_{2,1},{\cal G}_1} {\hat {\otimes }} {\cal
S}_{W,X_1,X_{2,2},{\cal G}_2}$ on $\Delta _1\times X_2$.

\par {\bf 35. Twisted Alexander-Spanier cohomologies.}
\par Let $G$ be a group satisfying Conditions 4$(A1,A2)$, which may be
in particular a wrap group for ${\cal A}_r$ pseudo-manifolds with
$2\le r\le 3$. For an open $U$ in $X$ denote by $A^m(U;G)$ a group
of all functions $f: U^{m+1}\to G$ with a pointwise multiplication
in $G$ as the group operation. Therefore, the functor $U\mapsto
A^m(U;G)$ is a presheaf on $X$ satisfying the condition:
\par $(S2)$ if $ \{ U_j: j \} $ is a family of open subsets of $X$
such that $\bigcup_j U_j=U$, then for a family of elements $s_j\in
A^m(U_j;G)$ such that $s_j|_{U_j\cap U_k} = s_k|_{U_j\cap U_k}$ for
each $j, k$ there exists $s\in A^m(U;G)$ such that $s|_{U_j}=s_j$
for each $j$. To satisfy this put $s_j=f_j: U_j^{m+1}\to G$ to be
functions here and $s=f$ is their combination such that
$f|_{U_j^{m+1}}=s_j$, while $f$ on $X^{m+1}\setminus (\bigcup_j
U_j^{m+1})$ is arbitrary. The property
\par $(S1)$ if $U=\bigcup_j U_j$, where $U_j$ is open in $X$
and $f, g\in A^{m+1}(U;G)$ coincide on $U_j^{m+1}$ for each $j$,
then $f=g$ on $\bigcup_j U_j^{m+1}$ is evident, since $f, g$ are
functions.
\par Recall that a family $\phi $ of closed subsets in $X$ is called
a family of supports, if it satisfies conditions $(SP1,SP2)$:
\par $(SP1)$ if $B$ is a closed subset in $C$, where
$C\in \phi $, then $B\in \phi $;
\par $(SP2)$ if $B_1,...,B_m\in \phi $, $m\in \bf N$, then
$\bigcup_{j=1}^m B_j\in \phi $.
\par The family $\phi $ of supports is called paracompactfying,
if satisfies two additional conditions:
\par $(SP3)$ each element in $\phi $ is a paracompact space;
\par $(SP4)$ each set from $\phi $ has a closed neighborhood
belonging to $\phi $.

\par The union $\bigcup_{C\in \phi} C =: {\sf E}(\phi )$ is called
a spread of $\phi $. Put $\Gamma _{\phi } ({\cal S}) := \{ s\in
{\cal S}(X): |s|\in \phi \} $ for a sheaf ${\cal S}$ on $X$, where
$|s| := \{ x\in X: s(x)\ne e \} $ denotes its support. Clearly
$\Gamma _{\phi }({\cal S})$ is a subgroup in ${\cal S}(X)$. For a
presheaf $A$ on $X$ put $A_{\phi } (X) := \{ s\in A(X): |s|\in \phi
\} $. For a presheaf $A$ on $X$ put $A_{\phi }(X) := \{ s\in A(X):
|s|\in \phi \} $.
\par Let now ${\sf A}^m(X;G)$ be a sheaf generated by the presheaf
$A^m(.;G)$. Define the differential $d: A^m(U;G)\to A^{m+1}(U;G)$ by
the formula:
\par $d f(x_0,...,x_{m+1}) = \sum_{j=0}^{m+1} (-1)^j f(x_0,...,{\hat
x}_j,...,x_{m+1})$, where $f: U^{m+1}\to G$ is an arbitrary
function. Then $f = \sum_{k=0}^{2^r-1} f_k i_k$, where $f_k\in {\hat
G}_k$, $\{ i_0,..,i_{2^r-1} \} $ are generators of ${\cal A}_r$,
$2\le r\le 3$. Hence $d$ is the homomorphism of presheaves and
$d^2=0$, since $df = \sum_{k=0}^{2^r-1} (df_k)i_k$.

Then twisted Alexander-Spanier cohomologies are defined as
\par $\mbox{ }_{AS}{\sf H}^m_{\phi } (X;G) = {\sf H}^m(A^*_{\phi
}(X;G))/A_0^*(X;G))$.

\par {\bf 36. Theorem.} {\it Let $A$ be a pre-sheaf on $X$ satisfying
Condition 35$(S2)$ and ${\cal S}$ be a sheaf generated by $A$, where
$\cal S$ and $A$ are twisted over $ \{ i_0,...,i_{2^r-1} \} $ with
$1\le r\le 3$. Then for each paracompatifying family $\phi $ of
supports in $X$ there exists the exact sequence
\par $e\to A_0(X)\to A_{\phi } (X) \stackrel{\theta}\longrightarrow
\Gamma _{\phi }({\cal S})\to e$, where $\theta : A(X)\to {\cal
S}(X)$ is the natural mapping of the presheaf into the generated by
it sheaf.}
\par {\bf Proof.} Consider $s\in \Gamma _{\phi }({\cal S})$
and a neighborhood $U$ of $|s|$ such that $cl (U)\in \phi $, where
$cl (U)$ denotes the closure of $U$ in $X$. Since $cl (U)$ is
paracompact find a locally finite covering $ \{ U_j: j \} $ of $cl
(U)$, where each $U_j$ is open in $X$ and for which there exists
$s_j\in A(U_j)$ such that $\theta (s_j) = s|_{U_j}$. Let $\{ V_j: j
\} $ be a refinement of $ \{ U_j: j \} $ such that $U\cap cl (V_j)
\subset U_j$. \par For $x\in X$ the set $J(x) := \{ j: x\in cl (V_j)
\} $ is finite, hence for each $x\in X$ there exists a neighborhood
$W(x)$ such that $W(x)\subset U_j$ and for each $j\in J(x)$ and
every $y\in W(x)$ there is the inclusion $J(y)\subset J(x)$.

\par For $j\in J(x)$ we get $\theta (s_j(x)) = s(x)$. Take $W(x)$
sufficiently small such that $s_j|_{W(x)}=: s_x$ does not depend on
$j\in J(x)$, since $J(x)$ is finite, consequently, $s_x\in A(W(x))$.

\par Let $x, y\in U$, $z\in W(x)\cap W(y)$ and $j\in J(z)$, where
$J(z)\subset J(x)\cup J(y)$. Then $s_x|_{W(x)\cap W(y)} =
s_y|_{W(x)\cap W(y)}$. Due to Condition $(S2)$ there exists $\beta
\in A(U)$ such that $\beta |_{W(x)} =s_x$ for each $x\in U$,
clearly, $\theta (\beta )= s|_U$.

\par Take now $C\in \phi $ such that $|s|\subset Int (C)$ and
$C\subset U$, where $Int (C)$ denotes the interior of $C$. If $x\in
C\setminus Int (C)$, then $\theta (\beta ) (x) =s(x)=0$. Therefore,
there exists a covering $\{ Q_j \} $ of $C\setminus Int (C)$ with
open in $X$ sets $Q_j$ such that $Q_j\subset U$ and $t|_{Q_j} = 0$
for every $j$.
\par Choose the open covering $ \{ Q_j \} \cup \{ Int (C),
X\setminus C \} $ of $X$ and elements $e\in A(Q_j)$, $\beta |_{Int
(C)} \in A (Int (C))$ and $e\in A(X\setminus C)$. Restrictions of
each two elements on a common part of their domains of definition
coincide. Hence due to Condition $(S2)$ such elements have a common
extension $q\in A(X)$ and inevitably $\theta (q) = s$ and $|q| =
|\theta (q)| = |s| \in \phi $. The sequence $e\to A_0(X)\to A_{\phi
} (X) \stackrel{\theta}\longrightarrow \Gamma _{\phi }({\cal S})\to
e$ is exact, since each subsequence $e\to A_{0,k}(X)\to A_{\phi ,k}
(X) \stackrel{\theta}\longrightarrow \Gamma _{\phi ,k}({\cal S})\to
e$ is exact, where ${\hat A}_{\phi } = \sum_{k=0}^{2^r-1} {\hat
A}_{\phi ,k} i_k$ and ${\hat \Gamma }_{\phi } = \sum_{k=0}^{2^r-1}
{\hat \Gamma }_{\phi ,k} i_k$, where each ${\hat A}_{\phi ,k}$ is
commutative and they are pairwise isomorphic for different $k$, as
well as ${\hat \Gamma }_{\phi ,k}$ are commutative and pairwise
isomorphic for different values of $k$, since the sheaf ${\cal S}$
is twisted over the group of standard generators $ \{
i_0,...,i_{2^r-1} \} $ of the Cayley-Dickson algebra ${\cal A}_r$.

\par Mention, that for a pre-sheaf $A$
satisfying Condition $(S1)$ we have $A_0(X)=e$.

\par {\bf 37. Corollary.} {\it Let  conditions of Theorem 36 be satisfied.
Then for a paracompactifying family $\phi $ of supports there exists
the natural isomorphism:
\par ${\sf H}^m_{\phi } (X;G) \cong {\sf H}^m(\Gamma _{\phi }
({\cal S}^*(X;G))$.}
\par {\bf Proof.} This follows immediately from Theorem 36 and \S
35.

\par {\bf 38. Twisted singular cohomologies.}
\par Let $\cal B$ be a locally finite twisted sheaf on $X$, that is a group
${\cal B}(U)$ satisfies Conditions 4$(A1,A2,C1,C2)$ for each $U$
open in $X$. For $U\subset X$ denote by $S^m(U;{\cal B})$ the group
of singular $m$-dimensional cochains of the space $U$ with
coefficients in $\cal B$. Each element $f\in S^m(U;{\cal B})$ is a
function posing for each $m$-dimensional simplex $\sigma : \Delta
^m\to U$ a section $f(\sigma )\in \Gamma (\sigma ^*({\cal B}))$,
where $\Delta ^m$ is a standard $m$-dimensional simplex. \par The
pre-sheaf $S^m(. ;{\cal B})$ satisfies Condition $(S2)$. The sheaf
$\cal B$ is locally constant, then the sheaf $\sigma ^*({\cal B})$
is constant on $\Delta ^m$, since the simplex $\Delta ^m$ is simply
connected, where $m\ge 1$. Therefore, there exists a usual
coboundary operator $d: S^m(U;{\cal B})\to S^{m+1}(U;{\cal B})$.
\par Consider the sheaf ${\cal S}^m(U;{\cal B})$ generated by a
pre-sheaf $U\mapsto S^m(U;{\cal B})$. Then the differential $d$ in
the pre-sheaf induces the differential in the sheaf. For a locally
constant sheaf ${\cal B}$ singular cohomologies with coefficients in
$\cal B$ and supports in the family $\phi $ are defined as $\mbox{
}_{\Delta } {\sf H}^m_{\phi } (X;{\cal B}) = {\sf H}^m(S^*_{\phi
}(X;{\cal B}))$. Since $\cal B$ is the twisted sheaf over $\{
i_0,...,i_{2^r-1} \} $, then $S^*_{\phi }(X;{\cal B})$ and
inevitably $\mbox{ }_{\Delta } {\sf H}^m_{\phi } (X;{\cal B})$ are
twisted over $ \{ i_0,...,i_{2^r-1} \} $.

\par Let ${\cal U} := \{ U_j: j \} $ be an open covering of $X$ and
let $S^*({\cal U};{\cal B})$ be a group of singular cochains defined
on singular simplices subordinated to the covering $\cal U$. With
the help of the subdivision we get, that the homomorphism $b_{\cal
U}: S^*(X;{\cal B})\to S^*({\cal U};{\cal B})$ induces the
isomorphism of cohomologies, consequently, the complex $K_{\cal U}
=ker b_{\cal U}$ is acyclic. On the other hand, $S_0^*(X;{\cal B}) =
\bigcup K^*_{\cal U} = ind-\lim K_{\cal U}$, hence ${\sf
H}^*(S_0^*(X;{\cal B})) = {\sf H}^* (ind-\lim K^*_{\cal U}) =
ind-\lim {\sf H}^*(K^*_{\cal U})=e$.
\par Thus for a paracompatifying family of supports from
the exactness of the sequence \par $e\to S_0^*\to S^*_{\phi }\to
\Gamma _{\phi }({\cal S}^*)\to e$ \\ and Theorem 36 it follows the
isomorphism $\mbox{ }_{\Delta } {\sf H}^m_{\phi } (X;{\cal B}) \cong
{\sf H}^m(\Gamma _{\phi }({\cal S}^*(X;{\cal B})))$.

\par {\bf 39. Twisted differential sheaves.}
\par A graded sheaf is a sequence $ \{ {\cal S}^m: m\in {\bf Z} \} $
of sheaves, which is called a differential sheaf if there are
homomorphisms \par $(1)$ $d: {\cal S}^m\to {\cal S}^{m+1}$ such that
$d^2=0$ for each $m$. This sheaf may be twisted over $\{
i_0,...,i_{2^r-1} \} $, where $2\le r\le 3$. In this case we suppose
that \par $(2)$ up to an automorphisms $\theta _m: {\cal S}^m\to
{\cal S}^m$ we have $\theta _{m+1}\circ d({\cal S}^m_k)\subset {\cal
S}^{m+1}_k$ for each $k=0,...,2^r-1$.

\par A differential sheaf ${\cal S}^*$
having ${\cal S}^m = 0$ for each $m<0$ and supplied with the
augmentation homomorphism $\varepsilon : {\cal B}\to {\cal S}^0$ is
called the resolvent of the sheaf, if the sequence
\par $e\to {\cal B} \stackrel{\varepsilon}\longrightarrow
{\cal S}^0 \stackrel{d}\longrightarrow {\cal S}^1
\stackrel{d}\longrightarrow {\cal S}^2\to ...$ \\
is exact.

\par The notion of differential graded pre-sheaves is formulated
analogously. If ${\cal S}^m$ is twisted, that is ${\hat {\cal S}}^m
= {\hat {\cal S}}^m_0i_0+...+{\hat {\cal S}}^m_{2^r-1}i_{2^r-1}$,
where ${\hat {\cal S}}^m_k$ and ${\hat {\cal S}}^m_j$ are pairwise
isomorphic and commutative for each $k\ne j$, then $Ker (d: {\cal
S}^m\to {\cal S}^{m+1})$ and $Im (d: {\cal S}^{m-1}\to {\cal S}^m)$
are twisted as well, since up to isomorphisms $\theta _m: {\cal
S}^m\to {\cal S}^m$ we have $\theta _{m+1}\circ d({\cal
S}^m_k)\subset {\cal S}^{m+1}_k$ for each $k=0,...,2^r-1$.
\par A sheaf of cohomologies (in another words a derivative sheaf)
is defined as ${\sf H}^m({\cal S}^*) = Ker (d: {\cal S}^m\to {\cal
S}^{m+1})/ Im (d: {\cal S}^{m-1}\to {\cal S}^m)$. If ${\cal S}^*$ is
generated by a differential pre-sheaf $S^*$, then ${\sf H}^m({\cal
S}^*)$ is generated by the pre-sheaf $U\mapsto {\sf H}^m(S^*(U))$.

\par For a sheaf ${\cal B}$ on topological space $X$ and an open
subset $U\subset X$ denote by ${\sf Y}^0(U;{\cal B})$ a set of all
mappings (may be discontinuous) $f: U\to \cal B$ such that $\pi
\circ f=id$ is the identity mapping on $U$, where $\pi : {\cal B}\to
X$ is the canonical projection. Thus ${\sf Y}^0(U;{\cal B}) =
\prod_{x\in U} {\cal B}_x$ and it is the group with the pointwise
group operation. Therefore, $U\mapsto {\sf Y}^0(U;{\cal B})$ is the
pre-sheaf satisfying Conditions $(S1,S2)$, hence it is the sheaf
which we denote by ${\cal Y}^0(X;{\cal B})$. If $\cal B$ is twisted,
then ${\cal Y}^0(U;{\cal B})$ is twisted as well.
\par The inclusion of all continuous sections of $\cal B$ into the
family of all sections not necessarily continuous induces the
augmentation homomorphism $\varepsilon : {\cal B}\to {\cal
Y}^0(X;{\cal B})$.
\par For a family $\phi $ of supports put ${\sf Y}^0_{\phi }
(X;{\cal B})= \Gamma _{\phi } ({\cal Y}^0(X;{\cal B})$. If $e\to
{\cal B}_1\to {\cal B}_2\to {\cal B}_3\to e$ is a short exact
sequence of sheaves (may be twisted), then the sequence of
pre-sheaves $e\to {\sf Y}^0(X;{\cal B}_1)\to {\sf Y}^0(X;{\cal
B}_2)\to {\sf Y}^0(X;{\cal B}_3)\to e$ is exact. If $f\in {\sf
Y}^0_{\phi }(X;{\cal B})$, then its support is $|f| := cl \{ x:
f(x)\ne e \} $. Therefore, $f$ is an image of a section $g$ of the
sheaf $\cal B$ such that $g$ is not necessarily continuous and
$g(x)=e$ if $f(x)=e$ for $x\in X$, hence $|g| = |f|\in \phi $.
\par Denote by ${\cal Z}^1(X;{\cal B})$ the cokernel of the
homomorphism $\varepsilon $ such that the sequence $e\to {\cal B}
\stackrel{\varepsilon}\longrightarrow {\cal Y}^0(X;{\cal B})
\stackrel{\partial }\longrightarrow {\cal Z}^1(X;{\cal B})$ is
exact. Define by induction the sheaves \par ${\cal Y}^m(X;{\cal
B})={\cal Y}^0(X;{\cal Z}^m(X;{\cal B}))$, ${\cal Z}^{m+1} (X;{\cal
B})= {\cal Z}^1(X;{\cal Z}^m(X;{\cal B}))$. \par If $\cal B$ is
twisted over $ \{ i_0,..,i_{2^r-1} \} $, then ${\cal Z}^1(X;{\cal
B})$ is twisted as well and by induction ${\cal Z}^m(X;{\cal B})$
and ${\cal Y}^m(X;{\cal B})$ are twisted for each $m\in \bf N$.
Therefore, the sequence $e\to {\cal Z}^m(X;{\cal
B})\stackrel{\varepsilon}\longrightarrow {\cal Y}^m (X;{\cal B})
\stackrel{\partial }\longrightarrow {\cal Z}^{m+1}(X;{\cal B})\to e$
is exact. Consider the composition $d= \varepsilon \circ
\partial $ for ${\cal Y}^m(X;{\cal B})\stackrel{\partial
}\longrightarrow {\cal Z}^{m+1} (X;{\cal B}) \stackrel{\varepsilon
}\longrightarrow {\cal Y}^{m+1}(X;{\cal B})$, then the sequence
\par $e\to {\cal B} \stackrel{\varepsilon}\longrightarrow {\cal
Y}^0(X;{\cal B})\stackrel{d}\longrightarrow {\cal Y}^1 (X;{\cal B})
\stackrel{d}\longrightarrow {\cal Y}^2(X;{\cal B})\to ...$ is exact.
Thus, ${\cal Y}^*(X;{\cal B})$ is the resolvent of the sheaf $\cal
B$, which is called the canonical resolvent.

\par {\bf 40. Proposition.} {\it The canonical resolvent of the twisted
sheaf $\cal B$ is fiberwise homotopically trivial.}
\par {\bf Proof.} Consider the homomorphism ${\sf Y}^0(U;{\cal
B})\to {\cal B}_x$ such that $U\ni x\mapsto f(x)\in {\cal B}_x$ for
each $f\in {\sf Y}^0(U;{\cal B})$ and $x\in U$. The direct limit by
neighborhoods of a point $x$ induces the homomorphism $\eta _x:
{\cal Y}^0(X;{\cal B})_x\to {\cal B}_x$, consequently, $\eta _x\circ
\varepsilon : {\cal B}_x\to {\cal B}_x$ is the identity isomorphism,
where $\eta _x\circ \varepsilon (z) = \eta _x(\varepsilon (z))$.
Define the homomorphism $\nu _x: {\cal Z}^1(X;{\cal B})\to {\cal
Y}^0(X;{\cal B})$ by the formula $\nu _x\circ \partial = 1 -
\varepsilon \circ \eta _x$ which defines $\nu _x$ in a unique way.
Therefore, there exists a fiber splitting
\par ${\cal Z}^m(X;{\cal B})_x \stackrel{\varepsilon}\longrightarrow
{\cal Y}^m(X;{\cal B}) \stackrel{\partial } \longrightarrow  {\cal
Z}^{m+1}(X;{\cal B})_x$ and  \par ${\cal Z}^m(X;{\cal B})_x
\stackrel{\eta _x}\longleftarrow {\cal Y}^m(X;{\cal B})
\stackrel{\nu _x} \longleftarrow  {\cal Z}^{m+1}(X;{\cal B})_x$. Put
$D_x := \nu _x\circ \eta _x: {\cal Y}^m(X;{\cal B})_x\to {\cal
Y}^{m-1}(X;{\cal B})_x$ for $m>0$. Therefore, $d\circ D_x + D_x\circ
d = \varepsilon \circ \partial \circ \nu _x \circ \eta _x + \nu _x
\circ \eta _x \circ \varepsilon \circ
\partial = \varepsilon \circ \eta _x + \nu _x\circ \partial =1$ on ${\cal
Y}^m(X;{\cal B})$ for $m>0$. At the same time on ${\cal Y}^0(X;{\cal
B})_x$ we have $D_x\circ d = \nu _x\circ \eta _x \circ \varepsilon
\circ \partial = \nu _x\circ \partial = 1 - \varepsilon \circ \eta
_x$. This means, that ${\cal Y}^*(X;{\cal B})_x$ is homotopically
fiberwise trivial resolvent.

\par {\bf 41. Remark.} The functor ${\cal Y}^0(X;{\cal B})$ is
exact by $\cal B$, hence ${\cal Z}^1(X;{\cal B})$ is also the exact
functor by $\cal B$. Using induction we get, that all functors
${\cal Y}^m(X;{\cal B})$ and ${\cal Z}^m(X;{\cal B})$ are exact by
$\cal B$. For an arbitrary family $\phi $ of supports on $X$ put
${\sf Y}^m_{\phi } (X;{\cal B}) := \Gamma _{\phi } ({\cal
Y}^m(X;{\cal B})) = {\sf Y}^0_{\phi }(X;{\cal Z}^m(X;{\cal B}))$.
Since the functor ${\sf Y}^0(X;*)$ is exact, then the functor ${\sf
Y}^m_{\phi }(X;{\cal B})$ is exact.
\par {\bf 42. Definition.} Cohomologies in $X$ with supports in $\phi $
with coefficients in $\cal B$ are defined as ${\sf H}^m_{\phi
}(X;{\cal B}) := {\sf H}^m ({\sf Y}^*_{\phi }(X;{\cal B}))$.
\par {\bf 42.1. Note.}
The sequence $e\to \Gamma _{\phi }({\cal B})\to \Gamma _{\phi }
({\cal Y}^0(X;{\cal B}))\to \Gamma _{\phi }({\cal Y}^1(X;{\cal B}))$
is exact, consequently, $\Gamma _{\phi }({\cal B}) \cong {\sf
H}_{\phi }^0 (X;{\cal B})$. If there is a short exact sequence of
twisted sheaves $e\to {\cal B}_1\to {\cal B}_2\to {\cal B}_3\to e$
on $X$, then it implies the exact sequence of cochain complexes
\par $e\to {\sf Y}^*_{\phi }(X;{\cal B}_1)\to
{\sf Y}^*_{\phi }(X;{\cal B}_2) \to {\sf Y}^*_{\phi }(X;{\cal B}_3)
\to e$, that in its turn induces the long exact sequence
\par $...\to {\sf H}^m_{\phi }(X;{\cal B}_1) \to
{\sf H}^m_{\phi }(X;{\cal B}_2) \to {\sf H}^m_{\phi }(X;{\cal B}_3)
\stackrel{\delta }\longrightarrow {\sf H}^{m+1}_{\phi }(X;{\cal
B}_1)\to ... $.

\par {\bf 43. Definition.} Let $G$ be a topological group
satisfying conditions 4$(A1,A2,C1,C2)$ such that $G$ is a
multiplicative group of the ring $\hat G$, where $1\le r\le 2$. Then
define the smashed product $G^s$ such that it is a multiplicative
group of the ring ${\hat G}^s := {\hat G}\otimes _l{\hat G}$, where
$l=i_{2^r}$ denotes the doubling generator, the multiplication in
${\hat G}\otimes _l{\hat G}$ is \par $(1)$ $(a+bl)(c+vl) = (ac - v^*
b) + (va+bc^*)l$ for each $a, b, c, v\in \hat G$, where $v^* = conj
(v)$.
\par A smashed product $M_1\otimes _lM_2$ of manifolds $M_1, M_2$
over ${\cal A}_r$ with $dim (M_1)= dim (M_2)$ is defined to be an
${\cal A}_{r+1}$ manifold with local coordinates $z=(x,yl)$, where
$x$ in $M_1$ and $y$ in $M_2$ are local coordinates.

\par {\bf 44. Theorem.} {\it There exists smashed products
${\cal S}^s := {\cal S}_1\otimes _l {\cal S}_2$ on $X=X_1=X_2$ and
${\hat {\cal S}}^s := {\cal S}_1{\hat \otimes }_l {\cal S}_2$ on
$X=X_1\times X_2$ over $\{ i_0,...,i_{2^{r+1}-1} \} $ of isomorphic
twisted sheaves ${\cal S}_1$ on $X_1$ and ${\cal S}_2$ on $X_2$ over
$\{ i_0,...,i_{2^r-1} \} $ with $X_1=X_2$, in particular of wrap
sheaves, where $1\le r \le 2$, $l=i_{2^r}$.}
\par {\bf Proof.} If ${\cal S}_j$ is a sheaf on a topological space
$X_j$ twisted over $ \{ i_0,...,i_{2^r-1} \} $, then ${\hat {\cal
S}}_j = {\hat {\cal S}}_{0,j}i_0\oplus ... \oplus {\hat {\cal
S}}_{2^r-1,j} i_{2^r-1}$, where ${\hat {\cal S}}_{k,j}(U) = {\cal
S}_{k,j}(U)\cup \{ 0 \} $ are commutative rings for each $U$ open in
$X_j$, ${\hat {\cal S}}_{k,j}$ are sheaves on $X_j$ pairwise
isomorphic for different values of $k$. Then for $X=X_1=X_2$ take
${\cal S}_x^s := ({\cal S}_1)_x \otimes _l ({\cal S}_2)_x$ for each
$x\in X$ in accordance with Definition 43, that defines the twisted
sheaf $\cal S$ on $X$ over $\{ i_0,...,i_{2^{r+1}-1} \} $ due to
Proposition 19 \cite{lulaswgof}. This sheaf $\cal S$ is the smashed
tensor product of sheaves. \par If $X = X_1\times X_2$, then take
${\hat {\cal S}}^s := {\cal S}_1 {\hat {\otimes }} _l {\cal S}_2 =
(\pi ^*_1{\cal S}_1)\otimes _l (\pi ^*_2{\cal S}_2)$, which is the
smashed complete tensor product of sheaves, where $\pi _1: X\to X_1$
and $\pi _2: X\to X_2$ are projections.

\par {\bf 45. Corollary.} {\it Let $X_2 = X_{2,1}\otimes _lX_{2,2}$
be the smashed product, where $X_1$ and $X_2$ are $H^t_p$ and
$H^{t'}_p$ pseudo-manifolds respectively over ${\cal A}_{r+1}$,
$1\le r\le 2$. Then the restriction the smashed complete tensor
product of wrap sheaves ${\cal S}_{W,X_1,X_{2,1},{\cal G}} {\hat
{\otimes }}_l {\cal S}_{W,X_1,X_{2,2},{\cal G}}$ on $\Delta _1\times
X_2$ is isomorphic with ${\cal S}_{W,X_1,X_2,{\cal G}^s}$, where
${\cal G}^s$ is the smashed tensor product ${\cal G}^s := {\cal G}
{\otimes }_l {\cal G}$ twisted over $ \{ i_0,...,i_{2^{r+1}-1} \} $
of a sheaf ${\cal G}$ twisted over $ \{ i_0,...,i_{2^r-1} \} $ on
$X_1$, $\Delta _1 := \{ (x,x): x\in X_1 \} $ is the diagonal in
$X_1^2$.}
\par {\bf Proof.} The smashed product of manifolds was described
in details in the proof of Theorem 20 \cite{lulaswgof}. Consider an
${\cal A}_r$ shadow of $X_1$ that exists, since ${\cal A}_{r+1} =
{\cal A}_r \oplus {\cal A}_rl$, where $l=i_{2^r}$. For each $U$ open
in $X_1$ there exists a group ${\cal G}(U)$, hence ${\cal
G}(U)\otimes _l {\cal G}(U)$ is defined due to Proposition 19
\cite{lulaswgof}, that gives the sheaf ${\cal G}^s$ on $X_1$. Then
wrap sheaves ${\cal S}_{W,X_1,X_{2,b},{\cal G}}$ over ${\cal A}_r$
are defined, where $b= 1, 2$. Thus the statement of this corollary
follows from Proposition 19 \cite{lulaswgof} and Theorem 44,
modifying the proof of \S 34 for the smashed complete tensor product
instead of complete tensor product so that ${\bf P}_{{\hat {\gamma
}},u}({\hat s}_{0,k+q}) = {\bf P}_{{\hat {\gamma }}_1,u_1}({\hat
s}_{0,k+q})\otimes _l {\bf P}_{{\hat {\gamma }}_2,u_2}({\hat
s}_{0,k+q})\in G^s$ with $E=E(N,G^s,\pi ,\Psi )$, where $G= {\cal
G}(U)$, $U=U_1=U_2$, consequently, $<{\bf P}_{{\hat {\gamma
}},u}>_{t,h} = <{\bf P}_{{\hat {\gamma }}_1,u_1}>_{t,H}\otimes _l
<{\bf P}_{{\hat {\gamma }}_2,u_2}>_{t,H}$.

\par {\bf 46. Corollary.} {\it Let $X_1= X_{1,1}\otimes _lX_{1,2}$ and
$X_2 = X_{2,1}\otimes _lX_{2,2}$ are smashed products, where $X_1$,
and $X_2$ are $H^t_p$ and $H^{t'}_p$ pseudo-manifolds respectively
over ${\cal A}_{r+1}$, $1\le r\le 2$. Then the wrap sheaf ${\cal
S}_{W,X_1,X_2,{\cal G}^s}$ is twisted over $ \{
i_0,...,i_{2^{r+1}-1} \} $ and is isomorphic with the smashed
complete tensor product of twice iterated wrap sheaves \par ${\cal
S}_{W,X_{1,2},X_{2,1}, {\cal S}_{W,X_{1,1},X_{2,1},{\cal G}} } {\hat
{\otimes }}_l {\cal S}_{W,X_{1,2},X_{2,2}, {\cal
S}_{W,X_{1,1},X_{2,2},{\cal G}} }$, \\ where ${\cal G}^s$ is the
smashed tensor product ${\cal G}^s := {\cal G} {\otimes }_l {\cal
G}$ of a twisted sheaf ${\cal G}$ over $ \{ i_0,...,i_{2^r-1} \} $
on $X_1$.}
\par {\bf Proof.} Consider projections $ \pi _{b,j}: X_b\to
X_{b,j}$, where $j, b=1, 2$. Each ${\cal A}_{r+1}$ manifold has the
shadow which is the ${\cal A}_r$ manifold, since ${\cal A}_{r+1}
={\cal A}_r \oplus {\cal A}_rl$. If $U$ is open in $X_{1,j}$, then
$\pi _{1,j}^{-1}(U)$ is open in $X_1$ and there exists a group
${\cal G}(\pi _{1,j}^{-1}(U))$, where $j=1, 2$. \par Hence there
exist the projection sheaves ${\cal G}_j = \pi ^{-1}_{1,j}{\cal G}$
on $X_{1,j}$ induced by ${\cal G}$ such that ${\cal G}_j(U) := {\cal
G}(\pi _{1,j}^{-1}(U))$. Denote ${\cal G}_j$ on $X_{1,j}$ also by
${\cal G}$, since ${\cal G}_j$ is obtained from $\cal G$ by taking
the specific subfamily of open subsets. For $U_1$ open in $X_{1,1}$
and $U_2$ open in $X_{1,2}$ take $U=U_1\times U_2$ open in $X_1$.
The family of all such subsets gives the base of the topology in
$X_1$. \par In accordance with Definition 43 there exists ${\hat
{\cal G}}(U)\otimes _l{\hat {\cal G}}(U) =: {\hat {\cal G}}^s(U)$,
that induces ${\cal G}^s$ on $X_1$ such that ${\hat {\cal G}}^s_x =
{\hat {\cal G}}_x\otimes _l{\hat {\cal G}}_x$ for each $x\in X_1$.
Therefore, every element $q+vl$ is in ${\hat {\cal G}}^s(U)$ for
each $q, v\in {\hat {\cal G}}(U)$. Thus the statement of this
corollary follows from \S 25, Theorems 20 \cite{lulaswgof} and 44.

\par {\bf 47.} Consider now the iterated wrap sheaf
${\cal S}_{W,X_1,X_2,{\cal G};b}$ of iterated wrap groups
$(W^ME)_{b,\infty ,H}$ with $b\in \bf N$ instead of wrap groups for
$b=1$ such that for its presheaf \par $(1)$ $F_b(U\times V) =
\prod_{s_{0,1},...,s_{0,k}\in M\subset U; y_0\in N\subset V} (W^{M,
\{ s_{0,q}: q=1,...,k \} } E;N,G(U),{\bf P})_{b;\infty ,H}$, \\
where $s_{U_2,U_1}: G(U_1)\to G(U_2)$ is the restriction mapping for
each $U_2\subset U_1$ so that the parallel transport structure for
$M\subset U$ is defined, where ${\cal G}$ is the sheaf on $X_1$,
$G(U)={\cal G}(U)$, pseudo-manifolds $X_1$ and $X_2$ and the sheaf
$\cal G$ are of class $H^{\infty }_p$ (see also \S 25).

\par {\bf Corollary.} {\it There exists a homomorphism of
of iterated wrap sheaves $\theta : {\cal S}_{W,X_1,X_2,{\cal G};a}
\otimes {\cal S}_{W,X_1,X_2,{\cal G};b} \to {\cal
S}_{W,X_1,X_2,{\cal G};a+b}$ for each $a, b\in \bf N$. Moreover, if
$\cal G$ is either associative or alternative, then $\theta $ is
either associative or alternative.}
\par {\bf Proof.} For pre-sheaves the mapping
\par $(2)$ $\theta : F_a(U\times
V)\otimes F_b(U\times V)\to F_{a+b}(U\times V)$ \\ is induced by
Formula 47$(1)$ and due to Theorem 21 \cite{lulaswgof}. Then $\theta
$ has the extension on the sheaf of iterated wrap groups, since
$({\cal S}_{W,X_1,X_2,{\cal G};a})_z= ind-\lim F_a(U\times V)$,
where the direct limit is taken by open subsets $U\times V$ for a
point $z=x\times y\in X_1^k\times X_2$, $x\in X_1^k$, $y\in X_2$,
such that $x\subset U$, $y\in V$, $U$ is open in $X_1$, $V$ is open
in $X_2$. \par The inductive limit topology in $({\cal
S}_{W,X_1,X_2,{\cal G};a})_z$ is the finest topology relative to
which each embedding $F_a(U\times V)\hookrightarrow ({\cal
S}_{W,X_1,X_2,{\cal G};a})_z$ is continuous. If $f\in ({\cal
S}_{W,X_1,X_2,{\cal G};a})_z$ and $g\in ({\cal S}_{W,X_1,X_2,{\cal
G};b})_z$, then there exist open $U_1\times V_1$ and $U_2\times V_2$
such that $f\in F_a(U_1\times V_1)$ and $g\in F_b(U_2\times V_2)$,
consequently, $f\in F_a(U\times V)$ and $g\in F_b(U\times V)$, where
$U=U_1\cup U_2$ and $V=V_1\cup V_2$, hence $\theta (f,g)\in
F_{a+b}(U\times V)$. From $(2)$ and the definition of the inductive
limit topology it follows, that $\theta $ is continuous, since on
iterated wrap groups $\theta $ is $H^{\infty }_p$ differentiable.

\par Moreover, in accordance with Theorem 21 \cite{lulaswgof}
$\theta $ is either associative or alternative if $\cal G$ is
associative or alternative.

\par {\bf 48. Note.} Let $\phi $ be a family of supports
in $X$ and $\cal B$ be a sheaf on $X$, where $\cal B$ may be
twisted. A sheaf $\cal B$ is called $\phi $-acyclic, if ${\sf
H}^b_{\phi }(X;{\cal B})=0$ for each $b>0$.

\par Let $\cal L^*$ be a resolvent of $\cal B$.
Put ${\cal Z}^b := Ker ({\cal L}^b\to {\cal L}^{b+1}) = Im ({\cal
L}^{b-1}\to {\cal L}^b)$, where ${\cal Z}^0 =\cal B$. An exact
sequence \par $(1)$ $e\to {\cal Z}^{b-1}\to {\cal L}^{b-1}\to {\cal
Z}^b\to e$ \\ induces an exact sequence \par $(2)$ $e\to \Gamma
_{\phi } ({\cal Z}^{b-1})\to \Gamma _{\phi }({\cal L}^{b-1}) \to
\Gamma
_{\phi }({\cal Z}^b)\to {\sf H}^1_{\phi }(X;{\cal Z}^{b-1})$. \\
Therefore, there exists the monomorphism \par $(3)$ ${\sf
H}^b(\Gamma _{\phi } ({\cal L}^*)) = \Gamma _{\phi }({\cal Z}^b)/Im
(\Gamma _{\phi }({\cal L}^{b-1}\to \Gamma _{\phi }({\cal Z}^b)) \to
{\sf H}^1_{\phi }(X;{\cal Z}^{b-1})$.
\par Moreover, the sequence $e\to {\cal Z}^{b-v}\to {\cal
L}^{b-v}\to {\cal Z}^{b-v+1}\to e$ induces the homomorphism:
\par $(4)$ ${\sf H}^{b-1}_{\phi }(X;{\cal Z}^{b-v+1})\to {\sf H}^v_{\phi
}(X;{\cal Z}^{b-v})$. \\ Define $\kappa $ as the composition
\par $(5)$ ${\sf H}^b(\Gamma _{\phi }({\cal L}^*))\to {\sf H}^1_{\phi
}(X;{\cal Z}^{b-1})\to {\sf H}^2_{\phi }(X;{\cal Z}^{b-2})\to ...
\to {\sf H}^b_{\phi }(X;{\cal Z}^0)$.
\par If all sheaves ${\cal L}^b$ are $\phi $-acyclic, then
$(3,4)$ are isomorphisms. We call $\kappa $ natural, if from the
commutativity of the diagram:
$$
\begin{array}{cccccc}
{\cal B}& \longrightarrow & {{\cal L}^*} \\
\downarrow\lefteqn{f}&&\downarrow\lefteqn{g} \\
{\cal E}& \longrightarrow & {{\cal M}^*}
\end{array}
$$
where $g$ is a homomorphism of resolvents the commutativity of the
diagram
$$
\begin{array}{ccccccccc}
{{\sf H}^b(\Gamma _{\phi }({\cal L}^*)}& \stackrel{\kappa
}{\longrightarrow}& {{\sf H}^b_{\phi }
(X;{\cal B})} \\
\downarrow\lefteqn{g^*}&&\downarrow\lefteqn{f^*}\\
{{\sf H}^b(\Gamma _{\phi }({\cal M}^*)}& \stackrel{\kappa
}{\longrightarrow}&{{\sf H}^b_{\phi }(X;{\cal E})}
\end{array}
$$
follows.  Thus we get the statement.

\par {\bf 48.1. Theorem.} {\it If ${\cal L}^*$ is the resolvent of
the sheaf $\cal B$, consisting of $\phi $-acyclic sheaves, then for
each $b\in \bf N$ the natural mapping
\par $\kappa : {\sf H}^b(\Gamma _{\phi }({\cal L}^*)\to {\sf
H}^b_{\phi }(X;{\cal B})$ is the isomorphism.}
\par In view of the latter theorem if
$g: {\cal L}^*\to {\cal M}^*$ is the homomorphism of two resolvents
of the sheaf $\cal B$ consisting of $\phi $-acyclic sheaves, then
the induced mapping ${\sf H}^b(\Gamma _{\phi }({\cal L}^*))\to {\sf
H}^b(\Gamma _{\phi }({\cal M}^*))$ is an isomorphism.

\par {\bf 48.2. Corollary.} {\it If $e\to {\cal L}^0\to {\cal
L}^1\to {\cal L}^2\to ...$ is an exact sequence of $\phi $-acyclic
sheaves, then the corresponding sequence $e\to \Gamma _{\phi }
({\cal L}^0)\to \Gamma _{\phi }({\cal L}^1)\to \Gamma _{\phi }({\cal
L}^2)\to ...$ is exact.}
\par {\bf Proof.} In view of Theorem 48.1
${\sf H}^b(\Gamma _{\phi }({\cal L}^*)) ={\sf H}^b_{\phi }(X;e)$. On
the other hand, ${\sf Y}^n_{\phi }(X;e)=e$, since ${\cal
Y}^0(X;e)=e$ and hence ${\cal Y}^n(X;e)=e$ for all $n$,
consequently, ${\sf H}^b(\Gamma _{\phi }({\cal L}^*))=e$ for each
$b$.

\par {\bf 49. Differential forms and twisted cohomologies
over octonions.} A bar resolution exists for any sheaf or a complex
of sheaves. Consider differential forms on $N$. In local coordinates
write a differential $k$-form as \par $(1)$ $w = \sum_J f_J(z)
dx_{b_1,j_1}\wedge dx_{b_2,j_2}\wedge ... \wedge dx_{b_k,j_k}$,
\\ where $f_J: N \to {\cal A}_r$, $z = (z_1,z_2,...)$ are local
coordinates in $N$, $z_b =x_{b,0}i_0 + x_{b,1} i_1 +...+
x_{b,{2^r-1}} i_{2^r-1}$, where $z_b\in {\cal A}_r$, $x_{b,j}\in \bf
R$ for each $b$ and every $j=0,1,...,2^r-1$, $J =(b_1,j_1;
b_2,j_2;...; b_k,j_k)$. For the sheaf ${\cal S}^k_{N,{\cal A}_r}$ of
germs of ${\cal A}_r$ valued $k$-forms on $N$ has a bar resolution:
\par $(2)$ $0\to {\cal S}^k_{N,{\cal A}_r}
\stackrel{\sigma }{\longrightarrow } {\cal S}^k_{N,A{\cal A}_r}
\stackrel{\sigma }{\longrightarrow } {\cal S}^k_{N,AB{\cal A}_r}
\stackrel{\sigma }{\longrightarrow } ...$,
\\ where ${\cal S}^k_{N,AB^m{\cal A}_r}$ denotes the sheaf of germs
of $AB^m{\cal A}_r$ valued $k$-forms on $N$.
\par Denote by ${\bf Z}(q,{\cal C}_r)$ the group analogous to
${\bf Z}({\cal C}_r)$ with $u\in {\cal C}_r$ replaced on $u^q$,
where $u^q$ is considered as equivalent with $(-u)^q$, $q\in \bf N$.
Therefore, the exponential sequence \par $(3)$ $0\to {\bf Z}({\cal
C}_r)_N\stackrel{\eta }\longrightarrow {\sf C}^{\infty }(N,{\cal
A}_r)\stackrel{\exp
}{\longrightarrow } {\sf C}^{\infty }(N,{\cal A}_r^*)\to 0$ \\
can be considered as a quasi-isomorphism:
$$
\begin{array}{ccccccccc}
{\bf Z}({\cal C}_r)_N & \stackrel{\eta }\longrightarrow & {\sf
C}^{\infty }(N,{\cal
A}_r) \\  \downarrow && \downarrow \lefteqn{\exp } \\
0 & \longrightarrow & {\sf C}^{\infty }(N,{\cal A}_r^*)
\end{array}
$$
between the complex ${\bf Z}({\cal C}_r)^{\infty }_D: {\bf Z}({\cal
C}_r)_N\to {\sf C}^{\infty }(N,{\cal A}_r)$ and the sheaf ${\sf
C}^{\infty }(N,{\cal A}_r^*)$ of germs of $C^{\infty }$ functions
from $N$ into ${\cal A}_r^*$ placed in degree one, that is ${\sf
C}^{\infty }(N,{\cal A}_r^*)[-1]$, where $\eta (z) =2\pi z$ for each
$z$ and $\exp (0)=1$ (see also \S 19), ${\cal A}_r$ is considered as
the additive group $({\cal A}_r,+)$, while ${\cal A}_r^*$ is the
multiplicative group $({\cal A}_r^*,\times )$. More generally this
gives the quasi-isomorphism:
\par $(4)$ ${{\bf Z}(1,{\cal C}_r)_N} {\longrightarrow} {\sf
C}^{\infty }(N,{\cal A}_r) {\stackrel{d}{\longrightarrow }} {{\cal
S}^1_{N,{\cal A}_r}} {\stackrel{d}{\longrightarrow }} {...}
{\stackrel{d}{\longrightarrow }} {{\cal S}^{q-1}_{N,{\cal A}_r}}$
and
\par ${0 \quad} {\longrightarrow\quad}  {{\sf C}^{\infty }(N,{\cal
A}_r^*) \quad} {\stackrel{dLn}{\longrightarrow }\quad} {{\cal
S}^1_{N,{\cal A}_r}} {\stackrel{d}{\longrightarrow }} {...}
{\stackrel{d}{\longrightarrow }} {{\cal S}^{q-1}_{N,{\cal A}_r}}$ \\
with vertical homomorphisms ${{\bf Z}(1,{\cal C}_r)_N}\to 0$, ${\sf
C}^{\infty }(N,{\cal A}_r)\stackrel{e}\longrightarrow {\sf
C}^{\infty }(N,{\cal A}_r^*)$, ${\cal S}^1_{N,{\cal A}_r}
\stackrel{id}\longrightarrow {\cal S}^1_{N,{\cal A}_r}$,...,${\cal
S}^{q-1}_{N,{\cal A}_r} \stackrel{id}\longrightarrow {\cal
S}^{q-1}_{N,{\cal A}_r}$ for $2\le q\in \bf N$, where $e(f):= \exp
(f)$ between a degree $q$ smooth twisted complex
\par $(5)$ ${\bf Z}({\cal C}_r)^{\infty }_D: {\bf Z}({\cal
C}_r)_N\to {\sf C}^{\infty }(N,{\cal
A}_r){\stackrel{d}{\longrightarrow }} {{\cal S}^1_{N,{\cal A}_r}}
{\stackrel{d}{\longrightarrow }} {...} {\stackrel{d}{\longrightarrow
}} {{\cal S}^{q-1}_{N,{\cal A}_r}}$ \\ and the complex ${\cal
S}^{<q}(N,{\cal A}_r)(dLn)[-1]$, where \par $(6)$ ${\cal
S}^{<q}(N,{\cal A}_r)(dLn): {\sf C}^{\infty }(N,{\cal
A}_r^*){\stackrel{dLn}{\longrightarrow }} {{\cal S}^1_{N,{\cal
A}_r}} {\stackrel{d}{\longrightarrow }} {...}
{\stackrel{d}{\longrightarrow }} {{\cal S}^{q-1}_{N,{\cal A}_r}}$.
\par The hypercohomology
$\mbox{ }_h{\sf H}^q(N,{\bf Z}({\cal C}_r)^{\infty }_D)$ of ${\bf
Z}({\cal C}_r)^{\infty }_D$ is a twisted non-commutative for $r=2$
and non-associative analog for $r=3$ of smooth Deligne cohomology of
$N$, since ${\cal A}_2 ={\bf H} = {\bf C}\oplus i_2{\bf C}$ and
${\cal A}_3 = {\bf O} = {\bf C}\oplus i_2{\bf C} \oplus i_4{\bf C}
\oplus i_6{\bf C}$ are quaternion and octonion algebras over $\bf R$
with the corresponding twisted structures causing twisted structures
of $AG$ and $BG$ as above. Thus hypercohomolgies have induced
twisted structures. We have that ${\cal S}^{<q}(N,{\cal A}_r)(dLn))$
is a trancation of the acyclic resolution $(6)$ of the constant
sheaf $({\cal A}_r^*)_N$. Therefore, the quasi-isomorphism $(5)$
implies
\par $(7)$ $\mbox{ }_h{\sf H}^b(N,{\bf Z}({\cal C}_r)^{\infty }_D)
\cong \mbox{ }_h{\sf H}^{b-1}({\cal S}^{<q}(N,{\cal A}_r)(dLn))$ for
each $b$ and $q$.
\par For the covering dimension $b=dim N$ (see \cite{eng})
there are the isomorphisms: \par $(8)$ $\mbox{ }_h{\sf H}^b({\cal
S}^{<b+1}(N,{\cal A}_r)(dLn)){\stackrel{e_N}{\longrightarrow }}
\mbox{ }_h{\sf H}^b(N,{\cal A}_r^*){\stackrel{t^b_N}{\longrightarrow
}} {\cal A}_r^*$, the composition of which is the isomorphism:
\par $(9)$ ${\sf T}^b_N: \mbox{ }_h{\sf H}^b({\cal
S}^{<b+1}(N,{\cal A}_r)(dLn)) {\longrightarrow } {\cal A}_r^*$.
\par There is useful the short exact sequence of complexes of
sheaves:
\par $(10)$ $0\to ({\cal A}_r^*)_N\to {\sf C}^{\infty } (N,{\cal A}_r^*)
\stackrel{dLn}{\longrightarrow } {\cal S}^1(N,{\cal A}_r)
\stackrel{d}{\longrightarrow } ... \stackrel{d}{\longrightarrow }
{\cal S}^q(N,{\cal A}_r)\stackrel{d}{\longrightarrow } {\cal
S}^{q+1,cl}(N,{\cal A}_r)\to 0$, \\
where ${\cal S}^{q+1,cl}(N,{\cal A}_r)$ denotes the sheaf of germs
of closed ${\cal A}_r$ valued $q+1$ forms on $N$.

\par {\bf 50. Remark.} Consider an open covering $ {\cal V} :=
\{ V_j: j\in J \} $ of a $H^{\infty }$ manifold $N$, denote by $
{\cal T}(E) := \{ g_j: g_j\in \Gamma (V_j,E), j\in J \} $ a family
of local trivializations of $E$, where $J$ is a set. If $V_k\cap
V_j\ne \emptyset $, then the quotient $g_{k,j} := g_k (1/g_j)$ is an
$H^{\infty }$ smooth ${\cal A}_r^*$-valued function on $V_k\cap
V_j$, where $1\le r\le 3$. If $1\le r\le 2$, then ${\cal A}_r$ is
associative and $g_{k,j}g_{j,l}=g_{k,l}$ on $V_k\cap V_j\cap V_l$,
when the latter set is not empty. \par For $r=3$ the octonion
algebra $\bf O$ is only alternative and generally $g_{k,j}g_{j,l}$
may be different from $g_{k,l}$. Already for quaternions and
moreover for octonions $Ln (xy)$ generally may be different from
$Ln(x)+Ln(y)$ for $x, y\in {\cal A}_r$ with $2\le r\le 3$ because of
non-commutativity. \par In view of Proposition 3.2
\cite{luoyst,luoyst2} for each $x, y\in {\cal A}_r$ there exists
$z\in {\cal A}_r$ such that
\par $(1)$ $e^xe^y = e^z = e^{a+b}e^{K(M,N)}$, where
$a=Re (x)$, $b=Re (y)$, $M= x-Re (x) =: Im (x)$, $N = Im (y)$, $K =
Im (z)$. As usually we denote by $ln $ the natural logarithmic
function in the commutative case $0\le r\le 1$, while $Ln$ denotes
the natural logarithmic function over ${\cal A}_r$ when $2\le r$
(see Section 3 in \cite{luoyst,luoyst2} and \cite{lufejms}). The
logarithmic function is defined on ${\cal A}_r\setminus \{ 0 \} $
for non-zero Cayley-Dickson numbers and has a non-commutative analog
of the Riemann surface so that $\exp $ and $Ln$ are ${\cal A}_r$
holomorphic. For each Cayley-Dickson number $v$ in the
multiplicative group ${\cal A}_r^* = {\cal A}_r\setminus \{ 0 \} $
there exists $x\in {\cal A}_r$ such that $e^x=v$. Then
\par $(2)$ $Ln (e^xe^y) = Ln (e^z) = a+b + K(M,N)$, where
\par $(3)$ $K(M,N)-M-N =: P(M,N)$ may be non-zero.
Express the real part as \par $(4)$ $Re (z) = (z + (2^r-2)^{-1} \{ -
z + \sum_{j=1}^{2^r-1} i_j (zi_j^*) \} )/2$, then
\par $(5)$ $Im (z) = z - Re (z) = (z -(2^r-2)^{-1} \{
- z + \sum_{j=1}^{2^r-1} i_j (zi_j^*) \} )/2$ \\
and fix these $z$-representations with which $M=M(x)$,  $N=N(y)$ and
$P(M,N)$ are locally analytic functions by $x$ and $y$. Put
\par $(6)$ $Ln (f_k) = w_k$ and
\par $(7)$ $Ln (g_{k,j}) = w_k - w_j + \nu _{k,j}$ and
\par $(8)$ $Ln (g_{k,l}) = Ln (g_{k,j}) + Ln (g_{j,l}) + \eta
_{k,j,l}$, \\
so that $w_k$ and $\nu _{k,j}$ and $\eta _{k,j,l}$ are $H^{\infty }$
differential $1$-forms. Then from $(6-8)$ it follows, that
\par $(9)$ $w_k - w_l + \nu _{k,l} = w_k - w_j + \nu _{k,j}
+ w_j - w_l + \nu _{j,l} + \eta _{k,j,k}$ and hence
\par $(10)$ $\eta _{k,j,l} = \nu _{k,l} - \nu _{k,j} - \nu _{j,l}$.
\\ Generally $\eta _{k,j,l}$ may be non-zero because of
non-commutativity or non-associativity.
\par In view of the alternativity of the octonion algebra $\bf O$
the identities $e^Me^N = e^K$, $e^M = e^Ke^{-N}$, $e^N = e^{-M}e^K$
and $e^{-K} = e^{-N}e^{-M}$ are equivalent, that leads to the
identities:
\par $(11)$ $M = K(K(M,N), - N)$, $N= K(-M, K(M,N))$, $K(M,N) = -
K(-N,-M)$, \\ where $M, N, K$ are purely imaginary octonions,
moreover, $K(M,0)=M$, $K(0,N)=N$, since $e^0=1$.
\par Let $E(N,{\cal A}_r^*,\pi ,\Psi )$ be an $H^{\infty }$
principal ${\cal A}_r^*$-bundle with transition functions $ \{
g_{k,j}: V_k\cap V_j\to {\cal A}_r^* : k, j \} $ and consider a
family $ \{ w_k, \nu _{k,j}, \eta _{k,j,l}: k, j, l \} $ of
$1$-forms related by Equations $(6-8)$ so that $w_j\in \Gamma (V_j,
{\cal S}^1_{N,{\cal A}_r}),$ $\nu _{k,j}\in \Gamma (V_k\cap V_j,
{\cal S}^1_{N,{\cal A}_r})$ for $V_k\cap V_j\ne \emptyset $, $\eta
_{k,j,l}\in \Gamma (V_k\cap V_j\cap V_l, {\cal S}^1_{N,{\cal A}_r})$
for $V_k\cap V_j\cap V_l\ne \emptyset $, where $k, j, l \in J$.
\par Consider a $C^{\infty }$ partition of unity $ \{ f_j: j\in J \} $
subordinated to the covering $\cal V$. Then \par $(12)$ $-w(x) =
|f_{j_0},f_{j_1},...,f_{j_n}, - w_{j_0}(x), -w_{j_1}(x),..., -
w_{j_n}(x)|$ and \par $(13)$ $ - \nu (x) = |f_{j_0}f_{k_0},
f_{j_1}f_{k_1},..., f_{j_n}f_{k_n}, - \nu _{j_0,k_0}(x), - \nu
_{j_1,k_1}(x),..., - \nu _{j_n,k_n}(x)|$ and \\ $(14)$ $ - \eta (x)
= |f_{j_0}f_{k_0}f_{l_0}, f_{j_1}f_{k_1}f_{l_1},...,
f_{j_n}f_{k_n}f_{l_n}, - \eta _{j_0,k_0,l_0}(x), - \nu
_{j_1,k_1,l_1}(x),..., - \nu _{j_n,k_n,l_n}(x)|$, \\ where $w_j(x)$
and $\nu _{j,k}(x)$ and $\eta _{j,k,l}(x)$ denote the restriction of
$w_j$ and $\nu _{j,k}$ and $\eta _{j,k,l}$ to $T_xN$ so that
$w_j(x)$ and $\nu _{j,k}(x)$ and $\eta _{k,j,l}(x)$ are $A{\cal
A}_r$-valued $1$-forms on $N$, \par $(15)$ $\pi _* (-w(x)) =
|f_{j_0},f_{j_1}(x),...,f_{j_n}(x); [w_{j_0}(x)-w_{j_1}(x) + \nu
_{j_0,j_1}(x)|...|w_{j_{n-1}}(x) - w_{j_n}(x) + \nu _{j_{n-1},j_n}]|$, \\
where $\pi : E{\cal A}_r\to B{\cal A}_r$ is the standard projection.
\par The principal $G$-bundle $E(N,G,\pi ,\Psi )$ is a pull-back
of the universal bundle $AG\to BG$ by a classifying mapping
$g_{E(N,G,\pi ,\Psi )}: N\to BG$. In terms of transition functions
\par $(16)$ $g_{E(N,{\cal A}_r^*,\pi ,\Psi )} =
|f_{j_0}(x),f_{j_1}(x),...,f_{j_n}(x);
[g_{j_0,j_1}(x)|g_{j_1,j_2}(x)|...|g_{j_{n-1},j_n}(x)]|$. Therefore,
\par $(17)$ $\pi _*(w) + d Ln (g_{E(N,{\cal A}_r^*,\pi ,\Psi )}) = 0$, \\
where for any differentiable function $g: U\to B{\cal A}_r^*$ we
have \par $g(x) = |f_0(x),f_1(x),...,f_n(x); [g_1(x)|...|g_n(x)]|$.
While
\par $(18)$ $d Ln (g(x)) := |f_0(x),f_1(x),...,f_n(x);
[d Ln (g_1(x))|... |d Ln (g_n(x))]|$.
\par Consider the total complex $(Tot^*(B_N^{*,<p}),D)$
of $B_N^{*,<p}$. Then a $(b-1)$-cocycle in the total complex is a
sequence $(g,w_1,...,w_{b-1})$, where $g\in H^{\infty
}(N,AB^{b-1}{\cal A}_r^*)$ and $w_j\in S^j_{AB^{b-1-j}{\cal A}_r}
(N)$ satisfying conditions:
\par $\sigma (g)=0$ which means that $g$ is a differentiable mapping
from $N$ into $B^{b-1}{\cal A}_r$; \par $\sigma (w_1) + dLn (g)=0$
means that $w_1$ is a connection on the differentiable principal
$B^{b-2}{\cal A}_r^*$-bundle over $N$ induced by $g$; \par $\sigma
(w_{j+1}) + (-1)^j dw_j =0$ serves as the definition of a
$(j+1)$-connection on a differentiable principal $B^{b-2}{\cal
A}_r^*$-bundle $E\to B$ associated with the mapping $g$ for $1\le
j\le b-2$. Then the sequence $(g,w_1,...,w_j)$ is called the
$j$-connection bar cocycle. \par There exists an equivalence
relation in the group of differentiable principal $B^{b-2}{\cal
A}_r^*$-bundles with $(b-1)$-connections which is induced by the
cohomology equivalence relation in the complex
$(Tot^*(B^{*,<b}_N),D)$. Thus ${\sf H}^{b-1}(Tot^*(B^{*,<b}),D)$ can
be identified with a group ${\sf E}(N,B^{b-2}{\cal A}_r^*,\nabla
^{b-1})$ of equivalence classes of differentiable principal
$B^{b-2}{\cal A}_r^*$-bundles with $(b-1)$-connections.
\par An assignment $(g,w_1,w_2,...,w_{b-1})\mapsto
(-1)^{b-1}dw_{b-1}$ induces a homomorphism $K: {\sf
E}(N,B^{b-2}{\cal A}_r^*,\nabla ^{b-1})\to S^b_{{\cal A}_r}(N)$
called the curvature of the $b$-connection
$(g,w_1,w_2,...,w_{b-1})$. The kernel $ker (K)$ is isomorphic to the
group ${\sf E}(N,B^{b-2}{\cal A}_r^*,\nabla ^{flat})$ of isomorphism
classes of differentiable principal $B^{b-2}{\cal A}_r^*$-bundles
with flat connections.

\par {\bf 51. Curvature of holonomy.} If $v, w\in T_0{\bf R}^n$, put
\par $(1)$ $\gamma _{v,w}(u) = 4uv$ for $0\le u\le 1/4$, $\gamma _{v,w}(u)=
v+4(u-1/4)w$ for $1/4 \le u\le 1/2$, $\gamma _{v,w}(u)= w-4(u-3/4)v$
for $1/2 \le u\le 3/4$, $\gamma _{v,w}(u) = 4(1-u)w$ for $3/4 \le u
\le 1$ and $\gamma ^s_{v,w}(u) := \gamma _{sv,sw}(u)$, where $0\le
u, s\le 1$. For a sequence of vectors ${\bf w} = (w_0,w_1,...,w_q)$
in $T_0{\bf R}^n$ with $q\in \bf N$ define a $(q+1)$-dimensional
parallelepiped $p[w_0,...,w_q]$ in the Euclidean space $\bf R^n$
with $q<n$ if $w_0,...,w_q$ are linearly independent. Then define
$\gamma _{w_0,w_1,w_2}(u_1,u_2) := \gamma _{\gamma
_{w_0,w_1}(u_1),w_2}(u_2)$ and by induction \par $(2)$ $\gamma
_{w_0,...,w_q}(u_1,...,u_q) = \gamma _{\gamma
_{w_0,...,w_{q-1}}(u_1,...,u_{q-1}), w_q}(u_q)$ and $\gamma ^s_{\bf
w}(u_1,...,u_q) := \gamma _{s\bf w}(u_1,...,u_q)$, where $0\le
u_1,...,u_q, s\le 1$. This gives the natural parametrization of the
parallelepiped $p[w_0,...,w_q]$ and the mapping $\gamma _{\bf w}:
\partial I^{q+1}\to {\bf R}^n$ which is continuous and piecewise
$C^{\infty }$. Denote by $e_j =(0,...,0,1,0,...,0)$ the standard
orthonormal basis in ${\bf R}^n$ with $1$ in the $j$-th place. Put
$Ln (diag (a_1,...,a_k)) := diag ( Ln (a_1),...,Ln (a_k))$, where
$Ln$ is the principal branch of the logarithmic function with $Ln
(1) =0$ and $diag (a_1,...,a_n)$ is the diagonal matrix with entries
$a_1,...,a_k\in {\cal A}_r^*$.
\par If $h$ is an $({\cal A}_r^*)^k$-valued $C^n$ holonomy or an
homomorphism for a wrap group $(W^ME)_{\infty ,H}$ with ${\hat M}$
being $H^{\infty }_p$ diffeomorphic with $\partial I^{m+1}$ and
$\psi = (y_1,...,y_n)$ is a coordinate system centered at $y$, $\psi
: V\to {\bf R}^n$, $V$ is an open neighborhood of a point $y$ in
$N$, then a curvature of $h$ at $y$ is a $q$-form
\par $(3)$ $K_y := \sum_{1\le j_1<...<j_q\le n}
K_{j_1,j_2,...,j_q}(y)dy_{j_1}\wedge dy_{j_2}\wedge ... \wedge
dy_{j_q}\in \Lambda ^qT_y^*N$, where
\par $(4)$ $K_{j_1,....,j_q}(y) = (-1)^q \lim_{s\to 0} Ln [h(\psi
^{-1}(\gamma ^s_{e_{j_1},...,e_{j_q}}))] s^{-q-1}$, \\
where $m\ge q$.

\par Consider the inversion $(w_j,w_{j+1})\mapsto (w_{j+1},w_j)$.
In view of Theorem 2 \cite{lulaswgof} for ${\hat M}$ being
$H^{\infty }_p$ diffeomorphic with $\partial I^{m+1}$ using the
iterated loops and the mapping $u_j\mapsto (1-u_j)$ we get, that
\par $(5)$ $K_y(w_{g(1)},...,w_{g(q+1)}) =
(-1)^{|g|} K_y(w_1,...,w_{q+1})$, \\
where $g\in S_{q+1}$, $S_q$ denotes the symmetric group of the set $
\{ 1,...,q \} $, $|g|=1$ for odd $g$, while $|g|=2$ for an even
transposition $g$.

\par {\bf 52. Remark.} Consider an $H^{\infty }$ manifold $N$ and a
pseudo-manifold $X$. A mapping $\gamma : X\to N$ is called piecewise
$C^{\infty }$ or $H^{\infty }$ smooth if it is continuous and the
restriction of $\gamma $ to each top dimensional simplex of $X$ is a
$C^{\infty }$ or $H^{\infty }$ mapping. A piecewise smooth mapping
$\gamma : X\to N$ is called an oriented singular pseudo-manifold
$q$-cycle, if $X$ is an oriented pseudo-manifold $q$-cycle. Denote
by ${\sf Z}^{\psi }_q(N):= {\sf Z}^{\psi }_q(X,N)$ the group of
oriented singular pseudo-manifold $q$-cycles in $N$.
\par If there exists an oriented pseudo-manifold with boundary
$(Y,\partial Y)$ with a pseudo-diffeomorphism $\eta : \partial Y\to
X$ and a piecewise smooth mapping $\zeta : Y\to N$ such that $\gamma
= \zeta |_{\partial Y}\circ \eta ^{-1}$, where $\gamma $ is an
oriented singular pseudo-manifold $q$-cycle, then $\gamma $ is
called an oriented singular pseudo-manifold $q$-boundary in $N$.
Denote by ${\sf B}^{\psi }_q(N) := {\sf B}^{\psi }_q(X,N)$ the group
of oriented singular pseudo-manifold $q$-boundaries in $N$.
\par Two oriented singular pseudo-manifold $q$-cycles $\gamma _j:
X_j\to N$, $j=1, 2$, are homologous, if there exists an oriented
$(q+1)$-dimensional pseudo-manifold with boundary $(Y,\partial Y)$
and a piecewise differentiable mapping $\zeta : Y\to N$ such that
$\partial Y $ is isomorphic with $X_1\cup X_2$ and $\zeta
|_{X_j}=\gamma _j$ up to an isomorphism $\partial Y \cong X_1\cup
X_2$ for $j=1, 2$.
\par Then there exists the group ${\sf H}^{\psi }_q(N) = {\sf Z}^{\psi
}_q(N)/{\sf B}^{\psi }_q(N)$ of homology classes of oriented
singular pseudo-manifold $q$-cycles in $N$, where the group
structure is given by the disjoint union.
\par Consider a twisted ${\cal A}_r$ analog of Cheeger-Simons
differential group functor consisting of pairs $(h,\alpha )\in Hom
({\sf Z}^{\psi }_q(N),{\cal A}_r^*)\times {\sf S}^{q+1}_{{\cal
A}_r}(N)_0$ satisfying the condition
\par $(CS)$ $h(\partial \eta ) =\exp ((-1)^q\int_{\eta }\alpha )$ for each
$\eta \in {\sf S}_{q+1}(N)$, \\ where ${\sf S}_q(N)$ is the group of
smooth singular $q$-chains in $N$, ${\sf S}^{q+1}_{{\cal A}_r}(N)_0$
denotes the group of closed differential ${\cal A}_r$-valued
$q$-forms on $N$ with $2\pi {\bf Z}({\sf C}_r)$-integral periods
belonging to ${\cal I}_r = \{ z\in {\cal A}_r: Re (z)=0 \} $, $1\le
r\le 3$. The Cheeger-Simons group ${\hat {\sf H}}^q_{\psi }(N,{\bf
Z}({\cal C}_r))$ of degree $q$ differential characters on $N$
consists of homomorphisms $h$ described above.

\par Suppose that $X$ is an $H^{\infty }_p$ pseudo-manifold.
Construct quotients ${\sf Z}^{\tilde {\psi }}_q(N)$ and
${\sf B}^{\tilde {\psi }}_q(N)$ as quotients of ${\sf Z}^{\psi
}_q(N)$ and ${\sf B}^{\psi }_q(N)$ by the equivalence relation:
\par $(E1)$ if $\gamma : X\to N$ is an oriented singular pseudo-manifold
$q$-cycle and $\xi $ is a homeomorphism of $X$ such that its
restrictions on all top dimensional simplices of a refinement of a
triangulation $\sf T$ of $X$ is an $H^{\infty }_p$ diffeomorphism,
then $\gamma \sim \gamma \circ \xi $ and as a class of equivalent
elements take $<\gamma
>_{\infty ,H}$ which is the closure relative to the $H^{\infty
}_p$-uniformity of the family of all such $\gamma \circ \xi $. In
view of the Morse and the Sard theorems (see \S \S II.2.10, 11
\cite{dubnovfom}) if $\delta \in <\gamma
>_{\infty ,p}$, then $\delta $ is homologous to $\gamma $.
Put ${\sf H}^{\tilde {\psi }}_q(N) := {\sf Z}^{\tilde {\psi }}_q(N)/
{\sf B}^{\tilde {\psi }}_q(N)$, then ${\sf H}^{\tilde {\psi }}_q(N)
\cong {\sf H}^{\psi }_q(N)$ are isomorphic.

\par {\bf 53. Higher twisted holonomies.} Suppose that $E(N,B{\cal
A}_r^*,\pi ,\Psi )$ is a differentiable principal $B{\cal
A}_r^*$-bundle with a classifying mapping $g: N\to B^q{\cal A}_r^*$
and a $q$-connection $(g,w_1,...,w_q)$, where $2\le r\le 3$.
Consider a $q$-dimensional orientable closed pseudo-manifold $X$
over ${\cal A}_r$ and $\gamma : X\to N$ an $H^{\infty }$ mapping. We
have that $B^q{\cal A}_r^*$ is $q$-connected and $g\circ \gamma :
X\to B^q{\cal A}_r^*$ is homotopic to a constant mapping. This
implies an existence of a differentiable mapping ${\overline {g\circ
\gamma }}: X\to AB^{q-1}{\cal A}_r^*$ with $\pi \circ {\overline
{g\circ \gamma }}=g\circ \gamma $, where $\pi : AB^{q-1}{\cal
A}_r^*\to B^q{\cal A}_r^*$. On the other hand, \par $\pi _*(\gamma
^*w_1 + d Ln{\overline {g\circ \gamma }}) = \pi _*\gamma ^*w_1 + d
Ln (g\circ \gamma ) = \gamma ^*(\pi _*w_1 + d Ln (g)) =0$, then
$(\gamma ^*w_1 + d Ln ({\overline {g\circ {\gamma }}})$ is a $B{\cal
A}_r$-valued $1$-form on $X$. \par The projection $\pi : A{\cal
A}_r\to B{\cal A}_r$ induces the surjective homomorphism $\pi _*:
{\sf S}^j_{A{\cal A}_r}(X)\to {\sf S}^j_{B{\cal A}_r}(X)$ for each
$j=1,2,...$. Therefore, there exists an $A{\cal A}_r$-valued
$1$-form ${\bar w}_j\in {\sf S}^j_{A{\cal A}_r}(X)$ satisfying the
equation:
\par $\pi _*{\bar w}_1 = \gamma ^*w_1 + d Ln
({\overline {g\circ {\gamma }}})$. Since $\sigma (\gamma ^*w_2 -
d{\bar w}_1) = \sigma \gamma ^*w_2 - d\gamma ^*w_1 = \gamma
^*(\sigma w_2 - dw_1) =0$, then $\gamma ^*w_2 - dw_1$ is an ${\cal
A}_r$-valued 2-form on $X$. \par By induction we get, that there
exists a differential $j$-form ${\bar w}_j\in {\sf S}^j_{A{\cal
A}_r}(X)$ such that $\pi _*{\bar w}_j = \gamma ^*w_j + (-1)^{j-1}
d\gamma ^*w_{j-1}$ for each $j=2,...,q$. We have that $\sigma
(\gamma ^*w_j + (-1)^{j-1} d{\bar w}_{j-1}) = \sigma \gamma ^*w_j +
(-1)^{j-1} d \gamma ^*w_{j-1} = \gamma ^*(\sigma w_j + (-1)^{j-1}
dw_{j-1})=0$, consequently, $\gamma ^*w_j + (-1)^{j-1} d{\bar
w}_{j-1}$ is an ${\cal A}_r$-valued $j$-form on $X$.
\par The holonomy of the $q$-connection $(g,w_1,...,w_q)$ along
$\gamma : X\to N$ is given by
\par $h(\gamma ) = \exp (\int_X (\gamma ^* w_q +(-1)^{q-1}d{\bar
w}_{q-1}))$. \\ If there is some other lift ${\hat w}_{j-1}$, then
${\hat w}_{j-1} = {\bar w}_{j-1} + v_{j-1}$, where $v_{j-1}$ is a is
an ${\cal A}_r$-valued $(j-1)$-form on $X$. Therefore,
\par $\int_X (\gamma ^*w_j + (-1)^{j-1}d{\hat w}_{j-1}) = \int_X
(\gamma ^*w_j + (-1)^{j-1}d{\bar w}_{j-1}) + (-1)^{j-1}
\int_Xdv_{j-1} = \int_X (\gamma ^*w_j + (-1)^{j-1}d{\bar w}_{j-1})$
in the considered here case of $X$ with $\partial X =0$.
\par This holonomy can be generalized in an abstract way
for an equivalence class $\eta $ of a $q$-connection
$(g,w_1,...,w_q)$ along a singular oriented pseudo-manifold $X$ of
dimension $q$ with an $H^{\infty }$ mapping $\gamma : X\to N$ such
that $h^{\eta }(\gamma )\in {\cal A}_r$. Define ${\sf
H}^{q+1}(X,{\bf Z}({\sf C}_r)(q+1)^{\infty }_D) := {\sf
H}^{q+1}(X\setminus S_X,{\bf Z}({\cal C}_r)(q+1)^{\infty}_D)$, where
$S_X$ is a singularity of $X$. If the dimension of $X$ is $q$, then
${\sf H}^{q+1}(X,{\bf Z}({\cal C}_r)(q+1)^{\infty }_D) \cong {\sf
H}^q(X,{\cal A}_r^*)$. Since $codim (S_X)\ge 2$, then ${\sf
H}^q(X,{\cal A}_r^*)$ has a fundamental class that induces an
integration along the fundamental class isomorphism and ${\sf
H}^q(X,{\cal A}_r^*) \cong {\cal A}_r^*$. Thus we get the
isomorphism $T^q_X: {\sf H}^{q+1}(X,{\bf Z}({\sf C}_r)(q+1)^{\infty
}_D)\to {\cal A}_r^*$. Therefore, $h^{\eta }(\gamma ) = T^q_X(\gamma
^*(\eta ))$ is the holonomy of a $q$-connection corresponding to an
element $\eta \in {\sf H}^{q+1}(N,{\bf Z}({\cal C}_r)(q+1)^{\infty
}_D)$ along $\gamma : X\to N$ for a singular oriented ${\cal A}_r$
pseudo-manifold $\phi : X\to N$ of a real dimension $q$, where $2\le
r\le 3$.

\par {\bf 54. Twisted cohomology.} Consider a twisted sheaf
$\cal B$ over $ \{ i_0,...,i_{2^r-1} \} $. Then a twisted analog of
an Alexander-Spanier (or of isomorphic \v{C}ech) cohomology with
coefficients in $\cal B$ and supports in the family $\phi $ is
$\mbox{ }_{AS}{\sf H}^*_{\phi } (X;{\cal B}) = {\sf H}^*(\Gamma
_{\phi }({\cal S}^*\otimes {\cal B}))$. \par Take an element $\eta
\in \mbox{ }_{AS} {\sf H}^q({\cal S}^{<q+1}_{N,{\cal A}_r}(dLn )) :=
\mbox{ }_{AS} {\sf H}^q(X,{\cal S}^{<q+1}_{N,{\cal A}_r}(dLn ))$,
where $N$ is of dimension $q$. Write it as $\eta = (g_{{\bf j}_q},
w_{{\bf j}_{q-1}},...,w_{{\bf j}_0} )$, where ${\bf j}_b :=
(j_0,...,j_b)$ is a multi-index. The ${\cal A}_r^*$-valued
$q$-cocycle $g_{{\bf j}_q}$ is cohomologous to zero, since $N$ is of
dimension $q$. Therefore, ${\sf H}^q({\sf C}^{\infty }_N({\cal
A}_r^*)) \cong {\sf H}^{q+1}(N,{\bf Z}({\sf C}_r)) \cong 0$. If
${\bar g}_{{\bf j}_{q-1}}$ is a $(q-1)$-cochain such that $\delta
({\bar g}_{{\bf j}_{q-1}}) = g^{-1}_{{\bf j}_q}$, then \par
$(g_{{\bf j}_q} \delta ({\bar g}_{{\bf j}_{q-1}}), dLn ({\bar
g}_{{\bf j}_{q-1}}) + w_{{\bf j}_{q-1}},...,w_{{\bf j}_0}) = (1, dLn
({\bar g}_{{\bf j}_{q-1}}) + w_{{\bf j}_{q-1}},...,w_{{\bf j}_0} )
=: \eta '$. Denote by $D$ the differential in the twisted complex of
${\cal S}^{<q+1}_{N,{\cal A}_r}(dLn )$, then the cocycle condition
$D(\eta ') =0$ leads to $\delta (d Ln ({\bar g}_{{\bf j}_{q-1}}) +
w_{{\bf j}_{q-1}}) =0$. Since ${\cal S}^1_{N,{\cal A}_r}$ is an
acyclic sheaf, then its twisted complex is exact and inevitably
there exists a $(q-2)$-cocycle ${\bar w}_{{\bf j}_{q-2}}$ for which
$\delta ({\bar w}_{{\bf j}_{q-2}}) = d Ln ({\bar g}_{{\bf j}_{q-1}})
+ w_{{\bf j}_{q-1}}$. Then $D(1, - {\bar w}_{{\bf j}_{q-2}},0,...,0)
+ \eta '$ is cohomolous to a cocycle having the form
$(1,0,{w'}_{{\bf j}_{q-1}},...,w_{{\bf j}_0})$. Continuing this
procedure gives a $(q-1)$-cochain ${\bar {\eta }} = ({\bar g}_{{\bf
j}_{q-1}}, {\bar w}_{{\bf j}_{q-2}},...,{\bar w}_{{\bf j}_0})$ of
the twisted complex ${\cal S}^{<q+1}_{N,{\cal A}_r}(dLn )$ such that
$\eta + D({\bar {\eta }}) =(1,0,...,0,{\hat w}_{{\bf j}_0})$. From
the cocycle condition $D(\eta + D({\bar {\eta }})) =0$ it follows
that ${\hat w}_{{\bf j}_0}$ is a global $q$-form on $N$ which we
will denote by $\hat w$. Put \par $(1)$ $T^q_N(\eta ) := \exp
((-1)^q\int_N{\hat w})$.
\par The mapping $T^q_N$ of Formula $(1)$ depends only on the
cohomology class of $\eta $, since if $\eta = D({\bar {\eta }})$,
then ${\hat w}=0$. Moreover, $T^q_N$ is independent from a choice of
the chain $\bar \eta $. Indeed, if $\tilde \eta $ is another
$(q-1)$-cochain with $\eta + D({\tilde {\eta }}) = (1,0,...,0, \nu
_{{\bf j}_0})$, then $\nu - {\hat w} $ is an $2\pi {\bf Z}({\sf
C}_r)$-integral $q$-form and inevitably $\exp (\int_N (\nu - {\hat
w})) =0$. \par Define the isomorphism $t^q_N: {\sf H}^q(N,{\cal
A}_r^*)\to {\cal A}_r^*$ as the restriction of $T^q_N$ to ${\sf
H}^q(N,{\cal A}_r^*)$ or it can be written as $t^q_N = T^q_N\circ
i^q_N$, where $i^q_N: \mbox{ }_{AS}{\sf H}^q(N,{\cal A}_r^*)\to
\mbox{ }_{AS} {\sf H}^m({\cal Y}^{<q+1}_{N,{\cal A}_r}(dLn ))$,
$i^q_N(g_{{\bf j}_q}) = (g_{{\bf j}_q},0,...,0)$ is the monomorphism
induced by the morphism of complexes of sheaves $({\cal A}_r^*)_N
\to {\cal S}^{<q+1}_{N,{\cal A}_r}(dLn )$.
\par Now construct an isomorphism $e^q_N: {\sf H}^q({\cal
S}^{<q+1}_{N,{\cal A}_r}(dLn )) \to {\sf H}^q(N,{\cal A}_r^*)$.
Consider the $q$-cocycle $\eta $ as above, then $dw_{{\bf j}_0} =0$,
since $N$ is of dimension $q$. This implies an existence of a
$(q-1)$-cochain $\bar \eta $ of the twisted complex of ${\cal
S}^{<q+1}_{N,{\cal A}_r}(dLn )$ such that $\eta + D({\bar \eta }) =
({\bar g}_{{\bf j}_q},0,...,0)$. The cocycle condition $D(\eta +
D({\bar \eta }))=0$ implies that ${\bar g}_{{\bf j}_q}$ is a locally
constant ${\cal A}_r^*$-valued cochain. Then the mapping $\eta
\mapsto {\bar g}_{{\bf j}_q}$ induces the isomorphism $e^q_N: \mbox{
}_{AS} {\sf H}^q({\cal S}^{<q+1}_{N,{\cal A}_r}(dLn )) \to \mbox{
}_{AS} {\sf H}^q(N,{\cal A}_r^*)$. This construction implies, that
$e^q_N$ is the inverse of $i^q_N$, hence $T^q_N = t^q_N\circ e^q_N$.

\par {\bf 55. Theorem.} {\it  For an $H^{\infty }$ manifold $N$ over
${\cal A}_r$ the mapping $\eta \mapsto h^{\eta }$ (see \S 54)
induces an isomorphism $h: {\sf H}^q({\cal S}^{<q+1}_{N,{\cal A}_r}
(d Ln))\to {\hat {\sf H}}^q(N,{\bf Z}({\sf C}_r))$.}
\par {\bf Proof.} It is sufficient to show, that the following
diagram with the upper row \par $0\to {\sf H}^q(N,{\cal A}_r^*)
\stackrel{i^p_N}{\longrightarrow} {\sf H}^q({\cal S}^{<q+1}_{N,{\cal
A}_r}(d Ln)) \stackrel{d}{\longrightarrow} {\sf S}^{q+1}_{{\cal
A}_r}(N)_0\to 0$ \\
and the lower row
\par $0\to Hom({\sf H}^{\psi }_q(N),{\cal A}_r^*)
\stackrel{{\hat i}^p_N}{\longrightarrow} {\hat {\sf H}}^q_{\psi
}(N,{\bf Z}({\sf C}_r)) \stackrel{\pi _2}{\longrightarrow} {\cal
S}^{q+1}_{{\cal A}_r}(N)_0\to 0$ \\
and with the vertical rows
\par ${\sf H}^q(N,{\cal A}_r^*)\stackrel{u}{\longrightarrow}
Hom({\sf H}^{\psi }_q(N),{\cal A}_r^*)$ and \par ${\sf H}^q({\cal
S}^{<q+1}_{N,{\cal A}_r}(d Ln))\stackrel{h}{\longrightarrow} {\hat
{\sf H}}^q_{\psi }(N,{\bf Z}({\sf C}_r))$ commutes. For each $\eta
\in {\sf H}^q ({\cal S}^{<q+1}_{N,{\cal A}_r}(d Ln))$ put $h(\eta )
:= (h^{\eta },K^{\eta })$, where $h^{\eta }$ is the holonomy of
$\eta $ and $K^{\eta }$ denotes the curvature of $\eta $. If $\eta =
(g_{j_q},w_{j_{q-1}},...,w_{j_0})$, then $K^{\eta } =
K(g_{j_q},w_{j_{q-1}},...,w_{j_0}) = dw_{j_0}$, consequently, the
right hand side of the above diagram commutes. \par The universal
coefficient theorem isomorphism $u: {\sf H}^q(N,{\cal A}_r^*)\to Hom
({\sf H}_q^{\psi } (N),{\cal A}_r^*)$ is induced by a pairing
assigning to a $H^{\infty }_p$ mapping $\gamma : M\to N$ and a
cohomology class $\eta \in {\sf H}^q(N,{\cal A}_r^*)$ an octonion or
quaternion number $t^q_M\circ \gamma ^*(\eta )$, where $t^q_M$
denotes the restriction of $T^q_M$ to ${\sf H}^q(N,{\cal A}_r^*)$.
Therefore, for $\eta \in {\sf H}^q(N,{\cal A}_r^*)$ and a $q$-cycle
$\sum_j n_j\gamma _j$, where $\gamma _j: M\to N$ we get $u(\eta
)(\sum_j n_j\gamma _j) = t^q_M(\prod_j \gamma _j^*(\eta )^{n_j})$.
\par From the equalities $h^{\eta }(\gamma ) = u(\eta )(\gamma )$
and $T^q_M=t^q_M\circ e^q_M$ and $e^q_M\circ i^q_M=id$ for an
arbitrary $H^{\infty }_p$ mapping $\gamma : M\to N$ it follows that
$h^{i^q_N(\eta )}(\gamma ) = T^q_M\gamma ^*i^q_N(\eta ) =
T^q_Mi^q_M\gamma ^*(\eta ) = t^q_Me^q_Mi^q_M\gamma ^*(\eta )=
t^q_M\gamma ^*(\eta ) = {\hat i}^q_Nu(\eta )(\gamma )$. Since $h$ is
the homomorphism, then the left hand side square of the diagram is
commutative as well.

\par {\bf 56. Remark.} In view of Theorem 55 every element of
${\hat {\sf H}}^q_{\psi }(N,{\bf Z}({\cal C}_r))$ is a holonomy
homomorphism. The operator $T^q_X$ in the definition of the holonomy
uses the integration which is invariant under the equivalence
relation 52$(E1)$. Then the quotient mapping ${\sf Z}^{\psi
}_q(N)\to {\sf Z}^{\tilde {\psi }}_q(N)$ induces an isomorphism
${\hat {\sf H}}^q_{\tilde {\psi }}(N,{\bf Z}({\sf C}_r))\to {\hat
{\sf H}}^q_{\psi }(N,{\bf Z}({\sf C}_r))$, where ${\hat {\sf
H}}^q_{\tilde {\psi }}(N,{\bf Z}({\sf C}_r))$ consists of pairs
$(h,v)\in Hom ({\sf Z}^{\tilde {\psi }}_q(N),{\cal A}_r^*)\times
{\cal S}^{q+1}_{{\cal A}_r}(N)_0$ so that $h(\partial \zeta )=\exp
((-1)^q\int_{\zeta }v)$ for each $\partial \zeta \in {\sf B}^{\tilde
{\psi }}_q(N)$.
\par A set theoretic inclusion $H^{\infty }_p(M,N)\to
{\sf Z}^{\psi }_q(M,N)$, where $q$ is a dimension of $M$, induces a
group homomorphism $\kappa: (W^MN)_{t,H}\to {\sf Z}^{\tilde \psi
}_q(M,N)$.
\par Denote by ${\cal L}^q_{N,{\cal A}_r^*}$ the sheaf associated
with the pre-sheaf \par $U\mapsto \{ \gamma \in Hom^{\infty
}((W^MN)_{\infty ,H},{\cal A}_r^*): supp (\gamma )\subset U \} $.
Section 53 and Theorem 55 imply that $K^h$ is an $2\pi {\bf Z}({\sf
C}_r)$-integral closed $(q+1)$-form on $N$.

\par {\bf 56.1. Lemma.} {\it For each $H^{\infty }_p$ mapping $\zeta :
Y\to N$, where $(Y,\partial Y)$ is a pseudo-manifold with boundary
$\partial Y$ being a pseudo-manifold over ${\cal A}_r$ and for each
extension ${\hat h}: {\sf Z}^{\tilde {\psi }}_{qb}(M,N)\to G^{kb}$
of an $H^{\infty }_p$ differentiable homomorphism $h:
(W^ME)_{b;\infty ,H}\to G^{kb}$ being an element of ${\hat {\sf
H}}^q_{\tilde {\psi }}(N,{\bf Z}({\cal C}_r))$ there is the
identity: \par ${\hat h}(\partial \zeta ) = \exp ((-1)^q\int_{\zeta
}K^h)$, where $b\in \bf N$, $E=E(N,G,\pi ,\Psi )$, $G$ is a
commutative subgroup in ${\cal A}_r^*$, $G$ is isomorphic with ${\bf
C}^*$.}

\par {\bf Proof.}  In view of the isomorphism ${\sf H}^q({\cal
S}^{<q+1}_{N,{\cal A}_r}(d Ln ))\to {\hat {\sf H}}^q_{\tilde {\psi
}}(N,{\bf Z}({\cal C}_r))$ each element of ${\hat {\sf H}}^q_{\tilde
{\psi }}(N,{\bf Z}({\cal C}_r))$ is a holonomy homomorphism. For
every pseudo-manifold with boundary $(Y,\partial Y)$ and an
$H^{\infty }_p$ mapping $\zeta : Y\to N$ take a partition $\cal T$
of $Y$ into small cubes $Q_j$. From the cancellation property of
holonomies we get that ${\hat h}(\partial \zeta ) = \prod_{Q_j\in
\cal T} {\hat h}(\gamma |_{Q_j})$, since $G$ is commutative. On the
other hand, ${\hat h}$ is an extension of $h$, hence $\prod_{Q_j\in
\cal T} {\hat h}(\gamma |_{Q_j}) = \prod_{Q_j\in \cal T} h(\gamma
|_{Q_j})$. Thus the proof reduces to the case of $Y$ being a $(q+1)$
dimensional cube in ${\cal A}_r^m\times \bf R$ such that $\partial
Y$ is embedded into ${\cal A}_r^m$ and has the real shadow $\partial
[0,1]^{q+1}$. \par If $\mu $ is a Borel measure on $Y$ relative to
which the Sobolev uniformity is given, then $\mu (Y_S)=0$, since
$codim (Y_S)\ge 2$, where $S_Y$ is the singularity of $Y$. Moreover,
a Lebesgue measure on ${\bf R}^{q+1}$ induces $\mu $ on $Y$ using
the fact that $Y\setminus Y_S$ is an $H^{t'}$-manifold with $t'
>[(q+1)/2] +1$.
\par For matrix-valued over ${\cal A}_r$ differential forms $w =
(w_{j,k}: j, k =1,..., m ) $ put $\int_{\xi }w =(\int_{\xi }w_{j,k}:
j, k =1,...,m )$, for diagonal matrices $(a_1,...,a_m)$ put $\exp
(a_1,...,a_m) := (e^{a_1},...,e^{a_m})$, if $a_1\ne 0$,...,$a_m\ne
0$, then $Ln (a_1,...,a_m) = (Ln (a_1),...,Ln (a_m))$.

\par Without loss of generality $h$ is additive and $\bf R$ homogeneous on
${\sf Z}_q(N)$. For each $n\in \bf N$ divide $[0,1]$ into $n$ small
subintervals, that induces a subdivision of $[0,1]^{q+1}$ into
$n^{q+1}$ cubes with vertices denoted by $v_{j_1,...,j_{q+1}}(n)$,
where $j_1,...,j_{q+1} =0,1,...,n$. Consider the wrap $\gamma
^n_{j_1,...,j_{q+1}} := \gamma _{e_1/n,...,e_{q+1}/n} +
v_{j_1,...,j_{q+1}}(n)$, where $e_1,...,e_{q+1}$ is the standard
basis of ${\bf R}^{q+1}$. \par Take $\xi \in H^{\infty }_p(Y,E)$
such that $\pi \circ \xi = \gamma $. Therefore, \par $\int_{\xi } K
= \int_Y \xi ^* K = \lim_{n\to \infty } \sum_{j_1,...,j_{q+1}} \xi
^* K(v_{j_1,...,j_{q+1}}(n))n^{-q-1} $ \\ $= (-1)^q \lim_{n\to
\infty } \sum_{j_1,...,j_{q+1}} \lim_{s\to 0} [Ln ~ h(\xi \circ
\gamma ^n_{j_1,...,j_{q+1}})] s^{-q-1} n^{-q-1}$, \\
where $\xi ^* K_y = \xi ^* K(y) dx_1\wedge ... \wedge dx_{q+1}$ for
each $y\in N$. Taking $s=1/n$ gives
\par $\lim_{n\to \infty } \lim_{s\to 0} \sum_{j_1,...,j_{q+1}}
[Ln ~ h(\xi \circ \gamma ^n_{j_1,...,j_{q+1}})] s^{-q-1} n^{-q-1}$
\par $ = \lim_{n\to \infty } \sum_{j_1,...,j_{q+1}} Ln ~ h(\xi \circ
\gamma ^n_{j_1,...,j_{q+1}})$ \\  $= \lim_{n\to \infty } Ln
(\prod_{j_1,...,j_{q+1}} h(\xi \circ \gamma ^n_{j_1,...,j_{q+1}}))
=\lim_{n\to \infty } Ln ~ h(\sum_{j_1,...,j_{q+1}} \xi \circ \gamma
^n_{j_1,...,j_{q+1}} )$ \\  $ = \lim_{n\to \infty } Ln h(\xi \circ
\gamma ^n_{0,...,0}) = Ln ~ h(\gamma )$, \\ since $h(\gamma
_1\lambda \lambda ^{-1} \gamma _2) = h(\gamma _1\gamma _2)$ and $G$
is commutative, where $\lambda : Y\to N$ is a path joining marked
points $y_1$ and $y_2$ of wraps $\gamma _1$ and $\gamma _2$, that is
$\gamma _j({\hat s}_{0,q})=y_j$ and $\lambda ({\hat s}_{0,q})=y_1$,
$\lambda ({\hat s}_{0,q+k}) = y_2$ while $\lambda ^{-1}({\hat
s}_{0,q})=y_2$ and $\lambda ^{-1}({\hat s}_{0,q+k})=y_1$ for each
$j=1, 2$ and $q=1,...,k$.

\par {\bf 57. Lemma.} {\it Suppose that $\phi : A\subset X$ is a
pointed inclusion of CW-complexes and $\theta : X\to X/A$ is the
quotient mapping. Let a group $G$ be twisted over $ \{
i_0,...,i_{2^r-1} \} $. Then $\theta _*: (W^ME;X,G,{\bf P})_{t,H}\to
(W^ME;X/A,G,{\bf P})_{t,H}$ is a principal $(W^ME;A,G,{\bf
P})_{t,H}$-bundle.}
\par {\bf Proof.} Let $G, E$ and $B$ be topological groups so that
$G$ acts effectively on $E$. Consider $U$ open in $B$ with $e\in U$.
Suppose that $\pi : E\to B$ is an open surjective mapping. Each
$G$-equivariant mapping $\xi : \pi ^{-1}\to G$ induces a local
trivialization of $\pi : E\to B$ over $U$. A group structure in $E$
induces a system of local trivializations of $E/B$. It is described
as follows. For each $v\in E$ take an open subset $U_v = \pi (v\pi
^{-1}(U))$ in $B$. Then the family $ \{ U_v: v\in E \} $ forms an
open covering of $B$. For each $v\in E$ there exists a
$G$-equivariant mapping $\xi _v: \pi ^{-1}(U_v) = v\pi ^{-1}(U)\to
G$ given by $\xi _v(x) = \xi (v^{-1}x)$.
\par Therefore, an open surjective mapping $\pi : E\to B$ is a principal
$G$-bundle if and only if there exists a neighborhood $U$ of the
unit element $e$ in $B$ and a $G$-equivariant mapping $\xi :
\pi^{-1}(U)\to G$.
\par Since the group $G$ is twisted, then due to Proposition 19 and
Theorem 20 \cite{lulaswgof} it is sufficient to prove this Lemma for
the commutative group $G_0$.

\par Consider a deformation retraction $\eta :
[0,1]\times V\to A$ of $V$ onto $A$, where $V$ is an open
neighborhood of $A$, put $U = \theta _*[(W^ME;V,G_0,{\bf
P})_{t,H}]$. A $(W^ME;A,G_0,{\bf P})_{t,H}$-equivariant mapping $\xi
: (\theta _*)^{-1} (U) = (W^ME;W,G_0,{\bf P})_{t,H} \to
(W^ME;A,G_0,{\bf P})_{t,H}$ is given by the formula $\xi (<{\bf
P}_{{\hat {\gamma }},v}>_{t,H}) = <\eta (1, {\bf P}_{{\hat {\gamma
}},v}>_{t,H}$ due to Propositions 7.1 and 13$(2)$ \cite{lulaswgof}.

\par {\bf 58. Theorem.} {\it For each connected smooth manifold
$N$, the homomorphism $\kappa $ induces an isomorphism
\par $\kappa _*: \pi _0((W^ME;N,{\cal A}_r^*,{\bf P})_{b;\infty ,H})
\to {\sf H}^{\tilde
{\psi }}_{qb}(N,{\bf Z}({\sf C}_r))$, where $1\le b\in \bf N$, $q$
is a dimension of $M$.}
\par {\bf Proof.} The uniform space
$H^{\infty }_p(M,E)$ is everywhere dense in the uniform space
$C^0(M,E)$ of all continuous mappings from $M$ into $E$, where $M$
is an $H^{\infty }_p$-pseudo-manifold. Therefore, there exists an
extension of $N\mapsto \pi _0((W^ME;N,{\cal A}_r^*,{\bf
P})_{b;\infty ,H})$ to a functor on the category of pointed
CW-complexes and pointed continuous mappings, that does not change a
homotopy type.
\par Recall a reduced homology theory. It is a functor ${\sf H}_*$
from the category of pointed CW-complexes and pointed continuous
mappings into the category of graded twisted groups satisfying the
properties $(H1-H4)$.
\par $(H1)$. For each pointed continuous mapping of CW-complexes
$f: X\to Y$ and $a\in \bf Z$, the induced homomorphism $f_*: {\sf
H}_a(X)\to {\sf H}_a(Y)$ depends only on the homotopy type of $f$.
\par $(H2)$. For each pointed CW-complex $X$ and $a\in \bf Z$ there
is a natural isomorphism
\par $\Sigma _X: {\sf H}_a(X)\to {\sf H}_{a+1}(\Sigma X)$, \\
where $\Sigma X$ is a reduced suspension of $X$.
\par $(H3)$. For each pointed inclusion $i: A\subset X$ of
CW-complexes and $a\in \bf Z$ the sequence \par ${\sf H}_a(A)
\stackrel{i_*}{\rightarrow} {\sf H}_a(X)\stackrel{g_*}{\rightarrow}
{\sf H}_a(X/A)$ is exact, \\ where $g: X\to X/A$ is the quotient
mapping.
\par $(H4)$. ${\sf H}_a(S^1)=e$ for $a\ne 1$ and ${\sf
H}_1(S^1)={\bf Z}$. These properties are standard and they are
demonstrated in Lemma 4.5 \cite{gajer} for commutative groups. Due
to Conditions 4$(A1,A2)$ on twisted groups we get the reduced
twisted homology theory.
\par In view of Lemma 57 $\pi _j((W^ME;A,G,{\bf P})_{t,H})
\to \pi _j((W^ME;X,G,{\bf P})_{t,H}) \to \pi _j((W^ME;X/A,G,{\bf
P})_{t,H})$ is a fragment of the long exact homotopy sequence of the
fibration $\theta _*$, where $G$ is the twisted group over $ \{
i_0,...,i_{2^r-1} \} $, $j=0, 1, 2,...$. Moreover, Conditions
$(H2,H3)$ follow from Lemma 57. Therefore, Properties $(H1-H4)$ for
twisted groups are direct consequences of the corresponding
properties for commutative groups. Though for the proof of this
theorem the case of commutative graded groups is sufficient.
\par Since $\kappa $ is a natural transformation of homology theory
and in view of Proposition 19 and Theorem 20 \cite{lulaswgof} this
induces the isomorphism $\kappa ^*$.

\par {\bf 59. Proposition.} {\it The curvature morphism
$K: {\cal L}^q_{N,{\cal A}_r^*}\to {\cal S}^{q+1,cl}_{N,{\cal A}_r}$
is an isomorphism.}
\par {\bf Proof.} The family ${\bf C}_M := {\bf R}\oplus M{\bf R}$
with $M\in {\cal A}_r$, $Re (M)=0$ and $|M|=1$ is such that its
union gives $\bigcup_M {\bf C}_M = {\cal A}_r$. In view of Theorem
55 and Lemma 56.1 $K$ is a monomorphism and an epimorphism from
${\cal L}^q_{N,{\cal A}_r^*}$ onto ${\cal S}^{q+1,cl}_{N,{\cal
A}_r}$. This gives the statement of this proposition.

\par {\bf 60. Theorem.} {\it  The restriction homomorphism
$\kappa ^*: Hom ({\sf Z}^{\tilde {\psi }}_{qb}(M,N),{\cal A}_r^*)\to
Hom ((W^MN)_{b;\infty ,H},{\cal A}_r^*)$ induces an isomorphism \par
${\hat {\kappa }}: {\hat {\sf H}}^{qb}_{\tilde {\psi }}(N,{\bf
Z}({\sf C}_r))\to Hom^{\infty }((W^MN)_{b;\infty ,H},{\cal A}_r)$,
where $1\le b\in \bf N$, $q$ is a covering dimension of $M$, $M$ and
$N$ are over ${\cal A}_r$.}
\par {\bf Proof.} Each homomorphism
$h: {\sf Z}^{\tilde {\psi }}_{qb}(M,N)\to {\cal A}_r^*$ is a
holonomy of ${\hat {\sf H}}^{qb}_{\tilde {\psi }}(N,{\bf Z}({\sf
C}_r))$, since the holonomy induces an isomorphism ${\sf
H}^{qb}({\cal S}^{<qb+1}_{N,{\cal A}_r}(d Ln ))\to {\hat {\sf
H}}^{qb}_{\tilde {\psi }}(N,{\bf Z}({\sf C}_r))$. Then the
restriction of $h$ on $(W^ME)_{b;\infty ,H}$ is of $H^{\infty }$
class, where $E=E(N,G,\pi ,\Psi )$ with $G={\bf C}_M$ (see \S 56.1).
Therefore, $\kappa ^*(h)\in Hom^{\infty }((W^MN)_{b;\infty ,H},{\cal
A}_r^*)$, since $(W^ME)_{b;\infty ,H}$ is the principal
$G^{bk}$-bundle over $(W^MN)_{b;\infty ,H}$.
\par Consider an extension ${\hat h}: {\sf Z}^{\tilde {\psi
}}_{qb}(M,N)\to {\cal A}_r^*$ of $h$, hence ${\hat h}\in {\hat {\sf
H}}^{qb}_{\tilde {\psi }}(N,{\bf Z}({\sf C}_r))$.
\par We have the locally analytic mapping $Ln$ from ${\cal A}_r^*$
onto ${\cal A}_r$. The group $(W^MN)_{b;\infty ,H}$ is commutative,
therefore instead of $Hom^{\infty }((W^MN)_{b;\infty ,H},{\cal
A}_r^*)$ we can consider the commutative additive group $Hom^{\infty
}((W^MN)_{b;\infty ,H},{\cal A}_r)$, where ${\cal A}_r$ is
considered as the additive group $({\cal A}_r,+)$. At the same time
the group ${\bf Z}({\cal C}_r)$ is commutative. Then $Ln (\kappa
^*(h))\in Hom^{\infty }((W^MN)_{b;\infty ,H},{\cal A}_r)$.
\par For each $\xi \in {\sf Z}^{\tilde {\psi }}_{qb}(M,N)$ there
exists $\zeta \in (W^MN)_{b;\infty ,H}$ and $\partial \eta \in {\sf
B}^{\tilde {\psi }}_{qb}(M,N)$ such that $\xi = \kappa (\zeta )
+\partial \eta $, since $\kappa _*: \pi _0((W^ME;N,{\cal A}_r^*,{\bf
P})_{b;\infty ,H})\to {\hat {\sf H}}^{qb}_{\tilde {\psi }}(N,{\bf
Z}({\sf C}_r))$ is an isomorphism due to Theorem 58. Then for each
extension ${\hat h}: {\sf Z}^{\tilde {\psi }}_{qb}(M,N)\to {\cal
A}_r^*$ of $h$ there is the identity:
\par ${\hat h}(\xi ) = h(\zeta ) + \exp ((-1)^{qb}\int_{\eta }K^h)$
\\ due to Section 53 and Lemma 59. Therefore, $h$ has a unique
extension $\hat h$. This implies that $\hat \kappa $ is an
isomorphism.

\par {\bf 61. Remark.} Mention that Theorems 55, 58 and 60 can be
proved in another way using the corresponding statements over $\bf
C$ and the twisted structure of sheaves over $ \{ i_0,...,i_{2^r-1}
\} $.


\begin{thebibliography}{99}

\bibitem{beggs} E.J. Beggs. "The de Rham complex of infinite
dimensional manifolds". Quart. J. Math. Oxford (2) {\bf 38} (1987),
131-154.

\bibitem{botttu} R. Bott, L.W Tu. "Differential forms in algebraic
topology" (New York: Springer-Verlag, 1982).

\bibitem{bredon} G.E. Bredon. "Sheaf theory" (New York: McGraw-Hill,
1967).

\bibitem{dingpan} Y.H. Ding, J.Z. Pang. "Computing degree
of maps between manifolds". Acta Mathem. Sinica. English Series.
{\bf 21: 6} (2005), 1277-1284.

\bibitem{dubnovfom} B.A. Dubrovin, S.P. Novikov, A.T. Fomenko.
"Modern geometry" (Moscow: Nauka, 1979).

\bibitem{ebi} D.G. Ebin, J. Marsden. "Groups of diffeomorphisms and
the motion of incompressible fluid". Ann. of Math. {\bf 92} (1970),
102-163.

\bibitem{eichh} J. Eichhorn. "The manifold structure of maps between
open manifolds". Ann. Glob. Anal. Geom. {\bf 3} (1993), 253-300.

\bibitem{eliass} H.I. Eliasson. "Geometry of manifolds of maps".
J. Differ. Geom. {\bf 1} (1967), 169-194.

\bibitem{emch} G. Emch. "M$\grave e$chanique quantique quaternionienne et
Relativit$\grave e$ restreinte". Helv. Phys. Acta {\bf 36} (1963),
739-788.

\bibitem{eng} R. Engelking. "General topology" (Moscow: Mir, 1986).

\bibitem{esnaualg} H. Esnault. "Algebraic theory of characteristic
classes of bundles with connection". In: Algebraic K-theory. Proc.
Symp. in Pure Math. {\bf 67} (1999), 13-25.

\bibitem{gajer} P. Gajer. "Higher holonomies,
geometric loop groups and smooth Deligne cohomology". In: "Advances
in geometry". Progr. in Math. {\bf  172}, 195-235 (Boston:
Birkh\"auser, 1999).

\bibitem{gajinvm} P. Gajer. "Geometry of Deligne cohomology".
Invent. Mathem. {\bf 127} (1997), 155-207.

\bibitem{guetze} F. G\"ursey, C.-H. Tze. "On the role of
division, Jordan and related algebras in particle physics"
(Singapore: World Scientific Publ. Co., 1996).

\bibitem{hamilt} W.R. Hamilton. "Selected papers. Optics. Dynamics.
Quaternions" (Moscow: Nauka, 1994).

\bibitem{harvey} F.R. Harvey. "Spinors and calibrations".
Perspectives in Mathem. {\bf 9} (Boston: Academic Press, 1990).

\bibitem{ish} C.J. Isham. "Topological and global aspects of
quantum theory". In: "Relativity, groups and topology.II" 1059-1290,
(Les Hauches, 1983). Editors: R. Stora, B.S. De Witt (Amsterdam:
Elsevier Sci. Publ., 1984).

\bibitem{kansol} I.L. Kantor, A.S. Solodovnikov.
"Hypercomplex numbers" (Berlin: Springer-Verlag, 1989).

\bibitem{kling} W. Klingenberg. "Riemannian geometry"
(Berlin: Walter de Gruyter, 1982).

\bibitem{lawmich} H.B. Lawson, M.-L. Michelson. "Spin geometry"
(Princeton: Princ. Univ. Press, 1989).

\bibitem{luwrgfbqo} S.V. Ludkovsky.
"Wrap groups of fiber bundles over quaternions and octonions". Los
Alam. Nat. Lab. {\bf math.FA 0802.0661}, 27 pages.

\bibitem{lulaswgof} S.V. Ludkovsky. "Structure of wrap groups
of quaternion and octonion fiber bundles". Los Alamos Nat. Lab. {\bf
math.GR 0804.4286}, 22 pages.
\bibitem{ludan} S.V. Ludkovsky. "Quasi-invariant measures on loop
groups of Riemann manifolds". Dokl. Akad. Nauk {\bf 370: 3} (2000),
306-308.

\bibitem{lupom} S.V. Ludkovsky. "Poisson measures for topological groups
and their representations".  Southeast Asian Bulletin of
Mathematics. {\bf 25} (2002), 653-680. (shortly in Russ. Math. Surv.
{\bf 56: 1} (2001), 169-170; previous versions: {\bf IHES/M/98/88},
38 pages, also Los Alamos Nat. Lab. {\bf math.RT/9910110}).

\bibitem{lufsqv} S.V. Ludkovsky.
"Functions of several Cayley-Dickson variables and manifolds over
them". J. Mathem. Sci. {\bf 141: 3} (2007), 1299-1330 (previous
variant: Los Alamos Nat. Lab. {\bf math.CV/0302011}).

\bibitem{norfamlud} S.V. Ludkovsky.
"Normal families of functions and groups of pseudoconformal
diffeomorphisms of quaternion and octonion variables". Sovrem.
Mathem. Fundam. Napravl. {\bf 18} (2006), 101-164 (previous variant:
Los Alam. Nat. Lab. math.DG/0603006).

\bibitem{luoyst} S.V. Ludkovsky, F. van Oystaeyen.
"Differentiable functions of quaternion variables". Bull. Sci. Math.
(Paris). Ser. 2. {\bf 127} (2003), 755-796.

\bibitem{luoyst2} S.V. Ludkovsky.
"Differentiable functions of Cayley-Dickson numbers and line
integration". J. Mathem. Sci. {\bf 141: 3} (2007), 1231-1298
(previous version: Los Alam. Nat. Lab. math.NT/0406048;
math.CV/0406306; math.CV/0405471).

\bibitem{lujmslg} S.V. Ludkovsky. "Stochastic
processes on geometric loop groups, diffeomorphism groups of
connected manifolds, associated unitary representations". J. Mathem.
Sci. {\bf 141: 3} (2007), 1331-1384 (previous version: Los Alam.
Nat. Lab. math.AG/0407439,  July 2004).

\bibitem{lufoclg} S.V. Ludkovsky.
"Geometric loop groups and diffeomorphism groups of manifolds,
stochastic processes on them, associated unitary representations".
In the book: "Focus on Groups Theory Research" (Nova Science
Publishers, Inc.: New York) 2006, pages 59-136.

\bibitem{lugmlg} S.V. Ludkovsky.
"Generalized geometric loop groups of complex manifolds, Gaussian
quasi-invariant measures on them and their representations". J.
Mathem. Sci. {\bf 122: 1} (2004), 2984-3011 (earlier version: Los
Alam. Nat. Lab. {\bf math.RT/9910086}, October 1999).

\bibitem{lufejms} S.V. Ludkovsky. "Quasi-conformal functions of quaternion and octonion
variables, their integral transformations". Far East J. of Math.
Sci. (FJMS) {\bf 28: 1} (2008), 37-88.

\bibitem{mensk} M.B. Mensky. "The paths group. Measurement. Fields.
Particles" (Moscow: Nauka, 1983).

\bibitem{miha} V.P. Mihailov. "Differential equations in
partial derivatives" (Moscow: Nauka, 1976).

\bibitem{michor} P.W.  Michor.  "Manifolds of Differentiable
Mappings" (Boston: Shiva, 1980).

\bibitem{milmorse} J. Milnor. "Morse theory" (Princeton, New Jersey:
Princeton Univ. Press, 1963).

\bibitem{milcw} J. Milnor. "On spaces having the homotopy type of
a CW-complex". Transactions of the A.M.S. {\bf 90} (1959), 272-280.

\bibitem{museya} N. Murakoshi, K. Sekigawa, A. Yamada.
"Integrability of almost quaternionic manifolds". Indian J. Mathem.
{\bf 42: 3}, 313-329 (2000).

\bibitem{nari} L. Narici, E. Beckenstein. "Topological vector spaces"
(New York: Marcel-Dekker Inc., 1985).

\bibitem{omo} H. Omori. "Groups of diffeomorphisms and their subgroups".
Trans. Amer. Math. Soc. {\bf 179} (1973), 85-122.

\bibitem{omori2} H. Omori. "Local structures of groups of
diffeomorphisms". J. Math. Soc. Japan {\bf 24: 1} (1972), 60-88.

\bibitem{pont} L.S. Pontrjagin. "Continuous groups"
(Moscow: Nauka, 1984).

\bibitem{seel} R.T. Seeley. "Extensions of $C^{\infty }$
Functions Defined in a Half Space". Proceed.  Amer.  Math.  Soc.
{\bf 15} (1964), 625-626.

\bibitem{souriau} J.M. Souriau. "Groupes differentiels"
(Berlin: Springer Verlag, 1981).

\bibitem{steen} N. Steenrod. "The topology of fibre budles"
(Princeton, New Jersey: Princeton Univ. Press, 1951).

\bibitem{sulwint} R. Sulanke, P. Wintgen. "Differentialgeometrie
und Faserb\"undel" (Berlin: Veb deutscher Verlag der Wissenschaften,
1972).

\bibitem{swan} R.C. Swan. "The Grothendieck Ring of a Finite Group".
Topology {\bf 2} (1963), 85-110.

\bibitem{swit} R.M.  Switzer.  "Algebraic Topology - Homotopy and Homology"
(Berlin: Springer-Verlag, 1975).

\bibitem{touger} J.C.  Tougeron.  "Ideaux de
Fonctions Differentiables" (Berlin: Springer-Verlag, 1972).

\bibitem{whiteh} J.H.C. Whitehead. "Combinatorial homotopy.I".
Bull. Amer. Mathem. Soc. {\bf 55} (1949), 213-245.

\bibitem{yano} K. Yano, M. Ako. "An affine connection in almost
quaternion manifolds". J. Differ. Geom. {\bf 3} (1973), 341-347.

\bibitem{zorich} V.A. Zorich. "Mathematical analysis", V. 2
(Moscow: Nauka, 1984).

\end{thebibliography}
\end{document}